# INTERVAL SEMIGROUPS

W. B. Vasantha Kandasamy
Florentin Smarandache

**2011**

# INTERVAL SEMIGROUPS

W. B. Vasantha Kandasamy
Florentin Smarandache

**2011**

# CONTENTS

# PREFACE

In this book we introduce the notion of interval semigroups using intervals of the form [0, a], a is real. Several types of interval semigroups like fuzzy interval semigroups, interval symmetric semigroups, special symmetric interval semigroups, interval matrix semigroups and interval polynomial semigroups are defined and discussed. This book has eight chapters.

The main feature of this book is that we suggest 241 problems in the eighth chapter. In this book the authors have defined 29 new concepts and illustrates them with 231 examples. Certainly this will find several applications.

The authors deeply acknowledge Dr. Kandasamy for the proof reading and Meena and Kama for the formatting and designing of the book.

<div align="right">
W.B.VASANTHA KANDASAMY<br>
FLORENTIN SMARANDACHE
</div>

# ~ DEDICATED TO ~

# Ayyankali

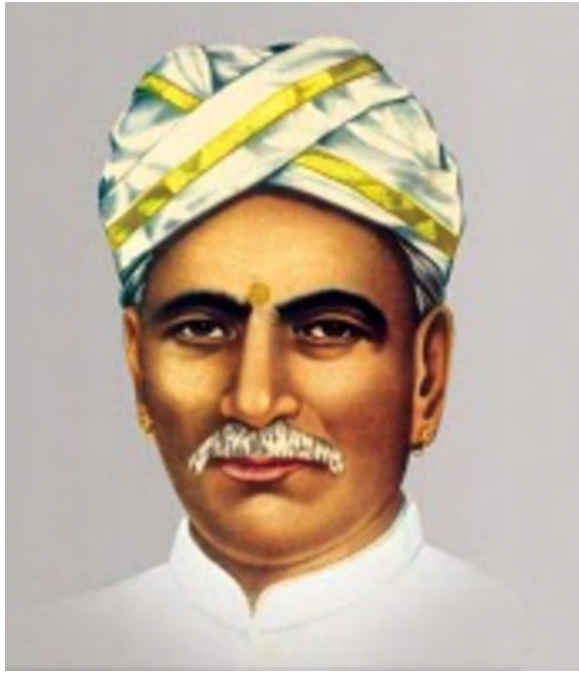

*Ayyankali* (1863–1941) was the first leader of Dalits from Kerala. He initiated several reforms to emancipate the lives of the Dalits. Ayyankali organized Dalits and fought against the discriminations done to Dalits and through his efforts he got the right to education, right to walk on the public roads and dalit women were allowed to cover their nakedness in public. He spearheaded movements against casteism.

**Chapter One**

# INTRODUCTION

We in this book make use of special type of intervals to build interval semigroups, interval row matrix semigroups, interval column matrix semigroups and interval matrix semigroups. We also introduce and study the Smarandache analogue of them.

The new notion of interval symmetric semigroups and special interval symmetric semigroups are defined and studied. For more about symmetric semigroups and their Smarandache analogue concepts please refer [9].

The classical theorems for finite groups like Lagrange theorem, Cauchy theorem and Sylow theorem are introduced in a special way and analyzed. Only under special conditions we see the notion of these classical theorems for finite groups can be extended interval semigroups.

The authors also introduce the notion of neutrosophic intervals and fuzzy intervals and study them in the context of interval semigroups. I denotes the indeterminate or inderminancy where $I^2 = I$ and $I + I = 2I$, $I + I + I = 3I$ and so on. For more about neutrosophy, neutrosophic intervals please refer [1, 3, 6-8].

Study of special elements like interval zerodivisors, interval idempotents, interval units, interval nilpotents are studied and their Smarandache analogue introduced [9].

Chapter Two

# INTERVAL SEMIGROUPS

In this chapter we for the first time introduce the notion of interval semigroups and describe a few of their properties associated with them. We see in general several of the classical theorems are not true in general case of semigroups. First we proceed on to give some notations essential to develop these new structures.

$$I(Z_n) = \{[0, a_m] \mid a_m \in Z_n\},$$
$$I(Z^+ \cup \{0\}) = \{[0, a] \mid a \in Z^+ \cup \{0\}\},$$
$$I(Q^+ \cup \{0\}) = \{[0, a] \mid a \in Q^+ \cup \{0\}\},$$
$$I(R^+ \cup \{0\}) = \{[0, a] \mid a \in R^+ \cup \{0\}\}$$
and $\quad I(C^+ \cup \{0\}) = \{[0, a] \mid a \in C^+ \cup \{0\}\}.$

**DEFINITION 2.1:** *Let $S = \{[0, a_i] \mid a_i \in Z_n; +\}$ S is a semigroup under addition modulo n. S is defined as the interval semigroup under addition modulo n.*

We will first illustrate this by some simple examples.

***Example 2.1:*** Let $S = \{[0, a_i] \mid a_i \in Z_6\}$, under additions is an interval semigroup. We see S is of finite order and order of S is six.

***Example 2.2:*** Let $S = \{[0, a_i] \mid a_i \in Z_{12}\}$ be an interval semigroup under addition modulo 12. This is also a interval semigroup of finite order.

Now we can define interval semigroup under addition using $Z^+ \cup \{0\}$, $Q^+ \cup \{0\}$, $R^+ \cup \{0\}$ and $C^+ \cup \{0\}$. All these interval semigroups are of infinite order.

We will illustrate these situations by some examples.

***Example 2.3:*** Let $S = \{[0, a_i] \mid a_i \in Z^+ \cup \{0\}\}$; S is an interval semigroup under addition. Clearly S is of infinite order.

***Example 2.4:*** Let $S = \{[0, a_i] \mid a_i \in Q^+ \cup \{0\}\}$; S is an interval semigroup under addition. Clearly S is of infinite order.

***Example 2.5:*** Let $S = \{[0, a_i] \mid a_i \in R^+ \cup \{0\}\}$; S is an interval semigroup under addition and is of infinite order.

***Example 2.6:*** Let $S = \{[0, a_i] \mid a_i \in C^+ \cup \{0\}\}$; S is an interval semigroup under addition and is of infinite order.

Thus we have seen examples of interval semigroups under addition, these are known as basic interval semigroups under addition.

We will now define polynomial interval semigroups and matrix interval semigroups defined using basic interval semigroups and then define their substructures.

**DEFINITION 2.2:** *Let $S = \{([0, a_1], [0, a_2], ..., [0, a_n]) \mid a_i \in Z_n\}$; S under component wise addition is an interval semigroup known as the row matrix interval semigroup.*

We can in the definition replace $Z_n$ by $Z^+ \cup \{0\}$ or $R^+ \cup \{0\}$ or $Q^+ \cup \{0\}$ or $C^+ \cup \{0\}$.

We will illustrate these by some examples.

**Example 2.7:** Let $S = \{([0, a_1], [0, a_2], [0, a_3], [0, a_4], [0, a_5]) \mid a_i \in Z_{12}; 1 \le i \le 5\}$ is a row matrix interval semigroup under addition.

**Example 2.8:** Let $P = \{([0, a_1], [0, a_2], [0, a_3], [0, a_4], [0, a_5], [0, a_6]) \mid a_i \in Z^+ \cup \{0\}; 1 \le i \le 6\}$ is a row matrix interval semigroup under addition. Clearly P is of infinite order.

**Example 2.9:** Let $S = \{([0, a_1], [0, a_2], [0, a_3], \ldots, [0, a_{12}]) / a_i \in Q^+ \cup \{0\}; 1 \le i \le 12\}$; S is a row matrix interval semigroup under addition and is of infinite order.

**Example 2.10:** Let $S = \{([0, a_1], [0, a_2], [0, a_3], \ldots, [0, a_{15}]) / a_i \in R^+ \cup \{0\}; 1 \le i \le 15\}$ be a row matrix interval semigroup under addition; T is of infinite order.

**Example 2.11:** Let $G = \{([0, a_1], [0, a_2]) \mid a_i \in C^+ \cup \{0\}; 1 \le i \le 2\}$; be a row matrix interval semigroup of infinite order.

Now we proceed onto define column matrix interval semigroup.

**DEFINITION 2.3:** *Let*

$$S = \left\{ \begin{bmatrix} [0, a_1] \\ [0, a_2] \\ [0, a_3] \\ \vdots \\ [0, a_n] \end{bmatrix} \middle| a_i \in Z_m; 1 \le i \le n \right\},$$

*S under addition modulo m is a semigroup defined as the column interval matrix semigroup under addition.*

We can replace $Z_n$ is definition 2.3 by $Z^+ \cup \{0\}$ or $Q^+ \cup \{0\}$, $R^+ \cup \{0\}$ or $C^+ \cup \{0\}$ and get column interval matrix semigroups under addition.

We will illustrate these situations by some examples.

*Example 2.12:* Let

$$S = \left\{ \begin{bmatrix} [0,a_1] \\ [0,a_2] \\ [0,a_3] \\ [0,a_4] \\ [0,a_5] \end{bmatrix} \middle| a_i \in Z_5; 1 \le i \le 5 \right\}$$

be a column interval matrix semigroup under addition. S is of finite order.

*Example 2.13:* Let

$$S = \left\{ \begin{bmatrix} [0,a_1] \\ [0,a_2] \\ [0,a_3] \\ \vdots \\ [0,a_{15}] \end{bmatrix} \middle| a_i \in Z^+ \cup \{0\}; 1 \le i \le 15 \right\}$$

be a column interval matrix semigroup under addition.

*Example 2.14:* Let

$$S = \left\{ \begin{bmatrix} [0,a_1] \\ [0,a_2] \\ [0,a_3] \\ \vdots \\ [0,a_{10}] \end{bmatrix} \middle| a_i \in Q^+ \cup \{0\}; 1 \le i \le 10 \right\},$$

be a column interval semigroup under addition of infinite order.

*Example 2.15:* Let

$$S = \left\{ \begin{bmatrix} [0, a_1] \\ [0, a_2] \\ [0, a_3] \\ \vdots \\ [0, a_6] \end{bmatrix} \middle| a_i \in R^+ \cup \{0\}; 1 \leq i \leq 6 \right\}$$

be a column interval semigroup under addition of infinite order.

*Example 2.16:* Let

$$P = \left\{ \begin{bmatrix} [0, a_1] \\ [0, a_2] \\ [0, a_3] \\ [0, a_4] \end{bmatrix} \middle| a_i \in C^+ \cup \{0\}; 1 \leq i \leq 4 \right\}$$

be a column interval semigroup under addition of infinite order.

Now we will define matrix interval semigroup.

**DEFINITION 2.4:** *Let $S = \{m \times n$ interval matrices with entries from $I(Z_n)\}$ be a $m \times n$ matrix interval semigroup under addition.*

We can replace $I(Z_n)$ in definition 2.4 by $I(Z^+ \cup \{0\})$ or $I(R^+ \cup \{0\})$ or $I(Q^+ \cup \{0\})$ or $I(C^+ \cup \{0\})$ and get $m \times n$ interval matrix semigroups under addition.

We will illustrate all these by some examples.

*Example 2.17:* Let

$$S = \left\{ \begin{bmatrix} [0, a_1] & [0, a_2] \\ [0, a_3] & [0, a_4] \\ [0, a_5] & [0, a_6] \end{bmatrix} \middle| a_i \in Z_{10}; 1 \leq i \leq 6 \right\}$$

be a 3 × 2 interval matrix semigroup under matrix addition modulo 10 of finite order.

*Example 2.18:* Let

$$S = \left\{ \begin{bmatrix} [0,a_1] & [0,a_2] & [0,a_3] \\ [0,a_4] & [0,a_5] & [0,a_6] \\ [0,a_7] & [0,a_8] & [0,a_9] \end{bmatrix} \middle| a_i \in Z_{42}; 1 \le i \le 9 \right\}$$

be a 3 × 3 interval square matrix semigroup of finite order under interval matrix addition modulo 42.

*Example 2.19:* Let

$$S = \left\{ \begin{bmatrix} [0,a_1] & [0,a_2] & [0,a_3] \\ [0,a_4] & [0,a_5] & [0,a_6] \\ \vdots & \vdots & \vdots \\ [0,a_{22}] & [0,a_{23}] & [0,a_{24}] \end{bmatrix} \middle| a_i \in Q^+ \cup \{0\}; 1 \le i \le 24 \right\}$$

be a 8 × 3 matrix interval semigroup under addition of infinite order.

*Example 2.20:* Let

$$P = \left\{ \begin{bmatrix} [0,a_1] & \cdots & [0,a_4] \\ [0,a_5] & \cdots & [0,a_8] \\ [0,a_9] & \cdots & [0,a_{12}] \\ [0,a_{13}] & \cdots & [0,a_{16}] \end{bmatrix} \middle| a_i \in R^+ \cup \{0\}; 1 \le i \le 16 \right\}$$

be a 4 × 4 square matrix interval semigroup of infinite order.

**DEFINITION 2.5:** *Let S be the matrix interval semigroup under addition. Let $M \subseteq S$ (M a proper subset of S), if M itself is a matrix interval semigroup under addition then we define M to be a matrix interval subsemigroup of S.*

We will illustrate this situation by some examples.

*Example 2.21:* Let

$$S = \left\{ \begin{bmatrix} [0,a_1] & [0,a_2] \\ [0,a_3] & [0,a_4] \end{bmatrix} \;\middle|\; a_i \in Z_5; 1 \le i \le 4 \right\}$$

be a square matrix interval semigroup under addition.

$$P = \left\{ \begin{bmatrix} [0,a] & [0,a] \\ [0,a] & [0,a] \end{bmatrix} \;\middle|\; a \in Z \right\} \subseteq S;$$

P is a square matrix interval subsemigroup of S.

*Example 2.22:* Let

$$S = \left\{ \begin{bmatrix} [0,a_1] & [0,a_2] \\ [0,a_3] & [0,a_4] \\ [0,a_5] & [0,a_6] \end{bmatrix} \;\middle|\; a_i \in Q^+ \cup \{0\}; 1 \le i \le 6 \right\}$$

be a $3 \times 2$ matrix interval semigroup under addition.

$$P = \left\{ \begin{bmatrix} [0,a] & [0,a] \\ [0,a] & [0,a] \\ [0,a] & [0,a] \end{bmatrix} \;\middle|\; a \in Z^+ \cup \{0\} \right\} \subseteq S$$

is a $3 \times 2$ matrix interval subsemigroup of S under addition.

*Example 2.23:* Let

$$M = \left\{ \begin{bmatrix} [0, a_1] \\ [0, a_2] \\ \vdots \\ [0, a_{12}] \end{bmatrix} \;\middle|\; a_i \in R^+ \cup \{0\} \right\}$$

be a $12 \times 1$ column matrix interval semigroup under addition. Choose

$$S = \left\{ \begin{bmatrix} [0, a] \\ [0, a] \\ \vdots \\ [0, a] \end{bmatrix} \;\middle|\; a \in R^+ \cup \{0\} \right\} \subseteq M;$$

S is a $12 \times 1$ column matrix interval subsemigroup of M.

*Example 2.24:* Let $M = \{([0, a_1], [0, a_2], \ldots, [0, a_{19}]) \mid a_i \in C^+ \cup \{0\}; 1 \leq i \leq 19\}$ be a $1 \times 19$ row matrix interval semigroup. Choose $W = \{([0, a], [0, a], \ldots, [0, a]) \mid a \in R^+ \cup \{0\}\} \subseteq M$; W is a $1 \times 19$ row matrix interval subsemigroup of M.

We can define ideals as in case of usual semigroups.

**DEFINITION 2.6:** *Let S be a matrix interval semigroup under addition. I a proper subset of S. I is said to be a matrix interval ideal of the semigroup S if (a) I is a matrix interval subsemigroup of S. (b) for each $s \in S$ and $a \in I$. $a + s$ and $s + a$ are in I.*

The notion of left and right ideal in a matrix interval semigroup can also be defined as a matter of routine.

The following theorem is obvious and the reader is expected to prove.

**THEOREM 2.1:** *Let S be a matrix interval semigroup. Every ideal I of S is a matrix interval subsemigroup of S but a matrix interval subsemigroup in general is not a matrix interval ideal of S.*

We will illustrate this by some simple examples.

*Example 2.25:* Let S be a $1 \times 5$ row matrix interval semigroup under addition with entries from $3Z^+ \cup \{0\}$. Choose I = $\{([0, a_1], [0, a_2], [0, a_3], [0, a_4], [0, a_5]) \mid a_i \in 6Z^+ \cup \{0\}\} \subseteq$ S is not a $1 \times 5$ matrix interval ideal of S.

We see some matrix interval subsemigroups of S are not in general matrix interval ideals. For take J = $\{([0, a], [0, a], \ldots, [0, a]) \mid a \in 9Z^+ \cup \{0\}\} \subseteq$ S. We see J is not a $1 \times 5$ matrix interval ideal, J is only a $1 \times 5$ matrix interval subsemigroup of S under addition.

*Example 2.26:* Let S = $\{([0, a_1], [0, a_2], \ldots, [0, a_{10}]) \mid a_i \in Z_{12}\}$ be a $1 \times 10$ row interval matrix semigroup. Consider any $1 \times 10$ row interval matrix subsemigroup, we see it cannot be an ideal.

For example take V = $\{([0, a_1], [0, a_2], \ldots, [0, a_{10}]) \mid a_i \in \{0, 6\} \subseteq Z_{12}\} \subseteq$ S; V is a $1 \times 10$ row interval matrix subsemigroup of S but is not an ideal. For take $2 \in Z_{12}$ we see V + 2 $\not\subseteq$ V for V + 2 = $\{([0, a_1], [0, a_2], \ldots, [0, a_{10}]) \mid a_i \in \{0, 6, 8, 2\}\} \neq$ V. Thus we see V has no $1 \times 10$ row interval matrix ideal.

We see it is difficult to get ideals in case of row matrix interval semigroups under addition, but however we have ideals in case of row matrix interval semigroups under multiplication.

We have a class of interval matrix semigroups under addition which have no ideals. We will call those interval matrix semigroups which have no subsemigroups as simple and those will have no ideals. Those interval matrix semigroups which has subsemigroups but no ideals as doubly simple.

*Example 2.27:* Let S = $\{([0, a_1], [0, a_2], \ldots, [0, a_9]) \mid a_i \in Z_5; 1 \leq i \leq 9\}$ be a interval matrix semigroup under addition, S has no interval matrix subsemigroups as well as S has no interval matrix ideals. Infact S is a doubly simple interval matrix as S has no interval matrix ideals and subsemigroups.

Now we proceed onto define interval matrix semigroup under multiplication.

**DEFINITION 2.7:** *Let S be a interval matrix semigroup under multiplication using $I(Z_n)$ or $I(Z^+ \cup \{0\})$ or $I(Q^+ \cup \{0\})$ or $I(R^+ \cup \{0\})$.*

We will illustrate this situation by examples.

*Example 2.28:* Let $V = \{([0, a_1], [0, a_2], \ldots, [0, a_9]) \mid a_i \in Z_{12}; 1 \leq i \leq 9\}$; V is an interval matrix semigroup under multiplication. Clearly V is a finite interval matrix semigroup.

*Example 2.29:* Let $V = \{([0, a_1], [0, a_2], \ldots, [0, a_{10}]) \mid a_i \in Z^+ \cup \{0\}; 1 \leq i \leq 10\}$ be an interval row matrix semigroup under multiplication. Clearly V is of infinite order.

*Example 2.30:* Let $V = \{([0, a_1], [0, a_2], \ldots, [0, a_{12}]) \mid a_i \in Q^+ \cup \{0\}; 1 \leq i \leq 12\}$ be an interval row matrix semigroup under multiplication.

We can define subsemigroups and ideals in case of these semigroups.

**DEFINITION 2.8:** *Let $V = \{([0, a_1], [0, a_2], \ldots, [0, a_n]) / a_i \in Q^+ \cup \{0\}$ (or $Z_n$ or $R^+ \cup \{0\}$ or $Z^+ \cup \{0\}\}$ be a row interval matrix semigroup under multiplication. Let $P = \{([0, a_1], [0, a_2], \ldots, [0, a_n])\} \subseteq V$, if P under the operations of V is matrix semigroup then we call P to be a row interval matrix subsemigroup of V under multiplication.*

We illustrate this by some examples.

*Example 2.31:* Let $V = \{([0, a_1], [0, a_2], \ldots, [0, a_8]) \mid a_i \in Z_{24}\}$ be a row matrix interval semigroup under multiplication.

Choose P = {([0, $a_1$], [0, $a_2$], …, [0, $a_8$]) | $a_i \in$ {0, 2, 4, 6, …, 22}; $1 \leq i \leq 8$} $\subseteq$ V; V is a row matrix interval subsemigroup of V.

*Example 2.32:* Let V = {([0, a], [0, a], [0, a], [0, a]) | a $\in Z_5 \setminus$ {0}} be a row matrix interval semigroup under multiplication. Clearly V has a row matrix interval subsemigroup.

In view of this we say a row matrix interval semigroup is simple if it has no proper row matrix interval subsemigroups. We have a large class of simple row matrix interval semigroups. We say proper if the row interval matrix semigroup is not {([0, 0], …, [0, 0])} or {([0, 1], [0, 1], …, [0, 1])}. These two semigroups will be known as improper row matrix interval subsemigroup or trivial row matrix interval subsemigroup.

*Example 2.33:* Let V = {([0, $a_1$], [0, $a_2$], …, [0, $a_9$]) | $a_i \in Z_{240}$; $1 \leq i \leq 9$} be a row matrix interval semigroup under multiplication.

P = {([0, $a_1$], [0, $a_2$], …, [0, $a_9$]) | $a_i \in$ {0, 10, 20, 30, 40, …, 230}; $1 \leq i \leq 9$} $\subseteq$ V is a row matrix interval subsemigroup of V.

**THEOREM 2.2:** *Let V = {([0, a], [0, a], …, [0, a]) | a $\in Z_p \setminus$ {0}} p a prime be a $1 \times n$ row interval matrix semigroup under multiplication. V is not a simple $1 \times n$ row interval matrix semigroup.*

The proof is left as an exercise for the reader.

**THEOREM 2.3:** *Let V = {([0, $a_1$], [0, $a_2$], …, [0, $a_n$]) | $a_i \in Z_n$; n not a prime} be a $1 \times n$ row interval matrix semigroup. V is not a simple $1 \times n$ row interval matrix semigroup (V has proper subsemigroup).*

Now we can define ideals of a $1 \times n$ row interval matrix semigroup as follows.

**DEFINITION 2.9:** *Let V be a 1 × n row interval matrix semigroup. P ⊆ V; be a proper subsemigroup. We say P is a 1 × n row interval matrix ideal of V if for all p ∈ P and v ∈ V, pv and vp are in P.*

We will illustrate this situation by some examples.

***Example 2.34:*** Let V = {([0, $a_1$], [0, $a_2$], [0, $a_3$], [0, $a_4$], [0, $a_5$]) / $a_i$ ∈ $Z_{12}$; 1 ≤ i ≤ 5} be a 1 × 5 row interval matrix semigroup. Choose P = {([0, $a_1$], [0, $a_2$], …, [0, $a_5$]) | $a_i$ ∈ {0, 2, 4, 6, 8, 10} ⊆ $Z_{12}$} ⊆ V to be a row interval matrix subsemigroup of V.

It is easily verified V is a row interval matrix ideal of V.

***Example 2.35:*** Let V = {([0, $a_1$], [0, $a_2$], …, [0, $a_9$]) | $a_i$ ∈ $Z^+$ ∪ {0}; 1 ≤ i ≤ 9} be a 1 × 9 row matrix interval semigroups. Choose P = {([0, $a_1$], [0, $a_2$], …, [0, $a_9$]) | $a_i$ ∈ $5Z^+$ ∪ {0}; 1 ≤ i ≤ 9} ⊆ V is a 1 × 9 row interval matrix ideal of V.

***Example 2.36:*** Let V = {([0, $a_1$], [0, $a_2$], …, [0, $a_n$]) | $a_i$ ∈ $Z_p$; 1 ≤ i ≤ n} be a 1 × n row interval matrix semigroup; p a prime. Clearly V has proper row interval matrix ideals. However (0) is a trivial row interval matrix ideal of V.

If V has no proper row interval matrix ideal then we call V to be a ideally simple row interval matrix semigroup. We have an infinite class of interval matrix semigroups which are not ideally simple row interval matrix semigroups.

**THEOREM 2.4:** *Let V = {([0, $a_1$], [0, $a_2$], …, [0, $a_n$]) / $a_i$ ∈ $Z_p$; p a prime; 1 ≤ i ≤ n} be a 1 × n row interval matrix semigroup. V is not an idealy simple 1 × n row interval matrix semigroup.*

*Since we need p to be prime, we have infinite number of interval row matrix semigroups which are not ideally simple (number of primes is infinite).*

How ever we give some more examples of these and subsemigroups in them under addition.

*Example 2.37:* Let

$$V = \left\{ \begin{bmatrix} [0,a_1] \\ [0,a_2] \\ \vdots \\ [0,a_{10}] \end{bmatrix} \middle| a_i \in Z_{20}; 1 \leq i \leq 10 \right\}$$

be a $10 \times 1$ column interval matrix semigroup. Clearly V is of finite order for V has only finite number of elements in them.

*Example 2.38:* Let

$$V = \left\{ \begin{bmatrix} [0,a_1] \\ [0,a_2] \\ [0,a_3] \\ [0,a_4] \\ [0,a_5] \\ [0,a_6] \end{bmatrix} \middle| a_i \in Z^+ \cup \{0\}; 1 \leq i \leq 6 \right\}$$

be a $6 \times 1$ column interval matrix semigroup under addition. Clearly V is of infinite order.

*Example 2.39:* Let

$$V = \left\{ \begin{bmatrix} [0,a_1] \\ [0,a_1] \\ [0,a_1] \\ [0,a_1] \end{bmatrix} \middle| a_1 \in Z_7 \right\}$$

be a $4 \times 1$ column interval matrix semigroup under addition. V is of order seven.

We have seen examples of these semigroups. The notion of subsemigroups can be defined as in case of row matrix interval semigroups. So we leave this simple task to the reader but give examples of these substructures.

*Example 2.40:* Let

$$V = \left\{ \begin{bmatrix} [0,a_1] \\ [0,a_2] \\ \vdots \\ [0,a_{12}] \end{bmatrix} \middle| a_i \in Z_{30}; 1 \le i \le 12 \right\}$$

be a $12 \times 1$ column interval matrix semigroup. Take

$$P = \left\{ \begin{bmatrix} [0,a_1] \\ [0,a_2] \\ \vdots \\ [0,a_{12}] \end{bmatrix} \middle| a_i \in \{0,2,4,6,8,10,12,\ldots,28\} \subseteq Z_{30} \right\} \subseteq V;$$

P is a $12 \times 1$ column interval matrix subsemigroup of V.

*Example 2.41:* Let

$$V = \left\{ \begin{bmatrix} [0,a_1] \\ [0,a_2] \\ \vdots \\ [0,a_{12}] \end{bmatrix} \middle| a_i \in Z^+ \cup \{0\}; 1 \le i \le 12 \right\}$$

be a $12 \times 1$ column interval matrix semigroup.

$$W = \left\{ \begin{bmatrix} [0,a_1] \\ [0,a_2] \\ \vdots \\ [0,a_{12}] \end{bmatrix} \middle| a_i \in 7Z^+ \cup \{0\}; 1 \le i \le 12 \right\} \subseteq V;$$

W is a 12 × 1 column interval matrix subsemigroup of V.

*Example 2.42:* Let

$$V = \left\{ \begin{bmatrix} [0,a_1] \\ [0,a_2] \\ \vdots \\ [0,a_7] \end{bmatrix} \middle| a_i \in Z_{36}; 1 \le i \le 7 \right\}$$

be a 7 × 1 column interval matrix semigroup.
Take

$$I = \left\{ \begin{bmatrix} [0,a_1] \\ [0,a_2] \\ \vdots \\ [0,a_7] \end{bmatrix} \middle| a_i \in \{0,2,4,...,34\} \subseteq Z_{36}; 1 \le i \le 7 \right\} \subseteq V;$$

I is a 7 × 1 column interval matrix subsemigroup of V.

*Example 2.43:* Let

$$V = \left\{ \begin{bmatrix} [0,a_1] \\ [0,a_2] \\ \vdots \\ [0,a_9] \end{bmatrix} \middle| a_i \in Z^+ \cup \{0\}; 1 \le i \le 9 \right\}$$

be a 9 × 1 column matrix interval semigroup. Take

$$I = \left\{ \begin{bmatrix} [0,a_1] \\ [0,a_2] \\ \vdots \\ [0,a_9] \end{bmatrix} \middle| a_i \in 7Z^+ \cup \{0\}; 1 \le i \le 9 \right\} \subseteq V,$$

I is a 9 × 1 column interval matrix subsemigroup of V of infinite order.

*Example 2.44:* Let

$$V = \left\{ \begin{bmatrix} [0, a_1] \\ [0, a_2] \\ \vdots \\ [0, a_{11}] \end{bmatrix} \middle| a_i \in Z_{12}; 1 \leq i \leq 11 \right\}$$

be a 11 × 1 column matrix interval semigroup.

$$I = \left\{ \begin{bmatrix} [0, a_1] \\ [0, a_2] \\ \vdots \\ [0, a_{11}] \end{bmatrix} \middle| a_i \in \{0, 3, 6, 9\} \subseteq Z_{12}; 1 \leq i \leq 11 \right\} \subseteq V;$$

is a 11 × 1 column matrix interval subsemigroup of V.

*Example 2.45:* Let

$$V = \left\{ \begin{bmatrix} [0, a_1] \\ [0, a_2] \\ \vdots \\ [0, a_8] \end{bmatrix} \middle| a_i \in Z_3; 1 \leq i \leq 8 \right\}.$$

V is a 8 × 1 column matrix interval semigroup.
Take

$$W = \left\{ \begin{bmatrix} [0, a] \\ [0, a] \\ \vdots \\ [0, a] \end{bmatrix} \middle| a \in Z_3 \right\} \subseteq V;$$

W is a 8 × 1 column matrix interval subsemigroup of V.

We can define m × n matrix interval semigroup.

Let V = {M = ($m_{ij}$) | $m_{ij}$ = [0, $a_{ij}$]; 1 ≤ i ≤ n and 1 ≤ j ≤ m, $a_{ij}$ ∈ $Z^+ \cup \{0\}$ (or $Z_n$ or $R^+ \cup \{0\}$ or $Q^+ \cup \{0\}$)} be a collection of n × m interval matrices. Define on V matrix addition i.e., if M = ([0, $a_{ij}$]) and N = ([0, $b_{ij}$]) then M + N = ([0, $a_{ij} + b_{ij}$])

V under interval matrix addition is a semigroup called the n × m matrix interval semigroup. If m = n then V can be a semigroup under multiplication as well as addition.

We will describe both the operation with some interval matrices. Let

$$A = \begin{bmatrix} [0,5] & [0,1] & [0,3] \\ [0,2] & [0,4] & [0,7] \\ [0,1] & [0,6] & [0,5] \\ [0,0] & [0,2] & [0,8] \end{bmatrix}$$

and

$$B = \begin{bmatrix} [0,1] & [0,2] & [0,3] \\ [0,4] & [0,5] & [0,6] \\ [0,7] & [0,8] & [0,1] \\ [0,2] & [0,4] & [0,5] \end{bmatrix}$$

be interval matrices with entries from $Z_9$. Now

$$A + B = \begin{bmatrix} [0,5] & [0,1] & [0,3] \\ [0,2] & [0,4] & [0,7] \\ [0,1] & [0,6] & [0,5] \\ [0,0] & [0,2] & [0,8] \end{bmatrix} + \begin{bmatrix} [0,1] & [0,2] & [0,3] \\ [0,4] & [0,5] & [0,6] \\ [0,7] & [0,8] & [0,1] \\ [0,2] & [0,4] & [0,5] \end{bmatrix}$$

$$= \begin{bmatrix} [0,6] & [0,3] & [0,6] \\ [0,6] & [0,0] & [0,4] \\ [0,8] & [0,5] & [0,6] \\ [0,2] & [0,6] & [0,4] \end{bmatrix}.$$

Clearly the product is not defined.
Now $[0, a] \times [0, b] = [0, ab]$.
 If we take
$$A = \begin{bmatrix} [0,5] & [0,7] \\ [0,1] & [0,4] \end{bmatrix}$$
and
$$B = \begin{bmatrix} [0,3] & [0,1] \\ [0,5] & [0,8] \end{bmatrix}$$

with entries from $Z^+ \cup \{0\}$, then

$$AB = \begin{bmatrix} [0,5] & [0,7] \\ [0,1] & [0,4] \end{bmatrix} \begin{bmatrix} [0,3] & [0,1] \\ [0,5] & [0,8] \end{bmatrix}$$

$$= \begin{bmatrix} [0,5][0,3]+[0,7][0,5] & [0,5][0,1]+[0,7][0,8] \\ [0,1][0,3]+[0,4][0,5] & [0,1][0,1]+[0,4][0,8] \end{bmatrix}$$

$$= \begin{bmatrix} [0,15]+[0,35] & [0,5]+[0,56] \\ [0,3]+[0,20] & [0,1]+[0,32] \end{bmatrix}$$

$$= \begin{bmatrix} [0,50] & [0,61] \\ [0,23] & [0,33] \end{bmatrix}.$$

Thus interval matrix addition and multiplication are well defined.
 Now we will examples of these structures.

***Example 2.46:*** Let
$$V = \left\{ \begin{bmatrix} [0,a_1] & [0,a_2] \\ [0,a_3] & [0,a_4] \\ [0,a_5] & [0,a_6] \\ [0,a_7] & [0,a_8] \end{bmatrix} \middle| a_i \in Z_{20}; 1 \le i \le 8 \right\}$$

be a 4 × 2 interval matrix semigroup under interval matrix addition.

*Example 2.47:* Let

$$V = \left\{ \begin{bmatrix} [0,a_1] & [0,a_4] & [0,a_7] & [0,a_{10}] \\ [0,a_2] & [0,a_5] & [0,a_8] & [0,a_{11}] \\ [0,a_3] & [0,a_6] & [0,a_9] & [0,a_{12}] \end{bmatrix} \middle| a_i \in Z^+ \cup \{0\}; 1 \leq i \leq 12 \right\}$$

be a 3 × 4 interval matrix semigroup under addition.

*Example 2.48:* Let

$$V = \left\{ \begin{bmatrix} [0,a_1] & [0,a_2] & [0,a_3] \\ [0,a_4] & [0,a_5] & [0,a_6] \\ [0,a_7] & [0,a_8] & [0,a_9] \end{bmatrix} \middle| a_i \in Z_{30}; 1 \leq i \leq 9 \right\}$$

be a 3 × 3 interval matrix semigroup under multiplication.

Thus we can as in case of other interval semigroups define interval matrix subsemigroups and ideals.

This task of defining and giving examples is left as an exercise for the reader.

Now having seen interval matrix semigroups we now put forth some of the important properties about these semigroups.

An interval matrix semigroup V is said to be an interval matrix Smarandache semigroup (interval matrix S-semigroup) or Smarandache matrix interval semigroup (S-matrix interval semigroup) if V has a proper subset P where P is a group under the operations of V. We say V is a Smarandache commutative matrix interval semigroup if every proper subset A of V which is a group under the operations of V is a commutative matrix interval group.

If only one subset A of V is a group and is commutative we call V to be a weakly commutative matrix interval S-semigroup.

*Example 2.49:* Let

$$V = \left\{ \begin{bmatrix} [0,a_1] & [0,a_2] \\ [0,a_3] & [0,a_4] \end{bmatrix} \middle| a_i \in Z_{12}; 1 \le i \le 4 \right\}$$

be the matrix interval semigroup under matrix multiplication.

Clearly V has atleast one subset

$$P = \left\{ \begin{bmatrix} [0,a] & [0,a] \\ [0,a] & [0,a] \end{bmatrix} \middle| a \in Z_{12} \right\} \subseteq V;$$

P is a matrix interval commutative subsemigroup of V; hence P is a weakly commutative matrix interval semigroup.

*Example 2.50:* Let $V = \{([0, a_1], [0, a_2], \ldots, [0, a_n]) \mid a_i \in Z^+ \cup \{0\}; 1 \le i \le n\}$ be a row matrix interval semigroup. Clearly V is a row matrix interval commutative semigroup.

*Example 2.51:* Let $V = \{([0, a_1], [0, a_2], \ldots, [0, a_n]) / a_i \in Z_7 \setminus \{0\}; 1 \le i \le n\}$ be a row matrix interval semigroup. Take $W = \{([0, a], [0, a], \ldots, [0, a]) \mid a \in Z_7 \setminus \{0\}\} \subseteq V$.

W is a row matrix interval group of V under the operations of V. Hence V is a Smarandache row matrix interval semigroup.

**THEOREM 2.5:** *Let $V = \{([0, a_1], [0, a_2], \ldots, [0, a_n]) / a_i \in Z_p \setminus \{0\}$; p is a prime; $1 \le i \le n\}$ be a row matrix interval semigroup. Take $W = \{([0, a], [0, a], \ldots, [0, a]) \mid a \in Z_p \setminus \{0\}\} \subseteq V$; W is a row matrix interval group. Hence V is a row matrix interval Smarandache semigroup.*

The proof is left as an exercise to the reader.

**THEOREM 2.6:** *Let $V = \{([0, a_1], [0, a_2], \ldots, [0, a_n]) / a_i \in Z_m\}$ m a composite number be a row matrix interval semigroup. If $Z_m$ is a Smarandache semigroup with $P \subseteq Z_m$; P a group then $W = \{([0, a], [0, a], \ldots, [0, a]) \mid a \in P \subseteq Z_m\} \subseteq V$; W is a interval* *group. Thus V is a row matrix interval Smarandache semigroup.*

This proof is also left as an exercise for the reader.

***Example 2.52:*** Let $Z_{30} = \{0, 1, 2, \ldots, 29\}$ be a semigroup under multiplication modulo 30. $V = \{([0, a_1], [0, a_2], \ldots, [0, a_9]) \mid a_i \in Z_{30} ; 1 \le i \le 9\}$ be a row interval matrix semigroup. $W = \{([0, a_1], [0, a_2], \ldots, [0, a_9]) \mid a_i \in \{0, 5, 10, 15, 20, 25\} \subseteq Z_{30}\} \subseteq V$. W is a row interval matrix ideal of V.

**THEOREM 2.7:** *Let $V = \{([0, a_1], [0, a_2], \ldots, [0, a_n]) / a_i \in Z_p \}$ be a row interval matrix semigroup under multiplication; V has no proper ideals.*

This proof is also left for the reader.

Let
$$V = \left\{ \begin{bmatrix} [0,a_1] & \cdots & [0,a_n] \\ [0,b_1] & \cdots & [0,b_n] \\ \vdots & & \vdots \\ [0,c_1] & \cdots & [0,c_n] \end{bmatrix} \middle| a_i, b_i, c_i \in Z_n; 1 \le i \le n \right\}$$

be the collection of all $n \times n$ interval square matrix. V is a square matrix interval semigroup under multiplication. (or addition, or used in the mutually exclusive sense).

***Example 2.53:*** Let
$$V = \left\{ \begin{bmatrix} [0,a_1] & [0,a_2] \\ [0,a_3] & [0,a_4] \end{bmatrix} \middle| a_i \in Z_4; 1 \le i \le 4 \right\}$$

be a $4 \times 4$ interval matrix semigroup under addition modulo 4. (V is also a $4 \times 4$ interval matrix semigroup under multiplication). For take

$$A = \begin{pmatrix} [0,3] & [0,1] \\ [0,2] & [0,2] \end{pmatrix}$$

and

$$B = \begin{pmatrix} [0,1] & [0,2] \\ [0,2] & [0,3] \end{pmatrix}$$

in V.

$$A \times B = \begin{pmatrix} [0,3] & [0,1] \\ [0,2] & [0,2] \end{pmatrix} \times \begin{pmatrix} [0,1] & [0,2] \\ [0,2] & [0,3] \end{pmatrix}$$

$$= \begin{pmatrix} [0,3][0,1]+[0,1][0,2] & [0,3][0,2]+[0,1][0,3] \\ [0,2][0,1]+[0,2][0,2] & [0,2][0,2]+[0,2][0,3] \end{pmatrix}$$

$$= \begin{pmatrix} [0,3]+[0,2] & [0,2]+[0,3] \\ [0,2]+[0,0] & [0,0]+[0,2] \end{pmatrix}$$

$$= \begin{pmatrix} [0,1] & [0,1] \\ [0,2] & [0,2] \end{pmatrix}.$$

We can define the notion of Smarandache Lagrange semigroup, Smarandache subsemigroup, Smarandache hypersubsemigroup Smarandache p-Sylow subgroup, Smarandache Cauchy elements of a S-semigroup and Smarandache coset in case of interval matrix semigroup in an analogous way [9].

We will illustrate these situations by examples for the definition is very similar to that of semigroups [9].

*Example 2.54:* Let

$$V = \left\{ \begin{bmatrix} [0,a_1] & [0,a_2] & [0,a_3] \\ [0,a_4] & [0,a_1] & [0,a_6] \\ [0,a_7] & [0,a_8] & [0,a_1] \end{bmatrix} \middle| a_i = 0; i = 2,3,4,6,7,8 \right\}$$

be a 3 × 3 interval matrix semigroup under multiplication.

$$W = \left\{ A = \begin{bmatrix} [0,a_1] & [0,a_2] & [0,a_3] \\ [0,a_4] & [0,a_5] & [0,a_6] \\ [0,a_7] & [0,a_8] & [0,a_9] \end{bmatrix} \middle| |A| \neq 0; a_1 \in \{1,3\} \subseteq Z_8; \; a_i = 0 \text{ if } i \neq 1 \right\}$$

⊆ V is a subgroup. Clearly V is a Smarandache matrix interval semigroup.

*Example 2.55:* Let

$$V = \left\{ \begin{bmatrix} [0,a_1] & [0,a_2] \\ [0,a_3] & [0,a_4] \end{bmatrix} \middle| a_i \in Z_9; 1 \leq i \leq 4 \right\}$$

be a square matrix interval semigroup.
Take

$$W = \left\{ \begin{bmatrix} [0,a_1] & [0,a_2] \\ [0,a_3] & [0,a_1] \end{bmatrix} \middle| \begin{array}{l} a_1 \in \{1,8\}, a_i = 0 \subseteq Z_9; \\ |A| \neq 0; 2 \leq i \leq 3 \end{array} \right\} \subseteq V;$$

V is a square matrix interval group under multiplication. Thus V is a Smarandache square matrix interval semigroup.

Now we proceed onto give examples of Smarandache matrix interval subsemigroup or matrix interval Smarandache semigroup.

*Example 2.56:* Let

$$V = \left\{ \begin{bmatrix} [0,a_1] & [0,a_2] & [0,a_3] & [0,a_4] \\ [0,a_5] & [0,a_6] & [0,a_7] & [0,a_8] \\ [0,a_9] & [0,a_{10}] & [0,a_{11}] & [0,a_{12}] \\ [0,a_{13}] & [0,a_{14}] & [0,a_{15}] & [0,a_{16}] \end{bmatrix} \middle| a_i \in Z_{11}; 1 \leq i \leq 16 \right\}$$

be a 4 × 4 matrix interval semigroup.
Take

$$W = \left\{ \begin{bmatrix} [0,a_1] & [0,a_2] & [0,a_3] & [0,a_4] \\ [0,a_5] & [0,a_1] & [0,a_7] & [0,a_8] \\ [0,a_9] & [0,a_{10}] & [0,a_1] & [0,a_{12}] \\ [0,a_{13}] & [0,a_{14}] & [0,a_{15}] & [0,a_1] \end{bmatrix} \middle| \begin{array}{l} |A| \neq 0; a_i = 0; \\ a_1 \in \{1,...,10\} \subseteq Z_{11} \\ i = 2,3,4,5,7,8,9, \\ 10,12,13,14,15 \end{array} \right\}$$

⊆ V; W is a interval matrix group.

Now take P = {All 4 × 4 square matrices with intervals of the form [0, $a_i$]; where $a_i \in Z_{11} \setminus \{0\}$} ⊆ V; P is a interval matrix subsemigroup of V and W ⊆ P; so P is a matrix interval Smarandache subsemigroup of V.

*Example 2.57:* Let

$$V = \left\{ \begin{bmatrix} [0,a_1] & [0,a_2] \\ [0,a_3] & [0,a_4] \end{bmatrix} \middle| a_i \in Z_{15}; 1 \leq i \leq 4 \right\}$$

be a interval matrix semigroup.
Take

$$W = \left\{ \begin{bmatrix} [0,a_1] & [0,a_2] \\ [0,a_3] & [0,a_4] \end{bmatrix} \middle| \begin{array}{l} a_1 \in \{0,3,6,9,12\} \subseteq Z_{15}; \\ a_i = 0; 2 \leq i \leq 4 \end{array} \right\} \subseteq V$$

to be a interval matrix subsemigroup of V.

Let

$$P = \left\{ A = \begin{bmatrix} [0,a_1] & [0,a_2] \\ [0,a_3] & [0,a_1] \end{bmatrix} \middle| \begin{array}{l} a_1 \in \{0,3,6,9,12\} \subseteq Z_{15}; \\ |A| \neq 0; a_2 = a_3 = 0 \end{array} \right\}$$

⊆ W; W is a interval matrix smarandache subsemigroup of V.

Now we proceed onto give examples of the notion of Smarandache interval matrix subsemigroup.

***Example 2.58:*** Let V = {set all 5 × 5 interval matrices with intervals of the form [0, $a_i$] with $a_i \in Z_{43}$} be a interval matrix semigroup.

Take W = {A / all 5 × 5 interval matrices with intervals of the form [0, $a_1$] with $a_1 \in Z_{43} \setminus \{0\}$ such that $|A| \neq 0$. A is a 5 × 5 diagonal interval matrix} $\subseteq$ V; W is a group under interval matrix multiplication. So V is a interval matrix Smarandache semigroup.

Further if we take P = {all 5 × 5 diagonal interval matrices with intervals of the form [0, $a_i$] with $a_i \in Z_{43} \setminus \{0\}$} $\subseteq$ V then P is a interval matrix subsemigroup of V.

We see W $\subseteq$ P and W is the largest interval matrix group present in P. Thus P is a matrix interval Smarandache subsemigroup of V.

***Example 2.59:*** Let

$$V = \left\{ \begin{bmatrix} [0,a_1] & [0,a_2] \\ [0,a_3] & [0,a_4] \end{bmatrix} \middle| a_i \in Z_{11}; 1 \leq i \leq 4 \right\}$$

be a matrix interval semigroup.

Let

$$P = \left\{ \begin{bmatrix} [0,a_1] & [0,a_2] \\ [0,a_3] & [0,a_1] \end{bmatrix} \middle| a_1 \in Z_{11}; |A| \neq 0; a_2 = a_3 = 0 \right\} \subseteq V;$$

P is the largest interval matrix group present in V; but V has no proper matrix interval subsemigroup which contains P.

***Example 2.60:*** Let V = {[0, $a_i$] | $a_i \in Z_{43}$} be a matrix interval semigroup. P = {[0, $a_i$] | $a_i \in Z_{43} \setminus \{0\}$} is the matrix interval

subgroup of V. Infact V has no proper matrix interval subsemigroup containing P.

***Example 2.61***: Let $V = \{[0, a_i] / a_i \in Z_6\}$ be the matrix interval semigroup. Take $W = \{[0, 1], [0, 5]\} \subseteq V$ is a interval subgroup of V.

Take $P = \{[0, 1], [0, 5], [0, 0]\} \subseteq V$ is a interval Smarandache subsemigroup we see P is a interval Smarandache hyper subsemigroup of V.

It is left for the reader to prove the following theorems.

**THEOREM 2.8:** *Let $V = \{$all $n \times n$ diagonal interval matrices with intervals of the form $[0, a_1]$, $a_1 \in Z_p$; all diagonal elements are the same$\}$ is a Smarandache simple interval matrix semigroup which is a Smarandache interval matrix semigroup.*

**THEOREM 2.9:** *Let V be a Smarandache matrix interval semigroup.*
*Every Smarandache matrix interval hyper subsemigroup is a Smarandache matrix interval subsemigroup but every Smarandache matrix interval subsemigroup in general is not a S-matrix interval hyper subsemigroup.*

Now we proceed onto give examples of Smarandache matrix interval Lagrange semigroup (S-matrix interval Lagrange semigroup).

***Example 2.62:*** Let $V = \{[0, a] \mid a \in Z_4\}$ be a S interval semigroup. $A = \{[0, 1], [0, 3]\} \subseteq V$ is a interval subgroup of V.
Clearly o(A) / o(V) so V is a S-matrix interval Lagrange semigroup.

***Example 2.63:*** Let $V = \{[0, a] \mid a \in Z_9\}$ be a S-matrix interval semigroup. Let $A = \{[0, 1], [1, 8]\} \subseteq V$ and $B = \{[0, 1], [0, 2], [0, 4], [0, 5], [0, 7], [0, 8]\} \subseteq V$ be subgroup of V.

We see both of them do not divide of the order of V. So V is not a S-interval Lagrange semigroup.

***Example 2.64:*** Let $V = \{[0, a] \mid a \in Z_{10}\}$ be S-interval semigroup. V is a S-weakly Lagrange interval semigroup.

The proof of the following theorem is left as an exercise for the reader.

**THEOREM 2.10:** *Every S-interval Lagrange semigroup is a S-interval weakly Lagrange semigroup.*

Next we proceed onto illustrate S-p-Sylow interval subgroup of a S-interval semigroup.

***Example 2.65:*** Let $V = \{[0, a] \mid a \in Z_{16}\}$ be a S-interval semigroup. $A = \{[0, 1], [0, 9]\} \subseteq V$ is a interval subgroup of V. $2 / o(V)$ but $2^2 / o(V)$, but V has S-2-Sylow interval subgroups of order 4 given by $B = \{[0, 6], [0, 2], [0, 4], [0, 8]\} \subseteq V; 4 / o(V)$.

We see in case of S-interval semigroup V we say if p is a prime such that $p / o(V)$ then we can have interval subgroup of order $p^\alpha$; where $p^\alpha / o(V)$, we call such intervals subgroups of the S-interval semigroup to be S-p-Sylow interval subgroups of V.

We give examples of S-Cauchy element of a interval semigroup.
We see a S-Cauchy element of a interval semigroup x of V is such that $x^t = 1$ and $t / o(V)$.

***Example 2.66:*** Let $V = \{[0, a] \mid a \in Z_{19}\}$ be a S-interval semigroup. Take $x = [0, 18] \in V; x^2 = ([0, 18])^2 = [0, 1]$.
We see $2 \nmid o(V)$. Thus x is not a S-Cauchy interval element of V.

We leave the proof of the following theorem of the reader.

**THEOREM 2.11:** *Let $V = \{[0, a] / a \in Z_p\}$; (p a prime) be the S-interval semigroup under multiplication. No element of V is a S-Cauchy element of V.*

The proof is obvious from the fact that no integer n can divide the prime p. Hence the claim.

# Chapter Three

# INTERVAL POLYNOMIAL SEMIGROUPS

In this chapter we introduce the notion of interval polynomial semigroups. We call a polynomial in the variable x to be an interval polynomial if the coefficients of x are intervals of the form $[0, a_i] / a_i \in Z_p$ (or $Z_n$ or $Z^+ \cup \{0\}$ or $Q^+ \cup \{0\}$ or $R^+ \cup \{0\}$.

$[0, 5] + [0, 7]x + [0, 2] x^3 + [0, 14] x^9 = p(x)$ is a interval polynomial in the variable x.

We now define interval polynomial semigroup under addition (or multiplication).

We just illustrate how interval polynomials are added.

Let $p(x) = [0, 2] + [0, 3] x^2 + [0, 7] x^7 + [0, 11] x^9$ and $q(x) = [0, 12] + [0, 7] x + [0, 14] x^3 + [0, 10] x^7 + [0, 5] x^8 + [0, 12] x^9 + [0, 5] x^{20}$.

$$p(x) + q(x) = ([0, 2] + [0, 3] x^2 + [0, 7] x^7 + [0, 11] x^9 + [0, 12] + [0, 7] x + [0, 14] x^3 + [0, 10] x^7 + [0, 5] x^8 + [0, 12] x^9 + [0, 5] x^{20})$$

$$= ([0, 2] + [0, 12]) + [0, 7] x + [0, 3] x^2 + [0, 14]x^3 + ([0, 7] x^7 + [0, 10] x^7) + [0, 5]x^8 + ([0, 11] x^9 + [0, 12] x^9) + [0, 5] x^{20}$$

$$= [0, 14] + [0, 7]x + [0, 3] x^2 + [0, 14] x^3 + [0, 17] x^7 + [0, 5] x^8 + [0, 23] x^9 + [0, 5] x^{20}.$$

Now we will just define interval polynomial multiplication.
$$p(x) = [0, 3] + [0, 5] x^2 + [0, 11] x^5$$
and
$$q(x) = [0, 8] + [0, 1] x + [0, 9] x^3.$$

$$p(x).q(x) = ([0, 3] + [0, 5] x^2 + [0, 11] x^5) ([0, 8] + [0, 1]x + [0, 9] x^3).$$

$$= [0, 3] [0, 8] + [0, 5] [0, 8] x^2 + [0, 11] [0, 8] x^5 + [0, 3] [0, 1] x + [0, 5] x^2 [0, 1] x + [0, 11] x^5 [0, 1] x + [0, 3] [0, 9] x^3 + [0, 5] x^2 [0, 9] x^3 + [0, 11] x^5 [0, 9] x^3.$$

$$= [0, 24] + [0, 40] x^2 + [0, 88] x^5 + [0, 3] x + [0, 5] x^3 + [0, 11] x^6 + [0, 27]x^3 + [0, 45] x^5 + [0, 99] x^8.$$

$$= [0, 24] + [0, 3]x + [0, 40]x^2 + [0, 32]x^3 + [0, 45]x^5 + [0, 11]x^6 + [0, 99] x^8.$$

Now having defined interval polynomial addition and multiplication we proceed onto define interval polynomial semigroup under these operations.

**DEFINITION 3.1:** *Let* $S = \left\{ \sum_{i=0}^{n} [0, a_i] x^i \,\middle|\, a_i \in Z_n \text{ (or } Z^+ \cup \{0\}$ *or* $R^+ \cup \{0\}$, *or* $Q^+ \cup \{0\}$, $C^+ \cup \{0\}$) *and x is a variable or*

*indeterminate} S under addition of interval polynomials is a semigroup defined as interval polynomial semigroup.*

We will illustrate this situation by some examples.

*Example 3.1:* Let

$$S = \left\{ \sum_{i=0}^{9} [0, a_i] x^i \;\middle|\; a_i \in Z^+ \cup \{0\} \right\}$$

be a interval polynomial semigroup under addition. Clearly the number of elements in S is infinite so S is an infinite order interval polynomial semigroup.

*Example 3.2:* Let

$$S = \left\{ \sum_{i=0}^{3} [0, a_i] x^i \;\middle|\; a_i \in Z_{11} \right\}$$

be a interval polynomial semigroup. Clearly S is of finite order.

We see clearly the interval polynomial semigroups given in the above examples are not compatible under multiplication.

*Example 3.3:* Let

$$S = \left\{ \sum_{i=0}^{7} [0, a_i] x^i \;\middle|\; a_i \in R^+ \cup \{0\} \right\}$$

be the interval polynomial semigroup. S is an infinite interval polynomial semigroup under addition.

Now we can define substructures for these structures.

**DEFINITION 3.2:** *Let S be a interval polynomial semigroup under addition. Suppose W ⊆ S be a proper subset of S and if W*

*is itself an interval polynomial semigroup under addition then we define W to be an interval polynomial subsemigroup of S.*

We will illustrate this situation also by some examples.

*Example 3.4:* Let

$$S = \left\{ \sum_{i=0}^{8} [0, a_i] x^i \,\Big|\, a_i \in Z^+ \cup \{0\} \right\}$$

be a interval polynomial semigroup under addition. Take

$$W = \left\{ \sum_{i=0}^{8} [0, a_i] x^i \,\Big|\, a_i \in 3Z^+ \cup \{0\} \right\} \subseteq S;$$

W is a interval polynomial subsemigroup of S under addition.

*Example 3.5:* Let

$$S = \left\{ \sum_{i=0}^{20} [0, a_i] x^i \,\Big|\, a_i \in Z_{12} \right\}$$

be a interval polynomial semigroup under addition. Take

$$W = \left\{ \sum_{i=0}^{10} [0, a_i] x^i \,\Big|\, a_i \in Z_{12} \right\} \subseteq S;$$

W is a interval polynomial subsemigroup of S under addition.
    Both S and W are of finite order.

   We can define ideals as in case of usual semigroups. If a interval polynomial semigroup S has no proper interval polynomial subsemigroups we call S to be a simple interval polynomial semigroup.

We now proceed onto define polynomial interval semigroup under multiplication.

**DEFINITION 3.3:** *Let*

$$V = \left\{ \sum_{i=0}^{\infty} [0,a_i]x^i \,\middle|\, a_i \in Z^+ \cup \{0\} \right.$$

*(or $Z_n$ or $R^+ \cup \{0\}$, or $Q^+ \cup \{0\}$, x a variable} be a collection of interval polynomials. If product is defined on V then V is a interval polynomial semigroup under multiplication.*

We will illustrate this situation by some examples.

*Example 3.6:* Let

$$S = \left\{ \sum_{i=0}^{\infty} [0,a_i]x^i \,\middle|\, a_i \in Z_8 \right\}$$

be a polynomial interval semigroup under multiplication. Clearly S is of infinite order.

*Example 3.7:* Let

$$S = \left\{ \sum_{i=0}^{\infty} [0,a_i]x^i \,\middle|\, a_i \in R^+ \cup \{0\} \right\}$$

be a interval polynomial semigroup under multiplication. Clearly S is of infinite order.

Substructure is defined as in case of usual semigroups.
However we will illustrate this situation by some examples.

*Example 3.8:* Let

$$S = \left\{ \sum_{i=0}^{\infty} [0,a_i]x^i \,\middle|\, a_i \in Q^+ \right\}$$

be a interval polynomial semigroup.
Take

$$P = \left\{ \sum_{i=0}^{\infty} [0, a_i] x^i \;\middle|\; a_i \in Z^+ \right\} \subseteq S;$$

P is a interval polynomial subsemigroup of S. Clearly P is a not as interval polynomial ideal of S.

*Example 3.9:* Let

$$S = \left\{ \sum_{i=0}^{\infty} [0, a_i] x^i \;\middle|\; a_i \in Z_{30} \right\}$$

be a interval polynomial semigroup under multiplication.
Take

$$T = \left\{ \sum_{i=0}^{\infty} [0, a_i] x^i \;\middle|\; a_i \in \{0, 2, 4, 6, ..., 26, 28\} \subseteq Z_{30} \right\} \subseteq S;$$

T is a interval polynomial subsemigroup of S. T is also a interval polynomial ideal of S.
   Thus we can have interval polynomial subsemigroups which are not interval ideals of the polynomial semigroup.

*Example 3.10:* Let

$$S = \left\{ \sum_{i=0}^{\infty} [0, a_i] x^i \;\middle|\; a_i \in Q^+ \cup \{0\} \right\} \subseteq S;$$

T is only a polynomial interval subsemigroup and is not a polynomial interval ideal of S.
   Infact S has no polynomial interval ideals but has infinite number of polynomial interval subsemigroups.
   Now having seen the two substructures we proceed on to consider finite polynomial interval semigroups which are Smarandache Lagrange polynomial interval semigroup, Smarandache polynomial interval semigroup and so on. It is pertinent to mention here that polynomial interval semigroup can contain Smarandache Cauchy elements.

We will illustrate this by some examples.

***Example 3.11:*** Let
$$S = \left\{\sum_{i=0}^{2}[0,a_i]x^i \;\middle|\; a_i \in Z_3\right\}$$

be a polynomial interval semigroup under addition.

S = {0, [0, 1] x, [0, 1] $x^2$, [0, 2] x, [0, 2] $x^2$, [0, 1], [0, 2], [0, 1] + [0, 1]x, [0, 1] + [0, 1] $x^2$, [0, 1] + [0, 2] x, [0, 1] + [0, 2] $x^2$, [0, 2] + [0, 1] x, [0, 2] + [0, 2] $x^2$ [0, 2] + [0, 1]$x^2$, [0, 2] + [0, 2]x, [0, 1] + [0, 1]x + [0, 1] $x^2$, …}.

We see
$$[0, 1] + [0, 1] + [0, 1] = 0$$
$$[0, 1] x + [0, 1]x + [0, 1]x = 0$$
$$[0, 2]x + [0, 2]x + [0, 2]x = 0$$

Thus we have several Cauchy elements, S is also a S-polynomial interval semigroup.

It is left as an exercise for the reader to find the order of S and find out whether the elements are S-Cauchy elements of S.

***Example 3.12:*** Let
$$S = \left\{\sum_{i=0}^{5}[0,a_i]x^i \;\middle|\; a_i \in Z_2\right\}$$

be a polynomial interval semigroup under addition.

We see S has several elements of finite order but one is to find the order of S. S is a S-polynomial interval semigroup. We see [0, 1] + [0, 1] = 0, ([0, 1]x + [0, 1]) + ([0, 1]x + [0, 1]) = 0 and so on.

The reader is left with the task of finding the order of S.
However S is a S-interval polynomial semigroup.

*Example 3.13:* Let

$$S = \left\{ \sum_{i=0}^{\infty} [0, a_i]x^i \;\middle|\; a_i \in Z_7 \right\}$$

be a polynomial interval semigroup. We see S is of infinite order (S be under addition or multiplication).
We cannot in this case define S – Cauchy element. However S is a S-polynomial interval semigroup under addition and S is a S-polynomial interval semigroup under multiplication.

*Example 3.14:* Let

$$S = \left\{ \sum_{i=0}^{8} [0, a_i]x^i \;\middle|\; a_i \in Z_8 \right\}$$

be a interval polynomial semigroup under addition. S is a S-interval polynomial semigroup. S is a S-commutative interval polynomial semigroup.

Further it is easily verified S is a S-weakly interval polynomial semigroup. For $[0, 1]x^i$ in S generates a cyclic group under addition where $1 \leq i \leq 8$.

*Example 3.15:* Let

$$S = \left\{ \sum_{i=0}^{6} [0, a_i]x^i \;\middle|\; x^7 = 1,\; x^8 = x,\text{ so on};\; a_i \in Z_6 \right\}$$

be a polynomial interval semigroup under multiplication.
For if $p(x) = [0, 1]x + [0, 5]x^5 + [0, 2]$ and $q(x) = [0, 4] + [0, 3]x^3$ in S the

$$\begin{aligned}
p(x)q(x) &= ([0, 2] + [0, 1]x + [0, 5]x^5) \times ([0, 4] + [0, 3]x^3) \\
&= [0, 2][0, 4] + [0, 2][0, 3]x^3 + [0, 1]x[0, 4] + \\
&\quad [0, 1]x[0, 3]x^3 + [0, 5]x^5[0, 4] + [0, 5]x^5 \cdot [0, 3]x^3 \\
&= [0,2] + [0,0]\,x^3 + [0,4]x + [0,3]x^4 + [0,2]\,x^5 + \\
&\quad [0,3]x. \\
&= [0,2] + [0,1]x + [0,3]x^4 + [0,2]x^5.
\end{aligned}$$

It is easily verified S is a S-interval polynomial weakly cyclic semigroup.

This simple result can be proved by the reader.

*Example 3.16:* Let

$$S = \left\{ \sum_{i=0}^{\infty} [0, a_i] x^i \;\middle|\; a_i \in Z_{12} \right\}$$

be interval polynomial semigroup under multiplication. S has interval polynomial ideals, for take

$$P = \left\{ \sum_{i=0}^{\infty} [0, a_i] x^i \;\middle|\; a_i \in \{0, 2, 4, 6, 8, 10\} \subseteq Z_{12} \right\} \subseteq S$$

is an interval polynomial ideal of S.

*Example 3.17:* Let

$$S = \left\{ \sum_{i=0}^{\infty} [0, a_i] x^i \;\middle|\; a_i \in Z_{12} \right\}$$

be a interval polynomial semigroup under addition. Clearly S has only interval polynomial subsemigroups and has no ideals.

$$P = \left\{ \sum_{i=0}^{\infty} [0, a_i] x^i \;\middle|\; a_i \in \{0, 2, 4, 6, 8, 10\} \subseteq Z_{12} \right\} \subseteq S$$

is not an interval polynomial ideal of S.

*Example 3.18:* Let

$$S = \left\{ \sum_{i=0}^{2} [0, a_i] x^i \;\middle|\; a_i \in Z_{12}; x^3 = x^0 = 1, \; x^4 = x \text{ and so on} \right\}$$

be a polynomial interval semigroup. Find order of S. Is S a S-Lagrange interval polynomial semigroup?

*Example 3.19:* Let

$$S = \left\{ \sum_{i=0}^{3} [0, a_i] x^i \;\middle|\; a_i \in Z_2; x^4 = x^0 = 1, x^5 = x \right\}$$

be polynomial interval semigroup under multiplication.

Clearly $S = \{[0, 1]x, 0, [0, 1], [0, 1]x^2, [0, 1]x^3, [0, 1] + [0, 1]x [0, 1] + [0, 1]x^2, [0, 1] + [0, 1]x^3, [0, 1]x + [0, 1]x^2, [0, 1]x + [0, 1]x^3, [0, 1]x^2 + [0, 1] + [0, 1]x + [0, 1]x^3, [0, 1] + [0, 1]x^2 + [0, 1]x^3, [0, 1]x + [0, 1]x^2 + [0, 1]x^3, [0, 1] + [0, 1]x + [0, 1]x^2 + [0, 1]x^3\}$, and o (S) = 16.

$T = \{[0, 1]x, [0, 1]x^2, [0, 1]x^3, [0, 1]\} \subseteq S$ is a interval polynomial subgroup of S.

$P = \{[0, 1]x^2, [0, 1]\} \subseteq S$ is also a interval polynomial subgroup of S. Thus S is a S-interval polynomial semigroup. Infact S is a commutative interval polynomial semigroup with identity [0, 1]. Further $I = \{0, [0, 1] + [0, 1]x + [0, 1]x^2 + [0, 1]x^3\} \subseteq S$ is a interval polynomial ideal of S.

Now having seen examples of polynomial interval semigroups we now proceed onto define symmetric interval semigroups or permutation interval semigroup or interval permutative semigroup in the following chapter.

**Chapter Four**

# SPECIAL INTERVAL SYMMETRIC SEMIGROUPS

In this chapter we for the first time introduce the notion of mapping of n row intervals $([0, a_1], …, [0, a_n])$ to itself. This forms the semigroup under the composition of mappings and is isomorphic with the symmetric semigroup $S(n)$.

We also define special interval symmetric group and study some properties related with them.

**DEFINITION 4.1**: *Let $X = \{[0, a_1], [0, a_2], …, [0, a_n]\}$ be a set of n distinct intervals. We say $\eta : X \to X$ is an interval mapping if $\eta ([0, a_i]) = [0, a_j]; 1 \le i, j \le n$.*
    *Let $S(X)$ denote the collection of all interval mappings of X to X. $S(X)$ under the composition of interval mappings is a semigroup defined as the interval symmetric semigroup.*

We will first illustrate this situation by some examples.

***Example 4.1:*** Let $X = \{[0, a_1], ]0, a_2]\}$ be the interval set $a_1 \neq a_2$.
The set of all maps of X to X are as follows:
$\eta_1: X \to X$ given by
$\quad \eta_1([0, a_1]) = [0, a_1]$ and $\eta_1([0, a_2]) = [0, a_2]$.
$\eta_2: X \to X$ is given by
$\quad \eta_2([0, a_1]) = [0, a_2]$ and $\eta_2([0, a_2]) = [0, a_1]$,
$\eta_3: X \to X$ is defined by
$\quad \eta_3([0, a_1]) = [0, a_1]$ and $\eta_3([0, a_2]) = [0, a_1]$.
$\eta_4: X \to X$ is such that
$\quad \eta_4([0, a_1]) = [0, a_2]$ and $\eta_4([0, a_2]) = [0, a_2]$.

Thus $S(X) = \{\eta_1, \eta_2, \eta_3, \eta_4\}$; and $S(X)$ under composition of maps is an interval symmetric semigroup.
$\quad$ Clearly $|S(X)| = 2^2 = 4$.

***Example 4.2:*** Let $X = \{[0, a_1], [0, a_2], [0, a_3]\}$; $a_i \neq a_j$ if $i \neq j$ $a_i > 0$; $1 \leq i \leq 3$. The maps of X to X is $S(X) = \{\eta_1, \eta_2, \eta_3, \eta_4, \eta_5, \eta_6, \eta_7, \ldots, \eta_{26}, \eta_{27}\}$.

$\eta_1([0, a_i]) = [0, a_i]$; $i = 1, 2, 3$;
$\eta_2([0, a_1]) = [0, a_2]$, $\eta_2([0, a_2]) = [0, a_3]$,
$\eta_2([0, a_3]) = [0, a_1]$, $\eta_3([0, a_1]) = [0, a_1]$;
$\eta_3([0, a_2]) = [0, a_3]$, $\eta_3([0, a_3]) = [0, a_2]$;
$\eta_4([0, a_1]) = [0, a_2]$, $\eta_4([0, a_2]) = [0, a_1]$;
$\eta_4([0, a_3]) = [0, a_3]$, $\eta_5([0, a_1]) = [0, a_3]$;
$\eta_5([0, a_2]) = [0, a_2]$, $\eta_5([0, a_3]) = [0, a_1]$;
$\eta_6([0, a_1]) = [0, a_3]$, $\eta_6([0, a_2]) = [0, a_1]$;
$\eta_6([0, a_3]) = [0, a_2]$, $\eta_7([0, a_1]) = [0, a_1]$;
$\eta_7([0, a_2]) = [0, a_1]$, $\eta_7([0, a_3]) = [0, a_1]$, $\ldots$, $\eta_{27}([0, a_3]) = [0, a_3]$
$\eta_{27}([0, a_1]) = [0, a_3]$ and $\eta_{27}([0, a_2]) = [0, a_3]$.

Thus $o(S(X)) = 27 = 3^3$.
We see $S(X)$ is a interval symmetric group of mappings of the interval set X to itself. We can in general say if $X = \{[0, b_1], [0, b_2], \ldots, [0, b_n]\}$ with $b_i \neq b_j$; $i \neq j$, $(1 \leq i, j \leq n)$ then $S(X)$ is the interval symmetric semigroup of order $n^n$.

We have the following interesting theorem and observations when we say the interval [0, a] is mapped on to [0, b] we mean the continuous interval segment 0 to a is mapped onto the continuous interval segment 0 to b. We see the map may contract or extend the interval for instance [0, 5] is mapped to [0, $\sqrt{2}$ ] then certainly a contraction has taken place or we can realize the map is not an embedding.

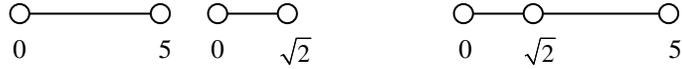

On the other hand if [0, $\sqrt{2}$ ] interval is mapped onto [0, 5] we can realize it as expansion. All these maps will be useful when we use the concept of finite element methods like, stiffness matrices or in any other applications.

However we see we have an isomorphism between S (n) and S(X) where X = {[0, $a_1$], [0, $a_2$], ..., [0, $a_n$]} all intervals are distinct $a_i > 0$ and $a_i \neq a_j$, if $i \neq j$; $1 \leq i, j \leq n$. After all $\eta_i$ : [0, $a_t$] $\to$ [0, $a_p$] are only maps. $1 \leq t, p \leq n$.

Keeping this in mind we have the following.

**THEOREM 4.1:** *Let S (n) be the symmetric semigroup on the set {1, 2, ..., n} and S (X) be the interval symmetric semigroup on the set X = {[0, $a_1$], [0, $a_2$], ..., [0, $a_n$]}, $a_i \neq a_j$, $i \neq j$, $a_i > 0$, $1 \leq i, j \leq n$. Then S (n) is isomorphic with S (X).*

*Proof* : We know the symmetric semigroup S(n) is of order $n^n$. Now the order of the interval symmetric semigroup S(X) where X = {[0, $a_1$], ..., [0, $a_n$]} is also of order $n^n$. Now if we put [0, $a_i$] = $x_i$, i = 1, 2, ..., n and S(n) is the set of all maps of (1, 2, ..., n) to itself and associate each i to $x_i$ ; $1 \leq i \leq n$, we see the one to one correspondence between the maps. Thus S(n) $\cong$ S(X).

We will illustrate this in case of S(2).

***Example 4.3:*** Let $S(X) = \{\eta_1, \eta_2, \eta_3, \eta_4\}$ where $\eta_i : \{1, 2\} \to \{1, 2\}$; $i = 1, 2, 3, 4$. $\eta_1(1) = 1$, $\eta_1(2) = 2$, $\eta_2(1) = 1$, $\eta_2(2) = 1$, $\eta_3(1) = 2$, $\eta_3(2) = 2$, $\eta_4(1) = 2$ and $\eta_4(2) = 1$ is the symmetric semigroup of order $2^2 = 4$.

Now $X = \{[0, a_1], [0, a_2]\}$ be the interval set $a_1 \ne a_2$ and $a_i > 0$; $i = 1, 2$. $S(X) = \sigma_1, \sigma_2, \sigma_3, \sigma_4\}$ where $\sigma_I : X \to X$; $i = 1, 2, 3, 4$.

$\sigma_1([0, a_1]) = [0, a_1]$, $\sigma_1 = ([0, a_2]) = [0, a_2]$
$\sigma_2([0, a_1]) = [0, a_1]$, $\sigma_2 = ([0, a_2]) = [0, a_1]$
$\sigma_3([0, a_1]) = [0, a_2]$, $\sigma_3 = ([0, a_2]) = [0, a_2]$
and $\sigma_4([0, a_1]) = [0, a_2]$, $\sigma_4 = ([0, a_2]) = [0, a_1]$.

$S(X)$ is the interval symmetric semigroup of order $2^2 = 4$. Now define a map $\mu : S(2) \to S(X)$ as follows.

$$\mu(\eta_1) = \sigma_1, \mu(\eta_2) = \sigma_2$$
$$\mu(\eta_3) = \sigma_3 \text{ and } \mu(\eta_4) = \sigma_4.$$

It is easily verified $\mu$ is a semigroup homomorphism, infact an isomorphism. Hence the claim.

We will enumerate some of the properties enjoyed by the interval symmetric semigroup.

***Example 4.4:*** Let $S(X)$ be the set of all maps from the three element interval set $X = \{[0, a_1], [0, a_2], [0, a_3]\}$ to itself. Clearly $S(X)$ is the semigroup under the operation of composition of map. Thus $S(X)$ is the symmetric interval semigroup of order $3^3 = 27$.

We see $S(X)$ is S-symmetric interval semigroup as it has 5 interval subgroups. For take $P_1 = \{\eta_1, \eta_2\}$ where

$$\eta_1 : \begin{cases} [0, a_1] \mapsto [0, a_1] \\ [0, a_2] \mapsto [0, a_2] \\ [0, a_3] \mapsto [0, a_3] \end{cases}$$

and

$$\eta_2 : \begin{cases} [0, a_1] \mapsto [0, a_1] \\ [0, a_2] \mapsto [0, a_3] \\ [0, a_3] \mapsto [0, a_2] \end{cases}$$

$P_2 = \{\eta_1, \eta_3\}$ where $\eta$ is given by and

$$\eta_3 : \begin{cases} [0, a_1] \mapsto [0, a_3] \\ [0, a_2] \mapsto [0, a_2] \\ [0, a_3] \mapsto [0, a_1] \end{cases}$$

$P_3 = \{\eta_1, \eta_4\}$ where $\eta_1$ is given above

$$\eta_4 : \begin{cases} [0, a_1] \mapsto [0, a_2] \\ [0, a_2] \mapsto [0, a_1] \\ [0, a_3] \mapsto [0, a_3] \end{cases}$$

$P_4 = \{\eta_1, \eta_5, \eta_6\}$ where $\eta_1$ is the identity map and

$$\eta_5 : \begin{cases} [0, a_1] \mapsto [0, a_2] \\ [0, a_2] \mapsto [0, a_3] \\ [0, a_3] \mapsto [0, a_1] \end{cases}$$

and

$$\eta_6 : \begin{cases} [0, a_1] \mapsto [0, a_3] \\ [0, a_2] \mapsto [0, a_1] \\ [0, a_3] \mapsto [0, a_2] \end{cases}$$

and $P_5 = \{\eta_1, \eta_2, \eta_3, \eta_4, \eta_5, \eta_6\}$.

Thus S(X) has $P_1$, $P_2$, $P_3$, $P_4$ and $P_5$ to be 5 interval subgroups of which $P_1$, $P_2$, $P_3$ and $P_4$ are cyclic. Thus S (X) is only a S-interval symmetric weakly cyclic semigroup as $P_5$ is not an abelian group.

In view of this we have the following theorems the proof of which are left as exercises for the reader.

**THEOREM 4.2:** *Let S (X) be a interval symmetric semigroup where $X = \{[0, a_1], [0, a_2], ..., [0, a_n]\}$, $a_i \neq a_j$; $i \neq j$, $1 < i, j < n$. S(X) is a S-weakly cyclic interval symmetric semigroup.*

**THEOREM 4.3:** *Let S (X) be a interval symmetric semigroup where $X = \{[0, a_1], [0, a_2], ..., [0, a_n]\}$, $a_i \neq a_j$, $i \neq j$, $1 < i, j < n$. S(X) is a S-interval symmetric semigroup.*

**THEOREM 4.4:** *Let S (X) be a interval symmetric semigroup where $X = \{[0, a_1], [0, a_2], ..., [0, a_n]\}$, $a_i \neq a_j$ if $i \neq j$, $a_i > 0$, $1 \leq i, j \leq n$. S(X) is only a Smarandache weakly commutative interval symmetric semigroup.*

*Proof*: Take in S (X), P the collection of all one to one mapping of X to itself, then P is a interval symmetric subgroup of S (X) but is not a commutative interval symmetric group.

Hence S (X) is only a S-weakly commutative interval symmetric semigroup.

We first proceed onto give the basic definition of interval symmetric group.

**DEFINITION 4.2:** *Let $X = \{[0, a_1], [0, a_2], ..., [0, a_n]\}$ be an interval set $S_X$ denote the set of all one to one maps of the interval set X. $S_X$ under the composition of mappings is a group, which will be known as the interval symmetric group.*

***Example 4.5:*** Let $S_X = \{\eta_1, \eta_2, \eta_3, ..., \eta_6\}$ where $X = \{[0, a_1], [0, a_2], [0, a_3]\}$ is the interval set. $S_X$ is the interval symmetric group with $\eta_1 ([0, a_i]) = [0, a_i]$; $1 \leq i \leq 3$. $\eta_2 ([0, a_1]) = [0, a_1]$, $\eta_2 ([0, a_2]) = [0, a_3]$ and $\eta_2 ([0, a_3]) = [0, a_2]$.

$\eta_3 ([0, a_1]) = [0, a_2]$, $\eta_3 ([0, a_2]) = [0, a_1]$
$\eta_3 ([0, a_3]) = [0, a_3]$, $\eta_4 ([0, a_1]) = [0, a_3]$
$\eta_4 ([0, a_2]) = [0, a_2]$, $\eta_4 ([0, a_3]) = [0, a_1]$
$\eta_5 ([0, a_1]) = [0, a_2]$, $\eta_5 ([0, a_2]) = [0, a_3]$
$\eta_5 ([0, a_3]) = [0, a_1]$, and $\eta_6 ([0, a_1]) = [0, a_3]$;
$\eta_6 ([0, a_2]) = [0, a_1]$ and $\eta_6 ([0, a_3]) = [0, a_2]$.

It is easily verified $S_X$ under the composition of maps is a group, called the interval symmetric group.

We have the following theorems which are left as exercises for the reader to prove.

**THEOREM 4.5:** *Let $S_X$ be the interval symmetric group on $X = \{[0, a_1], [0, a_2], ..., [0, a_n]\}$, of n distinct intervals $a_i > 0$, $S_n \cong S_X$ where $S_n$ is the symmetric group of degree n.*

**THEOREM 4.6:** *The S – interval symmetric semigroup $S(X)$ has its largest interval group $S_X$ to be contained in the proper interval subset $A = S_X \cup \{\sigma_1, \sigma_2, ..., \sigma_n\}$ where $\sigma_i([0, a_j]) = [0, a_i]$ for all $j = 1, 2, ..., n$. true for $i=1, 2, ..., n$, which is an interval symmetric subsemigroup of $S(X)$.*

We will illustrate this situation by an example.

*Example 4.6:* Let $X = \{[0, a_1], [0, a_2], [0, a_3], [0, a_4]\}$ be a interval set of cardinality four. $S(X)$ be the interval symmetric semigroup. Consider

$$A = S_X \cup \left\{ \begin{pmatrix} [0,a_1] & [0,a_2] & [0,a_3] & [0,a_4] \\ [0,a_2] & [0,a_2] & [0,a_2] & [0,a_2] \end{pmatrix}, \right.$$

$$\begin{pmatrix} [0,a_1] & [0,a_2] & [0,a_3] & [0,a_4] \\ [0,a_1] & [0,a_1] & [0,a_1] & [0,a_1] \end{pmatrix},$$

$$\begin{pmatrix} [0,a_1] & [0,a_2] & [0,a_3] & [0,a_4] \\ [0,a_3] & [0,a_3] & [0,a_3] & [0,a_3] \end{pmatrix},$$

$$\left. \begin{pmatrix} [0,a_1] & [0,a_2] & [0,a_3] & [0,a_4] \\ [0,a_4] & [0,a_4] & [0,a_4] & [0,a_4] \end{pmatrix} \right\}$$

Clearly A is a proper subset and is the interval symmetric subsemigroup of S(X). Further A is a S-hyper interval symmetric subsemigroup of S(X).

**COROLLARY 4.1:** *S(X) the S-interval symmetric semigroup is not a S-simple symmetric semigroup.*

We will be using the definitions of S-Lagrange semigroup and S-weakly Lagrange semigroup [9].

*Example 4.7:* Consider the interval symmetric semigroup S(X) where X = {[0, $a_1$], [0, $a_2$], [0, $a_3$], [0, $a_4$]} of 4 distinct intervals $a_i > 0$; i = 1, 2, 3, 4. Clearly order of S(X) is $4^4$. We see S(X) has $S_X$ to be interval subgroup of order 24. Clearly 24 $\nmid$ $4^4$. Thus S(X) is not a S-Lagrange interval symmetric semigroup. However S(X) has interval subgroups of order two and four which divides $4^4$. Hence S(X) is a S-interval symmetric weakly Lagrange semigroup.

**THEOREM 4.7:** *Let S(X) be a interval symmetric semigroup on n-distinct intervals, i.e. X = {[0, $a_1$], [0, $a_2$], …, [0, $a_n$]},*
    *(1) S(X) is a S-interval symmetric semigroup.*
    *(2) S(X) is not a S-interval symmetric Lagrange semigroup.*
    *(3) S(X) is a S-interval symmetric weakly Lagrange semigroup.*
    *(4) S(X) has Smarandache interval symmetric p-Sylow subgroups provided X has p number of distinct intervals and p is a prime. If the number of distinct intervals in X is a composite number say n and if p is a prime such that p/n then also S(X) has S-p-Sylow subgroups.*

**THEOREM 4.8:** *Let S (X) be a interval symmetric semigroup. S(X) has S-p-Sylow semigroups.*

Proof is direct and is left as an exercise to the reader.

Please refer [9] for Smarandache Cauchy elements of a semigroup.

*Example 4.8:* Let S(X) be a interval symmetric semigroup with X = {[0, $a_1$], [0, $a_2$], [0, $a_3$], [0, $a_4$], [0, $a_5$]} where X has 5-distinct intervals, o(S(X)) = $5^5$. We have $\eta \in$ S(X) such that $\eta^5$ = identity map and 5 | $5^5$. Thus S(X) has Smarandache Cauchy elements.

In view of this we have the following results.

**THEOREM 4.9:** *Let S(X) be a interval symmetric semigroup of order $n^n$ where X = {[0, $a_1$], [0, $a_2$], …, [0, $a_n$]}. S(X) has S-Cauchy elements.*

The proof is left as an exercise for the reader.

However we will illustrate this situation by an example.

*Example 4.9:* Let S(X) be a interval symmetric semigroup of order $6^6$. S(X) has S-Cauchy elements.
For take

$$\eta_t = \begin{pmatrix} [0,a_1] & [0,a_2] & [0,a_3] & [0,a_4] & [0,a_5] & [0,a_6] \\ [0,a_2] & [0,a_3] & [0,a_4] & [0,a_5] & [0,a_6] & [0,a_1] \end{pmatrix}$$

in S(X). Clearly $(\eta_t)^6$ = identity element of S(X). Thus S(X) has S-Cauchy elements.

$$\eta_p = \begin{pmatrix} [0,a_1] & [0,a_2] & [0,a_3] & [0,a_4] & [0,a_5] & [0,a_6] \\ [0,a_2] & [0,a_3] & [0,a_1] & [0,a_4] & [0,a_5] & [0,a_6] \end{pmatrix}$$

in S(X) is such that $(\eta_p)^3$ = identity elements of S(X). Thus S(X) has S-Cauchy elements.

In view of this we have the following theorem.

**THEOREM 4.10:** *Let S(X) be a interval symmetric semigroup of order $n^n$; n a composite number S(X) has S-Cauchy elements.*

*Proof:* Take

$$\eta_i = \begin{pmatrix} [0,a_1] & [0,a_2] & \cdots & [0,a_{n-1}] & [0,a_n] \\ [0,a_2] & [0,a_3] & \cdots & [0,a_n] & [0,a_1] \end{pmatrix}$$

in S (X). Clearly $(\eta_i)^n$ = identity element of S (X). Take p | n. p a prime then

$$\eta_t = \begin{pmatrix} [0,a_1] & [0,a_2] & \cdots & [0,a_{p-1}] & [0,a_p] & [0,a_{p+1}] & \cdots & [0,a_n] \\ [0,a_2] & [0,a_3] & \cdots & [0,a_p] & [0,a_1] & [0,a_{p+1}] & \cdots & [0,a_n] \end{pmatrix}$$

in S(X) is such that $(\eta_t)^p$ = identity element of S(X). Thus S(X) has S-Cauchy elements.

Cayley's theorem for S-semigroups can also be extended in case of S-interval semigroups.

Please refer for S-semigroup homomorphism and S-semigroup automorphism [9]. Since $S(X) \cong S(n)$ where S (n) is a permutation of (1, 2, …, n) and S(X) is a interval symmetric semigroup with n interval set $X = \{[0, a_1], [0, a_2], …, [0, a_n]\}$, $a_i > 0$ and $a_i \neq a_j$ if $i \neq j$; $1 \leq i, j \leq n$.

We can have Cayley's theorem for S-interval semigroups. Several interesting results in classical group theory can be proved for S-interval semigroups with appropriate modifications.

Now we proceed onto define special interval symmetric semigroups.

For now on wards by a special interval [a, b] we mean a < b, a ≠ 0, a and b positive integers.

We say $X = \{[a_1, b_1], [a_2, b_2], …, [a_n, b_n]\}$ is a special interval set if $a_i < b_i$, $a_i \neq 0$; $1 \leq i \leq n$ and each $a_j < b_k$ for every $1 \leq j, k \leq n$; that is all $a_i$'s are less than $b_k$ even if $i \neq k$. We call such interval collection to be special interval collection.

We will first illustrate special intervals.

***Example 4.10:*** Let $X = \{[7, 12], [5, 10], [3, 8]\}$, X is a special interval set. Suppose $X = \{[5, 10], [6, 7], [54, 5], [9, 12]\}$; X is not a special interval set as $7 < 9$ and $6 < 5$ so $a_i < b_j$ is not true for every $a_i$ and $b_j$.

Let $X = \{[a_1, b_1], [a_2, b_2]\}$ be the special interval set then the interval set generated by X denoted by $\langle X \rangle$ is $\{[a_1, b_1], [a_2, b_2], [a_2, b_1], [a_2, b_2]\}$ is an interval set.

Let $X = \{[a_1, b_1], [a_2, b_2], [a_3, b_3]\}$ be the special interval set, then the interval set generated by X denoted by $\langle X \rangle$ is $\{\{[a_1, b_1], [a_2, b_2], [a_3, b_3], [a_1, b_2], [a_1, b_3], [a_2, b_1], [a_2, b_3] [a_3, b_1], [a_3, b_2]\}$.

Thus we see if $X = \{[a_1, b_1], [a_2, b_2], \ldots, [a_n, b_n]\}$ is a special interval set then the interval set generated by X denoted by $\langle X \rangle$ is $\{[a_1, b_1], \ldots, [a_n, b_n], [a_1, b_2], [a_1, b_3], \ldots, [a_1, b_n], \ldots, [a_n, b_1], [a_n, b_2], \ldots, [a_n, b_{n-1}]\}$. Clearly the number of elements in $\langle X \rangle$ is $n^2$.

Now we proceed onto define the notion of special interval symmetric semigroup or interval special symmetric semigroup.

**DEFINITION 4.3:** *Let X be a special interval set $\langle X \rangle$ the interval set generated by X. $S(\langle X \rangle)$ set of all mappings of $\langle X \rangle$ to $\langle X \rangle$. $S(\langle X \rangle)$ is defined as the special interval symmetric semigroup or interval special symmetric semigroup.*

*The order of $S(\langle X \rangle) = \left(n^2\right)^{n^2}$.*

We will illustrate this situation by an example.

***Example 4.11:*** Let $X = \{[a_1, b_1], [a_2, b_2]\}$ be a special interval set with two distinct elements. $\langle X \rangle = \{[a_1, b_1], [a_2, b_2], [a_1, b_2], [a_2, b_1]\}$. Now the set of all maps from $\langle X \rangle$ to $\langle X \rangle$ denoted by $S(\langle X \rangle)$. We have $o(S(\langle X \rangle)) = 4^4$. $S(\langle X \rangle)$ is a special symmetric interval semigroup.

These structures will be useful in several applications.

We can as in case of interval semigroups build the concept of special interval groups also.

$S_{\langle X \rangle}$ = {all one to one maps of $\langle X \rangle$ to itself} is a special interval symmetric group. Clearly $S_X \subseteq S_{\langle X \rangle}$; likewise $S(X) \subseteq S(\langle X \rangle)$.

$S_X$ is a proper subgroup of $S_{\langle X \rangle}$ and $|S_X| = \lfloor |X| $ and $|S_{\langle X \rangle}| = \lfloor \langle X \rangle$. Further $|S(X)| = |X|^{|X|}$ and $|S(\langle X \rangle)| = |\langle X \rangle|^{|\langle X \rangle|}$.

The following results are important.

**THEOREM 4.11:** *Let $X = \{[a_1, b_1], [a_2, b_2], \ldots, [a_n, b_n]\}$ be a special interval set. $\langle X \rangle$ denotes the intervals generated by X. Let $S(\langle X \rangle)$ be the special symmetric interval semigroup; $S(\langle X \rangle)$ is a S-special interval symmetric semigroup.*

*Proof*: We know $S(X)$ is a interval symmetric semigroup. $S_X \subseteq S(X)$ so $S(X)$ is a S-interval symmetric semigroup as $S_X$ is the interval symmetric group. Similarly $S_{\langle X \rangle}$ is the special symmetric interval group of $S(\langle X \rangle)$ as the set of all maps of $\langle X \rangle$ to $\langle X \rangle$ and $S_{\langle X \rangle}$ is the set of all one to one maps of $\langle X \rangle$ into $\langle X \rangle$ $S(\langle X \rangle)$ is a S-special symmetric interval semigroup.

We have several interesting results associated with these special symmetric interval semigroups which we describe below.

**THEOREM 4.12:** *Let $S(\langle X \rangle)$ be a special symmetric interval semigroup $P \cong S(X)$ is a S-special symmetric interval subsemigroup of $S(\langle X \rangle)$ and is not a S-special symmetric interval hyper subsemigroup.*

*Proof:* We know $S(\langle X \rangle)$ to be a special symmetric interval semigroup which contains as isomorphic copy of $S(X)$ as a special symmetric interval subsemigroup.
For if $\{[a_1, b_1], [a_2, b_2], \ldots, [a_n, b_n]\} = X$ then $\langle X \rangle = \{[a_1, b_1], [a_2, b_2], [a_3, b_3], \ldots, [a_n, b_n], [a_2, b_1], [a_2, b_3], \ldots, [a_2, b_n], [a_3, b_1], \ldots, [a_n, b_{n-1}]\}$. $P$ = {all maps from $\langle X \rangle$ to $\langle X \rangle$ which fixes the intervals $\{[a_1, b_2], [a_1, b_3], \ldots, [a_n, b_{n-2}] [a_n, b_{n-1}]\}$. We say a map fixes intervals if the intervals are mapped to itself.

Clearly $P \cong S(X)$ and P is a special symmetric interval subsemigroup of $S(\langle X \rangle)$. Now if we take all one to one maps of $T = \{[a_1, b_1], [a_2, b_2], \ldots, [a_n, b_n], [a_1, b_2], \ldots, [an, b_{n-1}]\}$ fixing the interval elements $[a_1, b_2], \ldots, [a_n, b_{n-1}]$. We see T is a interval symmetric group and is isomorphic with the symmetric group $S_n$ as $|X| = n$. Thus $T \cong S_X \cong S_n$. Hence T is a interval symmetric group isomorphic to $S_X$. Thus P is a S-special symmetric interval subsemigroup of $S(\langle X \rangle)$. Now if P is to be a special symmetric interval hyper subsemigroup then P must contain the largest special interval symmetric group of $S(\langle X \rangle)$. We see the largest special interval symmetric group B of $S(\langle X \rangle)$ is the set of all one to one maps of $\langle X \rangle$ to $\langle X \rangle$. Thus clearly $B \not\subseteq P$ so P is not a special symmetric interval hyper subsemigroup of $S\langle X \rangle$. Hence the claim.

**THEOREM 4.13:** *The special interval symmetric semigroup $S(\langle X \rangle)$ is not a S-Lagrange special interval symmetric semigroup but is a S weakly Lagrange special interval symmetric semigroup.*

Proof is left as an exercise for the reader.

We will however illustrate this situation by an example.

***Example 4.12:*** Let $X = \{[a_1, b_1], [a_2, b_2], [a_3, b_3]\}$ be the interval set with 3 distinct elements. $\langle X \rangle = \{[a_1, b_1], [a_2, b_2], [a_3, b_3], [a_1, b_2], [a_1, b_3], [a_2, b_1], [a_2, b_3], [a_3, b_1], [a_3, b_2]\}$ be the special interval symmetric semigroup of X. We see $o(S(\langle X \rangle)) = 9^9$. We see $S(\langle X \rangle)$ has interval symmetric subgroup of order 6. But $6 \nmid 9^9$, so $S(\langle X \rangle)$ is not a S-Lagrange special interval semigroup. $S(\langle X \rangle)$ has some special interval subgroups which divides the order of $S(\langle X \rangle)$. For take the special interval group P generated by $\eta$ where $\eta =$

$\begin{pmatrix} [a_1,b_1] & [a_2,b_2] & [a_3,b_3] & [a_1,b_2] & [a_1,b_3] & [a_2,b_1] & [a_3,b_3] & [a_3,b_1] & [a_3,b_2] \\ [a_2,b_2] & [a_3,b_3] & [a_1,b_2] & [a_1,b_3] & [a_2,b_3] & [a_2,b_3] & [a_3,b_1] & [a_3,b_2] & [a_1,b_1] \end{pmatrix}$

clearly number of elements of P is 9 and $9 \mid (9)^9$. Hence $S(\langle X \rangle)$ is only a S-weakly Lagrange symmetric interval semigroup.

The following theorems are expected to be proved by the reader.

**THEOREM 4.14:** *Let $S(\langle X \rangle)$ be the special interval symmetric semigroup $S(\langle X \rangle)$ has S-p-Sylow subgroups.*

**THEOREM 4.15:** *Let $S(\langle X \rangle)$ be the special interval symmetric semigroup. $S(\langle X \rangle)$ has S-Cauchy elements.*

**THEOREM 4.16:** *Let $S(X)$ be the S-symmetric interval semigroup where the interval set X has p distinct elements; p a prime. Then $S(X)$ has S-non p-Sylow subgroup for all primes $t < p$.*

We will illustrate this situation by an example.

***Example 4.13:*** Let $X = \{[a_1, b_1], [a_2, b_2], [a_3, b_3], [a_4, b_4], [a_5, b_5]\}$ be the interval set with five distinct intervals. $S(X)$ is the interval symmetric semigroup. Clearly $o(S(X)) = 5^5$. Now $3 < 5$, 3 is a prime. $S(X)$ has a subgroup P of order three generated by

$$\eta = \begin{pmatrix} [a_1,b_1] & [a_2,b_2] & [a_3,b_3] & [a_4,b_4] & [a_5,b_5] \\ [a_2,b_2] & [a_3,b_3] & [a_1,b_1] & [a_4,b_4] & [a_5,b_5] \end{pmatrix}.$$

P = {identity element of $S(X)$, $\eta$, $\eta^2$, as $\eta^3$ = identity element of $S(X)$}. Thus $o(P) = 3$ and $3 \nmid 5^5$.

Consider

$$\sigma = \begin{pmatrix} [a_1,b_1] & [a_2,b_2] & [a_3,b_3] & [a_4,b_4] & [a_5,b_5] \\ [a_2,b_2] & [a_1,b_1] & [a_3,b_3] & [a_4,b_4] & [a_5,b_5] \end{pmatrix}$$

in $S(X)$. $\sigma$ generates a group T of order 2 as $\sigma^2$ = identity element of $S(X)$. Now $o(T) = 2$ and $2 \nmid 5^5$. Only subgroups of order 5 are S-5-Sylow subgroups of $S(X)$. Clearly $S(X)$ has no S-3-Sylow subgroups and S-2-Sylow subgroups.

**Chapter Five**

# NEUTROSOPHIC INTERVAL SEMIGROUPS

In this chapter we for the first time introduce the notion of neutrosophic interval semigroups and pure neutrosophic interval semigroups and discuss several interesting properties related with them.

We just briefly introduce the notations which we will be using in this chapter and as well as in the following chapters.

$Z_nI = \{tI \mid t \in Z_n = \{0, 1, 2, \ldots, n-1\}$, where I is the indeterminancy introduced in chapter one of this book. Now $PN(Z_nI) = \{[0, tI] \mid tI \in Z_nI\}$ denotes the pure neutrosophic modulo integer intervals. Likewise $PN(Z^+I) = \{[0, nI] \mid nI \in Z^+I \cup \{0\}\}$ is the pure neutrosophic integer intervals.

$PN(Q^+I) = \{[0, tI] \mid tI \in Q^+I \cup \{0\}\}$ is the pure neutrosophic rational intervals.

$PN(R^+I) = \{[0, tI] \mid tI \in R^+I \cup \{0\}\}$ is the pure neutrosophic real intervals.

PN ($C^+I$) = {[0, tI] | tI ∈ $C^+I \cup \{0\}$} is the pure neutrosophic complex intervals.

N ($Z_n$) = {[0, a] | a = x + yI; x, y ∈ $Z_n$ and I the interminary} is the neutrosophic modulo integers intervals,

N($Z^+ \cup \{0\}$) = {[0, a] | a = x + yI where x, y ∈ $Z^+ \cup \{0\}$} is the neutrosophic integer intervals;

N ($Q^+ \cup \{0\}$) = {[0, a] | a = x + yI where x, y ∈ $Q^+ \cup \{0\}$} is the neutrosophic rational intervals;

N ($R^+ \cup \{0\}$) = {[0, a] | a = x + yI where x, y ∈ $R^+ \cup \{0\}$} is the neutrosophic real intervals and

N ($C^+ \cup \{0\}$) = {[0, a] | a = x + yI where x, y ∈ $C^+ \cup \{0\}$} is the neutrosophic complex intervals.

It is to be noted that even if the mention of which type of interval is not made it can be easily understood by the context.

Now we proceed onto define the notion of neutrosophic interval semigroup.

**DEFINITION 5.1:** *Let S = {[0, tI] | tI ∈ $Z_nI$}; S under addition modulo n is a semigroup called the interval pure neutrosophic semigroup or pure neutrosophic interval semigroup.*

We can replace $Z_nI$ by $Z^+I$ or $Q^+I$ or $R^+I$ or $C^+I$.
We will illustrate this by examples.

*Example 5.1:* Let S = {[0, aI] | aI ∈ $Z_5I$} be a pure neutrosophic interval modulo integer semigroup under addition.

*Example 5.2:* Let S = {[0, aI] | aI ∈ $Z^+I \cup \{0\}$} be a pure neutrosophic interval integer semigroup under addition.
Clearly S is infinite order in example 5.2 and S in example 5.1 is of finite order.

*Example 5.3:* Let S = {[0, aI] | aI ∈ $Q^+I \cup \{0\}$} be the pure neutrosophic interval rational semigroup under addition of infinite order.

*Example 5.4:* Let S = {[0, aI] | aI ∈ $R^+I \cup \{0\}$} be the pure neutrosophic interval real semigroup under addition of infinite order.

*Example 5.5:* Let S = {[0, aI] | aI ∈ $C^+I \cup \{0\}$} be the pure neutrosophic complex interval semigroup under addition of infinite order.

Now having seen neutrosophic interval semigroups under addition of finite and infinite order we can define pure neutrosophic interval semigroup under multiplication by replacing addition by multiplication of intervals in a very simple way [0, aI], [0, bI] = [0, abI] and [0, aI] + [0, bI] = [0, (a + b)I]. In the case of multiplication we do not include $C^+I \cup \{0\}$ as $i^2 =$ –1 and [0, i] [0, i] = [0, -1], not defined.

We give only examples of pure neutrosophic interval semigroup under multiplication.

*Example 5.6:* Let P = {[0, aI] | aI ∈ $Z_8I$} be a pure neutrosophic interval modulo integer semigroup under multiplication.

*Example 5.7:* Let W = {[0, aI] | aI ∈ $Q^+I \cup \{0\}$} be a pure neutrosophic rational interval semigroup under multiplication. Clearly W is of infinite order.

*Example 5.8:* Let T = {[0, aI] | aI ∈ $Z^+I \cup \{0\}$} be a pure neutrosophic integer interval semigroup under multiplication. T is also of infinite order.

*Example 5.9:* Let B = {[0, aI] | aI ∈ $R^+I \cup \{0\}$} be a pure neutrosophic integer real interval semigroup under multiplication.
However we cannot define multiplication of pure neutrosophic complex interval semigroup as it is not compatible.
Now we can define pure neutrosophic interval subsemigroup, pure neutrosophic interval ideal and S-pure neutrosophic ideal and S-pure neutrosophic interval semigroup.

We will illustrate these situations before we define the notion of zero divisors, nilpotents and idempotents in pure neutrosophic interval semigroups.

*Example 5.10:* Let $S = \{[0, aI] \mid aI \in Z_9I\}$ be a pure neutrosophic modulo integer interval semigroup under multiplication. Take $P = \{[0, aI] \mid aI \in \{0, 3I, 6I\}\} \subseteq S$, P is a pure neutrosophic modulo integer interval ideal of S.

*Example 5.11:* Let $M = \{[0, aI] \mid aI \in Z^+I \cup \{0\}\}$ be a pure neutrosophic interval semigroup, clearly $P = \{[0, aI] \mid aI \in 5Z^+I \cup \{0\}\} \subseteq M$ is a pure neutrosophic interval subsemigroup of M. P is also a pure neutrosophic interval ideal of M.

*Example 5.12:* Let $S = \{[0, aI] \mid aI \in Z_{12}I\}$ be a pure neutrosophic interval semigroup under addition modulo 12. $P = \{[0, aI] \mid aI \in \{0, 2I, 4I, 6I, 8I, 10I\} \subseteq Z_{12}I\} \subseteq S$ is a pure neutrosophic interval subsemigroup of S. Clearly S is not a pure neutrosophic interval ideal of S.

*Example 5.13:* Let $T = \{[0, aI] \mid aI \in Z^+I \cup \{0\}\}$ be a pure neutrosophic interval semigroup under addition. Clearly T has no ideals only pure neutrosophic interval subsemigroups.

*Example 5.14:* Let $W = \{[0, aI] \mid aI \in Q^+I \cup \{0\}\}$ be a pure neutrosophic interval semigroup under addition (or multiplication). W has only pure neutrosophic interval subsemigroups and has no ideals.

*Example 5.15:* Let $N = \{[0, aI] \mid aI \in R^+I \cup \{0\}\}$ be a pure neutrosophic interval semigroup under addition (or multiplication). N has no ideals, N has only pure neutrosophic interval subsemigroups given by $P = \{[0, aI] \mid aI \in Z^+I \cup \{0\}\} \subseteq N$ or $T = \{[0, aI] \mid aI \in Q^+I \cup \{0\}\} \subseteq N$ and so on. However P and T are not ideals of N under addition or multiplication.

Every element in $S = \{[0, aI] \mid aI \in Z^+I \cup \{0\}\}$ is torsion free with respect addition (or multiplication). Likewise $W = \{[0, aI] \mid$ $aI \in Q^+I \cup \{0\}\}$ or $N = \{[0, aI] \mid aI \in R^+I \cup \{0\}\}$ under addition (or multiplication) is torsion free. Thus S, N or W does not contain any element of finite order.

Before we proceed to define neutrosophic interval semigroups which are not pure we discuss about zero divisors, idempotents and nilpotents in these structures. The pure neutrosophic interval semigroup built using $Z^+I \cup \{0\}$ or $Q^+I \cup \{0\}$ or $R^+I \cup \{0\}$ under addition or multiplication has no interval idempotents or interval zero divisors or interval nilpotents. However all pure neutrosophic interval semigroups built using $Z_nI$ have zero divisors, nilpotents and idempotents (n a composite number).

*Example 5.16:* Let $S = [0, aI] \mid aI \in Z_{15}I\}$ be a pure neutrosophic interval semigroup under multiplication. Consider $x = [0, 6I]$ in S. Clearly $[0, 6I][0, 6I] = [0, 36I] = [0, 6I]$. Thus S has non trivial interval idempotent. Also $y = [0, 10I]$ in S is such that $y^2 = [0, 10I][0, 10I] = [0, 100I] = [0, 10I]$, thus $[0, 10I]$ is an interval idempotent.

Also $t = [0, 11t]$ in S is such that $t^2 = [0, 11t][0, 11t] = [0, 121I] = [0, I]$ is a torsion element in S as $[0, I]$ in S acts as the unit. $w = [0, 4I]$ in S is such that $w^2 = [0, 4I][0, 4I] = [0, 16I] = [0, I]$ is a torsion element and is not a S-Cauchy element in S, as $2 \times 15$.

We see $[0, 3I].[0, 5I] = [0, 15I] = 0$ is the interval zero divisor in S. Also $[0, 5I][0, 6I] = [0, 30I] = 0$ is again a interval zero divisor in S. Further $[0, 3I][0, 10I] = [0, 30I] = 0$ is again a interval zero divisor in S.

Now this pure neutrosophic interval semigroup S has interval units, interval zero divisors and interval idempotents.

*Example 5.17:* Let $S = \{[0, aI] \mid aI \in Z_7I\}$ be a pure neutrosophic interval semigroup under multiplication S has no idempotents or zero divisors but has only units.

*Example 5.18:* Let $S = \{[0, aI] \mid aI \in Z_{30}I\}$ be a pure neutrosophic interval semigroup under multiplication. Clearly $o(S) = 30$. $[0, I]$ is the identity element of S.

S is a S-pure neutrosophic interval semigroup. We have elements which are units or self inversed. For $[0, 29I]$, $[0, 29I] = [0, I]$ where $[0, 29I] \in S$.

$[0, 6I] \in S$ is an idempotent of S as $[0, 6I] [0, 6I] = [0, 36I] = [0, 6I]$. We have $[0, 6I] [0, 10I] = [0, 0] = 0$ and $[0, 3I] [0, 10I] = [0, 0] = 0$ to be zero divisors. Several other zero divisors are present $[0, 2I], [0, 15I] = 0$ and so on.

Also $[0, 10I] [0, 10I] = [0, 100I] = [0, 10I]$ is an idempotent of S.

$[0, 16I] [0, 16I] = [0, 256I] = [0, 16I]$ is also an idempotent of S. This S has no nilpotents.

Now we proceed onto give examples of pure neutrosophic interval semigroups which has non trivial nilpotents.

*Example 5.19:* Let $S = \{[0, aI] \mid aI \in Z_{27}I\}$ be a pure neutrosophic interval semigroup. Consider $[0, 3I] \in S$. We see $[0, 3I]^3 = [0, 27I] = 0$ thus $[0, 3I]$ is a nontrivial nilpotent element of S.

Also $[0, 9I]^2 = [0, 0]$ is again a nontrivial nilpotent element of S. Further $[0, 18I] \in S$ is nilpotent of order two that is $[0, 18I]^2 = 0$. But $[0, 26I]$ in S is such that $[0, 26I]^2 = [0, I]$ is a unit in S. Further $[0, 4I] [0, 7I] = [0, 28I] = [0, I]$, $[0, 4I]$ is the inverse of $[0, 7I]$. Likewise $[0, 2I]$ is the inverse of $[0, 14I]$ as $[0, 2I] [0, 14I] = [0, I]$, $[0, 11I] [0, 5I] = [0, I]$ is again a unit in S.

Now we will proceed onto describe other properties about the pure neutrosophic interval semigroups built using $Z_nI$.

*Example 5.20:* Let $S = \{[0, aI] \mid aI \in Z_{12}I\}$ be a pure neutrosophic interval semigroup, S is a S-semigroup and S is also a S-Lagrange pure neutrosophic interval semigroup.

The pure neutrosophic interval subgroups of S are $P_1 = \{[0, I], [0, 7I]\}$, $P_2 = \{[0, I], [0, 11I]\}$, $P_3 = \{[0, 4I], [0, 8I]\}$, $P_4 =$

{[0, 3I], [0, 9I]}; $P_5 = \{0, I], [0, 5I]\}$ and $P_6 = \{[0, I], [0, 5I], [0, 7I], [0, 11I]\}$.

*Example 5.21:* Let $S = \{[0, aI] \mid aI \in Z_{10}I\}$ be a pure neutrosophic interval semigroup. It is easily proved S is a S-weakly Lagrange pure neutrosophic interval semigroup.

*Example 5.22:* Let $S = \{[0, aI] \mid aI \in Z_7I\}$ be a pure neutrosophic interval semigroup. It is easily verified S is S-simple pure neutrosophic interval semigroup.

*Example 5.23:* Let $P = \{[0, aI] \mid aI \in Z_4I\}$ be a pure neutrosophic interval semigroup. P is a S-Lagrange pure neutrosophic interval semigroup.

*Example 5.24:* Let $T = \{[0, aI] \mid aI \in Z_9I\}$ be a pure neutrosophic interval semigroup. T is not a S-Lagrange pure neutrosophic interval semigroup.

*Example 5.25:* Let $W = \{[0, aI] \mid aI \in Z_{11}I\}$ be a pure neutrosophic interval semigroup. W does not contain a S-pure neutrosophic interval hyper subsemigroup.

*Example 5.26:* Let $T = \{[0, aI] \mid aI \in Z_{23}I\}$ be a pure neutrosophic interval semigroup. $A = \{[0, I], [0, 22I]\}$ is a subgroup of order 2 in T. Thus T has S-non 2-Sylow pure neutrosophic interval subgroup.

**THEOREM 5.1:** *Let $S = \{[0, aI] \mid aI \in Z_p$; p a prime} be a pure neutrosophic interval semigroup, S has no pure neutrosophic interval hyper subsemigroup.*

The proof is left as an exercise for the reader.

*Example 5.27:* Let $T = \{[0, aI] \mid aI \in Z_{19}I\}$ be a pure neutrosophic interval semigroup. It is easily verified T has no S-Cauchy element.
    In view of this we have the following theorem the proof of which is left an as exercise for the reader.

**THEOREM 5.2:** *Let W = {[0, aI] | aI ∈ $Z_pI$; p a prime} be a pure neutrosophic interval semigroup. W has no S-Cauchy element.*

**THEOREM 5.3:** *Let W = {[0, aI] | aI ∈ $Z_pI$; p a prime} be a pure neutrosophic interval semigroup. Clearly W has a proper subgroup which cannot be properly contained in a proper pure neutrosophic interval subsemigroup.*

*Proof:* Given W = {[0, [0, I], [0, 2I], [0, 3I], …, [0, (p-1)I]} is a S-pure neutrosophic interval semigroup of order p, p a prime. The pure neutrosophic interval subgroups of W are $A_1$ = {[0, I] [0, (p-1)I]} and $A_2$ = {[0, I] [0, 2I], …, [0, (p-1)I]}. Clearly $A_1$ ∪ {0} is a pure neutrosophic interval subsemigroup of $Z_p$. Hence $A_1$ ∪ {0} is a S-pure neutrosophic interval subsemigroup of W. But $A_2$ cannot be strictly contained in any proper subset of W. So $A_2$ cannot be contained in a pure neutrosophic interval subsemigroup of W. Hence the claim.

It is left for the reader to prove the following theorem.

**THEOREM 5.4:** *Let S = {[0, aI] | aI ∈ $Z_n$} be a S-Lagrange pure neutrosophic interval semigroup then S is a S-weakly Lagrange pure neutrosophic interval semigroup.*

*Proof*: Clear from the very definition.

Now having seen properties about pure neutrosophic interval semigroup we now proceed onto define neutrosophic interval semigroup.

**DEFINITION 5.2:** *Let S = {[0, a + bI] | a, b ∈ $Z_n$}, under the operations of addition S is a semigroup called the neutrosophic interval semigroup.*

*We can replace $Z_n$ by $Q^+$ ∪ {0} or $R^+$ ∪ {0} or $Z^+$ ∪ {0} or $C^+$ ∪ {0}.*

*Example 5.28:* Let S = {[0, a + bI] | a, b ∈ $Z_9$} be a neutrosophic interval semigroup under addition. Clearly S is of finite order.

*Example 5.29:* Let $T = \{[0, a + bI] \mid a, b \in Q^+ \cup \{0\}\}$ be a neutrosophic interval semigroup under addition of infinite order.

*Example 5.30:* Let $W = \{[0, a + bI] \mid a, b \in Z^+ \cup \{0\}\}$ be a neutrosophic interval semigroup under addition. W is also an infinite semigroup.

It is important and interesting to note that the pure neutrosophic interval semigroups are always a proper subsemigroups of a neutrosophic interval semigroup.

*Example 5.31:* Let $S = \{[0, a + bI] \mid a, b \in Z_{10}\}$ be a neutrosophic interval semigroup. $T = \{[0, bI] \mid b \in Z_{10}\} \subseteq S$; is a neutrosophic interval subsemigroup. We see T is a pure neutrosophic interval semigroup, which is a subsemigroup of S.

*Example 5.32:* Let $W = \{[0, a + bI] \mid a, b \in Q^+ \cup \{0\}\}$ be a neutrosophic interval semigroup. Take $S = \{[0, aI] \mid a \in Q^+ \cup \{0\}\} \subseteq W$, S is a neutrosophic interval subsemigroup of W.

We can as in case of pure neutrosophic interval semigroup define interval subsemigroups and ideals of neutrosophic interval semigroups, which is left as an exercise for the reader.

*Example 5.33:* Let $S = \{[0, a + bI] \mid a, b \in Z_{20}\}$ be a neutrosophic interval semigroup. Now $T = \{[0, a + bI] \mid a, b \in \{0, 2, 4, 6, 8, 10, 12, 14, 16, 18\} \subseteq Z_{20}\} \subseteq S$ is a proper neutrosophic interval subsemigroup of S. $C = \{[0, a + bI] \mid a, b \in \{0, 5, 10, 15\} \subseteq Z_{20}\} \subseteq S$ is a proper neutrosophic interval subsemigroup of S.

*Example 5.34:* Let $T = \{[0, a + bI] \mid a, b \in Z^+ \cup \{0\}\}$ be a neutrosophic interval semigroup. Clearly T is not a S-neutrosophic interval semigroup. T has infinite number of neutrosophic interval subsemigroups given by $P_n = \{[0, a + bI] \mid a, b \in nZ^+ \cup \{0\}\} \subseteq T$; $n = 1, 2, \ldots, \infty$.

***Example 5.35:*** Let $S = \{[0, aI] \,/\, aI \in Z_{19}I\}$ be a S-pure neutrosophic interval semigroup. Clearly S has no S - pure neutrosophic interval hyper sub semigroup. $W = \{[0, a + bI] \,/\, a, b \in Z_{19}\}$, the neutrosophic interval semigroup contains S as a neutrosophic interval subsemigroup. Clearly S is a S-neutrosophic interval semigroup.

The properties of S-coset in case of pure neutrosophic interval semigroup and neutrosophic interval semigroup can be studied by the reader. For the notion of S-coset refer [9].

It is interesting to note that every neutrosophic interval semigroup contains an interval semigroup which is not neutrosophic.

In view of this we make the following definition.

**DEFINITION 5.3:** *Let $S = \{[0, a+bI] \,/\, a, b \in Z^+ \cup \{0\}\}$ (or $Z_n$ or $Q^+ \cup \{0\}$, or $R^+ \cup \{0\}$) be a neutrosophic interval semigroup. Take $W = \{[0, a] \,/\, a \in Z^+ \cup \{0\}$ (or (or $Z_n$ or $Q^+ \cup \{0\}$, or $R^+ \cup \{0\})\} \subseteq S$; W is a interval semigroup which is not a neutrosophic interval semigroup. We call W to be a pseudo neutrosophic interval subsemigroup of S.*

We will illustrate this by examples.

***Example 5.36:*** Let $S = \{[0, a + bI] \,|\, a, b \in Z_8\}$ be a neutrosophic interval semigroup. Clearly $W = \{[0, a] \,|\, a \in Z_8\} \subseteq S$ is a pseudo neutrosophic interval subsemigroup of S. Further W is isomorphic to $Z_8$. $P = \{[0, bI] \,|\, b \in Z_8I\} \subseteq S$; is a pure neutrosophic interval subsemigroup of S.

***Example 5.37:*** Let $S = \{[0, a + bI] \,|\, a, b \in Z_5\}$ be a neutrosophic interval semigroup. Clealry $o(S) = 5^2 = 25$. S has a pseudo neutrosophic interval subsemigroup of order 5 and a pure neutrosophic subsemigroup of order 5.

Further consider $[0, (1 + I)^2] = 1 + I + 2I = 1 + 3I = [0, 1 + 3I]$

$$([0, 2+2I])^2 \quad = \quad [0, 4 + 4I + 8I]$$
$$= \quad [0, 4 + 2I]$$

$(0, 1 + I) \ (0, 2+3I) = [0, 2 + 2I + 3I + 3I]$
$\phantom{(0, 1 + I) \ (0, 2+3I)} = [0, 2 + 3I]$

$(0, 2 +3I) \ (0, 3I+2I) = [0, 6 + 9I + 4I + 6I]$
$\phantom{(0, 2 +3I) \ (0, 3I+2I)} = [0, 1 + 4I]$

$(0, 1 + 4I) \ (0, I+4) = [0, I + 4I + 4 + 16I]$
$\phantom{(0, 1 + 4I) \ (0, I+4)} = [0, I + 4]$

$(0, 3 + 3I)^2 = [0, 9 + 9I + 18I]$
$\phantom{(0, 3 + 3I)^2} = [0, 4 + 2I]$

$[0, 2 + 3I]^2 = [0, 4 + 9I + 12I]$
$\phantom{[0, 2 + 3I]^2} = [0, 4 + I]$

and so on.

As in case of pure neutrosophic interval semigroups we can find interval zero divisors, interval idempotents, interval units and interval nilpotents in case of neutrosophic interval semigroups.

**Example 5.38:** Let $S = \{[0, a + bI] \,|\, a, b \in Z_6\}$ be a neutrosophic interval semigroup. S has interval zero divisors given by

$$[0, 3I] \ [0, 2I] = [0, 0] = 0.$$
$$[0, 3] \ [0, 2] = [0, 0] = 0.$$
$$[0, 3I]^2 = [0, 3I] \text{ and}$$
$$[0, 3]^2 = [0, 3]$$

are neutrosophic interval idempotents as well as real interval idempotent in S respectively.
Consider

$[0, 3 + 3I] \ [0, 2 + 2I] = [0, 6 + 6I + 6I + 6I]$
$\phantom{[0, 3 + 3I] \ [0, 2 + 2I]} = [0, 0] = 0.$

Thus S has non trivial interval zero divisors.
$[0, 3I] \ [0, 4] = [0, 0]$ is again a interval zero divisor.

We see [0, 5] [0, 5] = [0, 1] and [0, 5I] [0, 5I] = [0, I].

The reader can study these special elements in case of neutrosophic interval semigroups. Compare (i) the special elements in neutrosophic interval semigroups and those in the pure neutrosophic interval semigroups. (ii) the special substructures like ideals, hyper subsemigroups in these two structures.

**Chapter Six**

# NEUTROSOPHIC INTERVAL MATRIX SEMIGROUPS AND FUZZY INTERVAL SEMIGROUPS

This chapter has three sections. In section one the notion of neutrosophic interval matrix semigroups and pure neutrosophic interval matrix semigroups are introduced. In section two pure neutrosophic interval semigroup polynomials and neutrosophic interval polynomial semigroups are introduced. In the final section the notion of fuzzy interval semigroups are described and studied.

## 6.1 Pure Neutrosophic Interval Matrix Semigroups

In this section we for the first time introduce the notion pure neutrosophic of interval matrix semigroup and neutrosophic interval matrix semigroup and discuss a few of their properties.

We will be using the notations given in chapter five of this book.

**DEFINITION 6.1.1:** *Let $A = \{([0, a_1], [0, a_2], ..., [0, a_n]) \mid a_i \in Z_nI$ (or $Z^+I \cup \{0\}$ or $Q^+I \cup \{0\}$, $R^+I \cup \{0\})\}$ be a row interval neutrosophic matrix. Define usual addition on A so that A becomes a semigroup. A is defined as the pure neutrosophic interval row matrix semigroup under addition. When product (multiplication) is used instead of addition then we call A to be a pure neutrosophic interval row matrix semigroup under multiplication.*

We will illustrate both situations by some examples.

*Example 6.1.1:* Let $A = \{([0, a_1], [0, a_2], ..., [0, a_9]) \mid a_i \in Z^+I \cup \{0\}\}$ be a pure neutrosophic interval row matrix under addition.
If
$$x = ([0, a_1], [0, a_2], ..., [0, a_9])$$
and
$$y = ([0, b_1], [0, b_2], ..., [0, b_9])$$
are in A then

$$x + y = ([0, a_1], [0, a_2], ..., [0, a_9]) + ([0, b_1], [0, b_2], ..., [0, b_9])$$
$$= ([0, a_1 + b_1], [0, a_2 + b_2], ..., [0, a_9 + b_9])$$

is in A so A is a pure neutrosophic interval row matrix under addition. Clearly A has infinite number of elements in it.

*Example 6.1.2:* Let $X = \{([0, a_1], [0, a_2], [0, a_3], [0, a_4], [0, a_5]) \mid a_i \in Z_3 I; 1 \leq i \leq 5\}$ be a pure neutrosophic interval row matrix semigroup of finite order under addition.

It is also a pure neutrosophic interval row matrix semigroup of finite order under multiplication.

*Example 6.1.3:* Let $X = \{([0, a_1], [0, a_2]) \mid a_i \in Q^+I \cup \{0\}\}$ be a pure neutrosophic interval row matrix semigroup of infinite order under addition or multiplication (or used in the mutually exclusive sense).

***Example 6.1.4:*** Let $S = \{([0, a_1], [0, a_2], [0, a_3], [0, a_4], \ldots, [0, a_{12}]) \mid a_i \in R^+I \cup \{0\}\}$ be a pure neutrosophic interval row matrix semigroup under addition (or multiplication). Clearly S is of infinite order we can define subsemigroups and ideals as in case of other pure neutrosophic semigroups.

Here we only describe them by some examples.

***Example 6.1.5:*** Let $T = \{([0, a_1], [0, a_2], \ldots, [0, a_7]) \mid a_i \in Z^+I \cup \{0\}\}$ be a pure neutrosophic interval row matrix semigroup under multiplication. Take $W = \{([0, a_1], [0, a_2], \ldots, [0, a_7]) \mid a_i \in 3Z^+I \cup \{0\}\} \subseteq T$; W is a pure neutrosophic interval row matrix subsemigroup of T. It can also be verified that W is a pure neutrosophic interval row matrix ideal of T.

***Example 6.1.6:*** Let $S = \{([0, a_1], [0, a_2], [0, a_3]) \mid a_i \in Z_7I, 1 \leq i \leq 3\}$ be a pure neutrosophic row matrix interval semigroup under multiplication.
    Consider $T = \{([0, a], [0, a], [0, a]) \mid a \in Z_7I\} \subseteq S$; T is a pure neutrosophic row matrix interval subsemigroup. T is a S-pure neutrosophic row matrix interval subsemigroup as $A = \{([0, a], [0, a], [0, a]) \mid a \in Z_7I \setminus \{0\}\} \subseteq T$ is a pure neutrosophic row matrix interval group under multiplication.

***Example 6.1.7:*** Let $P = \{([0, a_1], [0, a_2], [0, a_3], \ldots, [0, a_8]) \mid a_i \in Z^+I \cup \{0\}; 1 \leq i \leq 8\}$ be a pure neutrosophic row matrix interval semigroup. Take $S = \{([0, a_1], [0, a_2], \ldots, [0, a_8]) \mid a_i \in 5Z^+I \cup \{0\}\} \subseteq P$; S is a pure neutrosophic matrix interval subsemigroup as well as pure neutrosophic matrix interval ideal of P.

Here in this example pure neutrosophic matrix interval subsemigroup of P is an ideal of P. However P is not a S-pure neutrosophic row matrix interval semigroup.

***Example 6.1.8:*** Let $T = \{([0, a_1], [0, a_2], \ldots, [0, a_{12}]) \mid a_i \in Q^+I \cup \{0\}, 1 \leq i \leq 12\} \subseteq T$ be a pure neutrosophic row matrix interval semigroup. Choose $W = \{([0, a_1], [0, a_2], \ldots, [0, a_{12}]) \mid$ $a_i \in Z^+I \cup \{0\}$, $1 \le i \le 12\} \subseteq T$ is a pure neutrosophic row matrix interval subsemigroup of T. Clearly W is not a pure neutrosophic row matrix interval ideal of T. T is a S-pure neutrosophic row matrix interval semigroup for $S = \{([0, a_1], [0, a_2], \ldots, [0, a_{12}]) \mid a_i \in Q^+I, 1 \le i \le 12\} \subseteq T$ is a pure neutrosophic row matrix interval group.

***Example 6.1.9:*** Let $V = \{([0, a_1], [0, a_2], \ldots, [0, a_{18}]) \mid a_i \in R^+I \cup \{0\}, 1 \le i \le 18\}$ be a pure neutrosophic row matrix interval semigroup. V is a S-pure neutrosophic row matrix interval semigroup. V has no pure neutrosophic row matrix interval ideal but has pure neutrosophic row matrix interval subsemigroup.

Now having seen pure neutrosophic row-matrix interval semigroups we now proceed onto define pure neutrosophic column matrix interval semigroups. However it is important to mention that no pure neutrosophic row matrix interval semigroup under addition has non trivial pure neutrosophic row matrix interval ideals, built using $Z_nI$ or $Z^+I \cup \{0\}$ or $R^+I \cup \{0\}$ or $Q^+I \cup \{0\}$.

We can define pure neutrosophic column matrix interval semigroup under addition.

However multiplication cannot be defined on pure neutrosophic column matrix interval semigroup as it is not compatible.

We give examples of them. We cannot define multiplication of pure neutrosophic column matrix interval semigroup.

***Example 6.1.10:*** Let

$$S = \left\{ \begin{bmatrix} [0,a_1] \\ [0,a_2] \\ \vdots \\ [0,a_9] \end{bmatrix} \middle| a_i \in Z_n; 1 \le i \le 9 \right\}$$

be a pure neutrosophic column matrix interval semigroup under addition of finite order.

*Example 6.1.11:* Let

$$Q = \left\{ \begin{bmatrix} [0,a_1] \\ [0,a_2] \\ \vdots \\ [0,a_6] \end{bmatrix} \middle| a_i \in Z^+I \cup \{0\}; 1 \leq i \leq 6 \right\}$$

be a pure neutrosophic column matrix interval semigroup under addition. Q is of infinite order.

*Example 6.1.12:* Let

$$T = \left\{ \begin{bmatrix} [0,a_1] \\ [0,a_2] \\ \vdots \\ [0,a_{18}] \end{bmatrix} \middle| a_i \in Q^+I \cup \{0\}; 1 \leq i \leq 18 \right\}$$

be a pure neutrosophic column matrix interval semigroup under addition.

$$S = \left\{ \begin{bmatrix} [0,a_1] \\ [0,a_2] \\ \vdots \\ [0,a_{18}] \end{bmatrix} \middle| a_i \in Z^+I \cup \{0\} \right\} \subseteq T;$$

S is a pure neutrosophic column matrix interval subsemigroup under addition of T.

*Example 6.1.13:* Let

$$W = \left\{ \begin{bmatrix} [0,a_1] \\ [0,a_2] \\ \vdots \\ [0,a_{25}] \end{bmatrix} \middle| a_i \in R^+I \cup \{0\}; 1 \leq i \leq 25 \right\}$$

be a pure neutrosophic column interval matrix semigroup under addition.

$$T = \left\{ \begin{bmatrix} [0,a_1] \\ [0,a_2] \\ \vdots \\ [0,a_{25}] \end{bmatrix} \middle| a_i \in Q^+I \cup \{0\}; 1 \le i \le 25 \right\} \subseteq W;$$

T is a pure neutrosophic column interval matrix subsemigroup under addition. T is a pure neutrosophic column interval matrix subsemigroup and is not a pure neutrosophic column interval matrix ideal.

*Example 6.1.14:* Let

$$X = \left\{ \begin{bmatrix} [0,a_1] \\ [0,a_2] \\ \vdots \\ [0,a_8] \end{bmatrix} \middle| a_i \in Z_8I; 1 \le i \le 8 \right\}$$

be a pure neutrosophic column interval matrix semigroup under addition. Take

$$W = \left\{ \begin{bmatrix} [0,a_1] \\ [0,a_2] \\ \vdots \\ [0,a_8] \end{bmatrix} \middle| a_i \in \{0, 2I, 4I, 6I\} \subseteq Z_8I; 1 \le i \le 8 \right\} \subseteq X,$$

W is a only a pure neutrosophic column interval matrix subsemigroup under addition of X. Clearly W is not pure neutrosophic column interval matrix ideal of X.

Now we can put forth the question does there exist a S-pure neutrosophic column interval matrix semigroup built using $Z_nI$. We leave this open question.

Now we proceed onto define pure neutrosophic interval matrix semigroup using addition.

**DEFINITION 6.1.2:** *Let P = {Set of all n × m interval matrices with intervals of the form [0,$a_i$] where $a_i \in Z_nI$ (or $Q^+I \cup \{0\}$ or $R^+I \cup \{0\}$ or $Z^+I \cup \{0\}$} (m ≠ n), P under interval matrix addition is a semigroup called the pure neutrosophic n × m interval matrix semigroup.*

Clearly when m ≠ n we cannot define semigroup under multiplication. When m = n we can have pure neutrosophic interval matrix semigroup under addition or under multiplication.

We will illustrate this situation by some examples.

*Example 6.1.15:* Let

$$X = \left\{ \begin{bmatrix} [0,a_1] & [0,a_2] \\ [0,a_3] & [0,a_4] \\ [0,a_5] & [0,a_6] \end{bmatrix} \,\middle|\, a_i \in Z_{18}I; 1 \le i \le 6 \right\}$$

be a pure neutrosophic 3 × 2 interval matrix semigroup under addition. Clearly X is of finite order.

*Example 6.1.16:* Let

$$W = \left\{ \begin{bmatrix} [0,a_1] & [0,a_2] & [0,a_3] \\ [0,a_4] & [0,a_5] & [0,a_6] \end{bmatrix} \,\middle|\, a_i \in Q^+I \cup \{0\}; 1 \le i \le 6 \right\}$$

be a pure neutrosophic 2 × 3 interval matrix semigroup. Clearly W is of infinite order.

*Example 6.1.17 :* Let

$$S = \left\{ \begin{bmatrix} [0,a_1] & [0,a_2] \\ [0,a_3] & [0,a_4] \\ [0,a_5] & [0,a_6] \\ [0,a_7] & [0,a_8] \\ [0,a_9] & [0,a_{10}] \\ [0,a_{11}] & [0,a_{12}] \end{bmatrix} \,\middle|\, a_i \in R^+I \cup \{0\}; 1 \le i \le 12 \right\}$$

be a pure neutrosophic $6 \times 2$ interval matrix semigroup under addition. Clearly S is of infinite order. However we cannot define product semigroup as product is undefined in S.

We can define only one substructure in this case namely pure neutrosophic interval matrix subsemigroup.

We will illustrate this situation by some examples.

*Example 6.1.18:* Let

$$S = \left\{ \begin{bmatrix} [0,a_1] & [0,a_2] \\ [0,a_3] & [0,a_4] \\ [0,a_5] & [0,a_6] \\ [0,a_7] & [0,a_8] \end{bmatrix} \,\middle|\, a_i \in Z_{20}; 1 \le i \le 8 \right\}$$

be a pure neutrosophic $4 \times 2$ interval matrix semigroup.
Take

$$W = \left\{ \begin{bmatrix} [0,a_1] & 0 \\ 0 & [0,a_2] \\ [0,a_3] & 0 \\ 0 & [0,a_4] \end{bmatrix} \,\middle|\, a_i \in Z_{20}; 1 \le i \le 4 \right\} \subseteq S;$$

W is a pure neutrosophic interval matrix subsemigroup of S. However S has no proper pure neutrosophic interval matrix ideals.

*Example 6.1.19:* Let

$$W = \left\{ \begin{bmatrix} [0,a_1] & [0,a_2] \\ [0,a_3] & [0,a_4] \\ [0,a_5] & [0,a_6] \\ [0,a_7] & [0,a_8] \\ [0,a_9] & [0,a_{10}] \end{bmatrix} \middle| a_i \in Q^+I \cup \{0\}; 1 \le i \le 10 \right\}$$

be a pure neutrosophic 5 × 2 interval matrix semigroup under multiplication.
Take

$$S = \left\{ \begin{bmatrix} [0,a_1] & [0,a_2] \\ [0,a_3] & [0,a_4] \\ [0,a_5] & [0,a_6] \\ [0,a_7] & [0,a_8] \\ [0,a_9] & [0,a_{10}] \end{bmatrix} \middle| a_i \in Z^+ \cup \{0\} \right\} \subseteq W.$$

S is only a pure neutrosophic interval matrix subsemigroup of W. We see S is not an ideal.
   Thus several other interesting properties about these structures can be studied.

Now we proceed onto call a pure neutrosophic m × n interval matrix semigroup to be a pure neutrosophic square interval matrix semigroup if m = n. In case of pure neutrosophic interval square matrix semigroup we can define both addition or multiplication. Thus if product is defined we can have the concept of ideals.

We will only illustrate this situation.

***Example 6.1.20:*** Let

$$V = \left\{ \begin{bmatrix} [0,a_1] & [0,a_2] \\ [0,a_3] & [0,a_4] \end{bmatrix} \text{ where } a_1, a_2, a_3, a_4 \in Z_5I \right\}$$

be a pure neutrosophic interval square matrix semigroup under addition. V is of finite order.
   Take

$$U = \left\{ \begin{bmatrix} [0,a] & [0,b] \\ [0,a] & [0,b] \end{bmatrix} \middle| a, b \in Z_5I \right\} \subseteq V,$$

U is a pure neutrosophic interval square matrix subsemigroup of V under addition.
   Clearly V is a S-pure neutrosophic interval square matrix semigroup for

$$P = \left\{ \begin{bmatrix} 0 & 0 \\ 0 & 0 \end{bmatrix}, \begin{bmatrix} [0,I] & [0,I] \\ [0,I] & [0,I] \end{bmatrix}, \begin{bmatrix} [0,2I] & [0,2I] \\ [0,2I] & [0,2I] \end{bmatrix}, \right.$$

$$\left. \begin{bmatrix} [0,3I] & [0,3I] \\ [0,3I] & [0,3I] \end{bmatrix}, \begin{bmatrix} [0,4I] & [0,4I] \\ [0,4I] & [0,4I] \end{bmatrix} \right\}$$

$\subseteq$ V is a pure neutrosophic interval square matrix group under addition. Hence the claim.

   Now we can define product of pure neutrosophic interval square matrices. When a product is defined on these pure neutrosophic square interval matrices we can have ideals in these structures.

***Example 6.1.21:*** Let V = {all 5 × 5 neutrosophic interval matrices with intervals of the form [0, $a_i$]; $a_i \in Z^+I \cup \{0\}$} be a pure neutrosophic interval matrix semigroup under multiplication.

Take W = {All 5 × 5 interval neutrosophic matrices with intervals of the form [0, $a_i$]; $a_i \in 5Z^+ I \cup \{0\}\} \subseteq V$. Clearly W is a pure neutrosophic interval square matrix subsemigroup of V. Infact W is also a pure neutrosophic interval square matrix ideal of V.

Thus every pure neutrosophic interval square matrix subsemigroup built using $Z^+I \cup \{0\}$ is a pure neutrosophic interval square matrix ideal of $Z^+I \cup \{0\}$.

*Example 6.1.22:* Let T = {6 × 6 interval matrices with intervals of the form [0, $a_i$] where $a_i \in R^+I \cup \{0\}\}$ be a pure neutrosophic interval square matrix semigroup under multiplication. T has pure neutrosophic interval square matrix subsemigroups and and pure neutrosophic interval square matrix right ideals or left ideals. Let P = {6 × 6 intervals matrices with intervals of the form [0, $a_i$] where $a_i \in Q^+I \cup \{0\}\} \subseteq T$ is a pure neutrosophic interval square matrix subsemigroup and is clearly not a pure neutrosophic interval square matrix ideal of T.

*Example 6.1.23:* Suppose T given in example 6.1.22 is taken as a pure neutrosophic square matrix interval semigroup under addition. Then we see T has only subsemigroups and no ideals.

Now we proceed onto define the notion of restricted ideals.

**DEFINITION 6.1.3:** *Let S be a pure neutrosophic matrix interval semigroup and $P \subseteq S$ be a pure neutrosophic matrix interval subsemigroup of S. Suppose $T \subseteq P$ be a proper pure neutrosophic interval matrix subsemigroup of P and T is also an ideal of P then we define T to be a pure neutrosophic matrix interval restricted ideal of S restricted to P.*

We will illustrate this situation by some examples.

*Example 6.1.24:* Let S = {all 8 × 8 neutrosophic interval matrices with intervals of the form [0, $a_i$] where $a_i \in R^+I \cup \{0\}\}$ $R^+I \subseteq S$ be a pure neutrosophic interval matrix subsemigroup of S. P = {8 × 8 neutrosophic interval matrices with entries from $Z^+I \cup \{0\}\} \subseteq S$ is a pure neutrosophic interval matrix subsemigroup of S.

Consider T = {All 8 × 8 neutrosophic interval matrices with entries from $5Z^+ \cup \{0\} \subseteq P$; T is a pure neutrosophic interval matrix subsemigroup of P as well as T is a pure neutrosophic interval matrix ideal of S restricted to P.

***Example 6.1.25:*** Let S = {all 5 × 5 neutrosophic interval matrices with intervals of the form m [0,$a_i$] where $a_i \in Q^+I \cup \{0\}$} be a pure neutrosophic interval matrix semigroup. Choose W = {all 5 × 5 neutrosophic interval matrices with intervals of the form [0, $a_i$] where $a_i \in Z^+I \cup \{0\}\} \subseteq S$; W is only a pure neutrosophic interval matrix subsemigroup of S and is not an ideal of S.

Now consider B = {all 5 × 5 neutrosophic interval matrices with intervals of the form [0, $a_i$] where $a_i \in 7Z^+I \cup \{0\}\} \subseteq W$, B is a pure neutrosophic matrix interval subsemigroup of both W and S but B is a pure neutrosophic matrix interval ideal of S relative to the pure neutrosophic interval matrix subsemigroup W of S.

Now having seen such examples, call a pure neutrosophic interval matrix semigroup S to be ideally relatively simple if S has no pure neutrosophic interval matrix ideal relative to any subsemigroup.

We have a large class of ideally relatively pure neutrosophic interval matrix semigroup, which is evident from the following theorem, the proof of which is left as an exercise for the reader.

**THEOREM 6.1.1:** *Let S = {all 1 × n neutrosophic interval matrices with intervals of the form [0, $a_i$]; $a \in Z_pI \setminus \{0\}$ p a prime} be a pure neutrosophic interval matrix semigroup. Here S = {([0, $a_i$], [0, $a_i$], …, [0, $a_i$]) | $a_i \in Z_pI$}.*
*S is a ideally relatively simple pure neutrosophic interval matrix semigroup.*

We will illustrate this situation by an example.

*Example 6.1.26:* Let

$$S = \left\{ \begin{bmatrix} [0,a] & [0,a] \\ [0,a] & [0,a] \end{bmatrix} \middle| a \in Z_{19}I \setminus \{0\} \right\}$$

be a pure neutrosophic interval matrix semigroup. Clearly S is a ideally relatively simple pure neutrosophic matrix interval semigroup.

Also the class of all pure neutrosophic interval matrix semigroups built using $Z_nI$ or $Z^+I \cup \{0\}$ or $R^+I \cup \{0\}$ or $Q^+I \cup \{0\}$ under addition are ideally relatively simple pure neutrosophic interval matrix semigroups.

**COROLLARY 6.1.1:** *Let*

$$A = \left\{ \begin{bmatrix} [0,a] & [0,a] & \cdots & [0,a] \\ [0,a] & [0,a] & \cdots & [0,a] \\ \vdots & \vdots & \cdots & \vdots \\ [0,a] & [0,a] & \cdots & [0,a] \end{bmatrix} \right\}$$

*be n × n interval matrices with a $\in Z_pI$, p a prime} A is an ideally relatively simple pure neutrosophic interval matrix semigroup.*

Thus we see the ideally relatively simple pure neutrosophic matrix interval semigroups depends also on the operation which is defined on the pure neutrosophic matrix interval semigroup.

Now we can have neutrosophic column interval matrix semigroup and neutrosophic interval matrix semigroup.

We will only give examples of them. Properties associated with these structures will also be discussed.

*Example 6.1.27:* Let

$$V = \left\{ \begin{bmatrix} [0,a_1] \\ [0,a_2] \\ [0,a_3] \\ [0,a_4] \end{bmatrix} \middle| a_i \in \{x+yI/x, y \in Z_9\}; 1 \le i \le 4 \right\}$$

be a neutrosophic interval matrix semigroup. Choose

$$W = \left\{ \begin{bmatrix} [0,a_1] \\ [0,a_2] \\ [0,a_3] \\ [0,a_4] \end{bmatrix} \middle| a_i \in \{x+yI\} | x, y \in \{0,3,6\} \subseteq Z_9; 1 \le i \le 4 \right\} \subseteq V$$

W is a neutrosophic column interval matrix subsemigroup of V of finite order.

*Example 6.1.28:* Let

$$V = \left\{ \begin{bmatrix} [0,a_1] \\ [0,a_2] \\ \vdots \\ [0,a_{10}] \end{bmatrix} \middle| a_i \in \{x+yI/x, y \in Z^+ \cup \{0\}\}; 1 \le i \le 10 \right\}$$

be a neutrosophic column interval matrix semigroup of infinite order. Consider

$$P = \left\{ \begin{bmatrix} [0,a_1] \\ [0,a_2] \\ \vdots \\ [0,a_{10}] \end{bmatrix} \middle| a_i \in \{x+yI/x, y \in 3Z^+ \cup \{0\}\}; 1 \le i \le 10 \right\} \subseteq V;$$

P is a neutrosophic column interval matrix subsemigroup of V of infinite order.

We see P is not a neutrosophic column interval matrix ideal of V of infinite order.

*Example 6.1.29:* Let

$$P = \left\{ \begin{bmatrix} [0,a_1] \\ [0,a_2] \\ [0,a_3] \\ [0,a_4] \end{bmatrix} \middle| a_i \in \{x+yI\,/\,x,y \in R^+ \cup \{0\}\}; 1 \leq i \leq 4 \right\}$$

be a neutrosophic column interval matrix semigroup.
  Take

$$W = \left\{ \begin{bmatrix} [0,a_1] \\ [0,a_2] \\ [0,a_3] \\ [0,a_4] \end{bmatrix} \middle| a_i \in \{x+yI\,/\,x,y \in Z^+ \cup \{0\}\} \right\} \subseteq P;$$

W is a neutrosophic column interval matrix subsemigroup of P. Infact W is not an ideal. The neutrosophic column interval matrix semigroups are always only under addition and so it is not possible for them to contain ideals.

*Example 6.1.30:* Let

$$S = \left\{ \begin{bmatrix} [0,a_1] & [0,a_2] \\ [0,a_3] & [0,a_4] \\ [0,a_5] & [0,a_6] \end{bmatrix} \middle| a_i \in \{x+yI\,/\,x,y \in Z_{12}\}; 1 \leq i \leq 6 \right\}$$

be a neutrosophic matrix interval semigroup.
  Take

$$T = \left\{ \begin{bmatrix} [0,a_1] & 0 \\ 0 & [0,a_2] \\ [0,a_3] & 0 \end{bmatrix} \middle| a_i \in \{x+yI/x, y \in Z_{12}\}; 1 \le i \le 3 \right\} \subseteq S;$$

T is a neutrosophic matrix interval subsemigroup of S and is not an ideal of S. However S is of finite order. Take

$$B = \left\{ \begin{bmatrix} [0,a_1] & [0,a_2] \\ [0,a_3] & [0,a_4] \\ [0,a_5] & [0,a_6] \end{bmatrix} \middle| a_i \in Z_{12}I; 1 \le i \le 6 \right\} \subseteq S,$$

B is an interval matrix subsemigroup however B is not a neutrosophic interval matrix ideal of S, thus B is only a pseudo neutrosophic interval matrix subsemigroup of S.

However S has no ideals.

*Example 6.1.31:* Let V =

$$\left\{ \begin{bmatrix} [0,a_1] & [0,a_2] & [0,a_3] & [0,a_4] & [0,a_5] \\ [0,a_6] & [0,a_7] & [0,a_8] & [0,a_9] & [0,a_{10}] \end{bmatrix} \middle| \begin{array}{l} a_i \in \{x+yI \mid x, y \in Z^+ \cup \{0\} \\ 1 \le i \le 10 \end{array} \right\}$$

be a neutrosophic interval matrix semigroup. Take W =

$$\left\{ \begin{bmatrix} [0,a_1] & [0,a_2] & [0,a_3] & [0,a_4] & [0,a_5] \\ [0,a_6] & [0,a_7] & [0,a_8] & [0,a_9] & [0,a_{10}] \end{bmatrix} \middle| a_i \in Z^+I \cup \{0\}; 1 \le i \le 10 \right\}$$

$\subseteq$ V; W is a interval matrix subsemigroup of V. Take

$$B = \left\{ \begin{bmatrix} [0,a_1] & [0,a_2] & [0,a_3] & [0,a_4] & [0,a_5] \\ [0,a_6] & [0,a_7] & [0,a_8] & [0,a_9] & [0,a_{10}] \end{bmatrix} \right.$$

$a_i \in Z^+ \cup \{0\}$; $1 \le i \le 10\} \subseteq$ V; B is a interval matrix subsemigroup of W.

Thus B is a pseudo neutrosophic interval matrix subsemigroup of V.

However V has no neutrosophic interval matrix ideals.

*Example 6.1.32:* Let

$$V = \left\{ \begin{bmatrix} [0,a_1] & [0,a_2] \\ [0,a_3] & [0,a_4] \end{bmatrix} \middle| a_i \in \{x + yI / x, y \in Z_{12}\} \right\}$$

be a neutrosophic interval matrix semigroup under multiplication. Clearly V has pseudo neutrosophic interval matrix subsemigroups, neutrosophic interval matrix subsemigroups, pure neutrosophic interval matrix subsemigroups and neutrosophic matrix interval ideals.

Take

$$P = \left\{ \begin{bmatrix} [0,a_1] & [0,a_2] \\ [0,a_3] & [0,a_4] \end{bmatrix} \middle| a_i \in Z_{12}I; 1 \le i \le 4 \right\} \subseteq V,$$

P is a pure neutrosophic interval matrix subsemigroup of V and is not an ideal.

Consider

$$B = \left\{ \begin{bmatrix} [0,a_1] & [0,a_4] \\ [0,a_2] & [0,a_3] \end{bmatrix} \middle| \begin{array}{c} a_i \in \{x + yI / x, y \in \{0,2,\ldots,10\} \subseteq Z_{12} \\ 1 \le i \le 3 \end{array} \right\}$$

$\subseteq V$; B is a neutrosophic matrix interval subsemigroup of V.

$$C = \left\{ \begin{bmatrix} [0,a_1] & [0,a_2] \\ [0,a_3] & [0,a_4] \end{bmatrix} \middle| a_i \in Z_{12}; 1 \le i \le 4 \right\} \subseteq V;$$

C is a pseudo neutrosophic matrix interval subsemigroup of V. B is also a neutrosophic matrix interval ideal of V.

*Example 6.1.33:* Let V = {3 × 3 neutrosophic interval matrices with intervals of the form [0, a + bI] | a, b ∈ $Z^+ \cup \{0\}$}} be a

neutrosophic matrix interval semigroup under addition. Clearly V has no neutrosophic matrix interval ideals.

*Remark:* If V is a neutrosophic interval matrix semigroup under multiplication certainly V has non trivial neutrosophic interval matrix ideals, V built using $Z^+ \cup \{0\}$.

We give examples of neutrosophic interval matrix ideals of neutrosophic interval matrix semigroup under multiplication.

***Example 6.1.34:*** Let V = {([0,$a_1$], [0,$a_2$], …, [0,$a_n$]) | $a_i$ = {x + yI | x, y ∈ $Q^+ \cup \{0\}$}} be a neutrosophic interval row matrix semigroup. Take W = {([0,$a_1$], 0, …, 0) | $a_1$ = $x_1$ + $y_1$ I; $x_1$, $y_1$ ∈ $Q^+ \cup \{0\}$} ⊆ V; W is a neutrosophic interval row matrix ideal of V.

***Example 6.1.35***: Let V = {([0,$a_1$], [0,$a_2$], [0,$a_3$], [0,$a_4$], [0,$a_5$]) | $a_i$ ∈ N ($Z^+ \cup \{0\}$)} be a neutrosophic interval row matrix semigroup.
$S_1$ = {([0,$a_1$], 0, 0, 0, 0) | $a_1$ ∈ N $Z^+ \cup \{0\}$)} ⊆ V is a neutrosophic interval row matrix ideal of V.
$S_2$ = {(0, [0,$a_2$], 0, 0, 0) | $a_2$ ∈ N $Z^+ \cup \{0\}$)} is a neutrosophic interval row matrix ideal of V.
Thus we have 5 ideals with only one coordinate non zero rest zero.
$H_1$ = {([0,$a_1$], [0,$a_2$], 0, 0, 0) | $a_1$, $a_2$ ∈ N $Z^+ \cup \{0\}$)} is again a neutrosophic interval row matrix ideal of V.
Likewise we have $_5C_2$ = 10 ideals with only two coordinates non zero and three of the coordinates zero.
Let $P_1$ = {([0,$a_1$], [0,$a_2$], [0,$a_3$], 0, 0) | $a_i$ ∈ N $Z^+ \cup \{0\}$)} ⊆ V be a neutrosophic interval row matrix ideal of V. We have $_5C_3$ = 10, neutrosophic ideals with three of the coordinates non zero and two coordinates zero.
Further $T_1$ = {([0,$a_1$], [0,$a_2$], [0,$a_3$], [0,$a_4$], 0) ⊆ V be a neutrosophic interval row matrix ideal of V, clearly we have 5 such ideals which has four coordinates of the row interval matrix to be non zero and one coordinate of the row is zero.
Thus we have atleast $5C_1 + 5C_2 + 5C_3 + 5C_4$ = 5 + 10 + 10 + 5 = 30 nontrivial neutrosophic interval matrix ideals.

In view of this we can have the following theorem which we except the reader to prove.

**THEOREM 6.1.2:** *Let $V = \{([0, a_1], [0, a_2], ..., [0, a_n]) \mid a_i = \{x_i + y_iI\}$ where $x_i, y_i \in Z_t$ (or $Z^+ \cup \{0\}$ or $Q^+ \cup \{0\}$ or $R^+ \cup \{0\}$); $1 \leq i \leq n\}$ be a neutrosophic interval row matrix semigroup. V has atleast $(nC_1 + nC_2 + nC_3 + ... + nC_{n-2} + nC_{n-1})$ number of distinct neutrosophic interval row matrix ideals.*

It is important to note that neutrosophic interval column matrix semigroups have no ideals. Likewise neutrosophic interval s × t (s ≠ t) matrix semigroups have no proper ideals. However if s = t = m the neutrosophic interval m × m matrix semigroup has non trivial proper ideals.

It is pertinent to mention here that neutrosophic matrix interval semigroup in general may or may not satisfy the properties of usual S-semigroups.

We define the determinant of neutrosophic interval n × n square matrix as follows:

$$\text{Let } |A| = \left\| \begin{matrix} [0,a_1] & [0,a_2] \\ [0,a_3] & [0,a_4] \end{matrix} \right\|$$
$$= [0, a_1][0, a_4] - [0, a_2][0, a_3]$$
$$= [0, a_1 a_4] - [0, a_2 a_3]$$
$$= [0, |a_1 a_4 - a_2 a_3|]$$

If $|a_1 a_4 - a_2 a_3| \neq 0$ then we define A to be non singular.

Let
$$A = \begin{bmatrix} [0,a_1] & [0,a_2] & [0,a_3] \\ [0,a_4] & [0,a_5] & [0,a_6] \\ [0,a_7] & [0,a_8] & [0,a_9] \end{bmatrix}$$

be a neutrosophic matrix interval semigroup.

$$|A| = \begin{Vmatrix} [0,a_1] & [0,a_2] & [0,a_3] \\ [0,a_4] & [0,a_5] & [0,a_6] \\ [0,a_7] & [0,a_8] & [0,a_9] \end{Vmatrix}$$

$$= [0, a_1] \begin{Vmatrix} [0,a_5] & [0,a_6] \\ [0,a_8] & [0,a_9] \end{Vmatrix} -$$

$$[0, a_2] \begin{Vmatrix} [0,a_4] & [0,a_6] \\ [0,a_7] & [0,a_9] \end{Vmatrix} +$$

$$[0, a_3] \begin{Vmatrix} [0,a_4] & [0,a_5] \\ [0,a_7] & [0,a_8] \end{Vmatrix}$$

= $[0, a_1] |[[0, a_5] [0, a_9] - [0, a_6] [0, a_8]]| - [0, a_2] |[[0, a_4] [0, a_9] - [0, a_6] [0, a_7] | + [0, a_3] |[0, a_4] [0, a_8] - [0, a_5] [0, a_7]]$

= $[0, a_1] |[0, a_5 a_9] - [0, a_6 a_8] | - [0, a_2] |[0, a_4 a_9] - [0, a_6 a_7]| + [0, a_3] |[0, a_4 a_8] - [0, a_5 a_7]|$

= $[0, | a_1 a_5 a_9 - a_1 a_6 a_8|] - [0, | a_2 a_4 a_9 - a_2 a_6 a_7| ] + [0, |a_3 a_4 a_8 - a_3 a_5 a_7|$

= $[0, |(a_1 a_5 a_9 - a_2 a_6 a_7 + a_3 a_4 a_8) - (a_1 a_6 a_8 + a_2 a_6 a_7 + a_3 a_5 a_7)| ]$

$|A| \neq 0$ if $|(a_1 a_5 a_9 - a_2 a_6 a_7 + a_3 a_4 a_8) - (a_1 a_6 a_8 + a_2 a_6 a_7 + a_3 a_5 a_7| \neq 0$.

Likewise we proceed to find the determinant of order four, order five, etc.

In view of this we have the following interesting theorem.

**THEOREM 6.1.2:** *Let V = {All n × n interval matrices with intervals of the form [0, $a_i$] = [0, $x_i$ + $y_iI$] where $x_i$, $y_i$ ∈ $Z_n$ (or $Z^+$ ∪ {0} or $Q^+$ ∪ {0} or $R^+$ ∪ {0})} be a neutrosophic matrix interval semigroup under product.*

*Let A = {The group generated by all n × n non singular interval matrices in V}* ⊆ *V. V is a S-neutrosophic interval matrix semigroup.*

Proof is direct and is left as an exercise for the reader to prove.

However we can have S-neutrosophic interval m × n matrix semigroup under addition m ≠ n.

*Example 6.1.36:* Let

$$V = \left\{ \begin{bmatrix} [0,a_1] & [0,a_2] \\ [0,a_3] & [0,a_4] \\ [0,a_5] & [0,a_6] \end{bmatrix} \middle| a_i = x_i + y_i I; x_i, y_i \in Z_{10} \right\}; 1 \leq i \leq 6$$

be a neutrosophic interval matrix semigroup under addition. Take

$$P = \left\{ \begin{bmatrix} [0,a_1] & [0,a_2] \\ [0,a_3] & [0,a_4] \\ [0,a_5] & [0,a_6] \end{bmatrix} \middle| \begin{array}{c} a_i = x_i + y_i I; x_i, y_i \in \{0,\ldots,8\} \subseteq Z_{10} \\ 1 \leq i \leq 6 \end{array} \right\}$$

⊆ V; P is a neutrosophic interval matrix group under addition. Thus P is a S-neutrosophic interval matrix semigroup.

*Example 6.1.37:* Let A = {All 10 × 10 neutrosophic interval matrices with intervals of the form [0, $a_i$] with entries from N($Z^+$ ∪ {0})} be a neutrosophic interval matrices semigroup. Take B = {All 10 × 10 neutrosophic interval matrices with intervals of the form [0, $a_i$] with entries from N ($5Z^+$ ∪ {0})} ⊆A, B is an ideal of A.

Infact A has infinite number of neutrosophic matrix interval ideals.

*Example 6.1.38:* Let V = {([0, a], [0, a], [0, a], [0, a], [0, a], [0, a]) | a ∈ $Z_pI$} be a neutrosophic row matrix interval semigroup

under multiplication (p a prime). V is a Smarandache maximal interval semigroup.

***Example 6.1.39:*** Let A = {([0, a], [0, a], [0, a], [0, a]) | a ∈ $Z_{10}I$} be a S-neutrosophic row matrix interval semigroup under multiplication. M = {([0, a], [0, a], [0, a], [0, a]) | a ∈ {2I, 4I, 6I, 8I} ⊆ $Z_{10}I$} be the neutrosophic row matrix interval subgroup of A.

([0, 6I], [0, 6I], [0, 6I], [0, 6I]) is the identity of M under multiplication modulo 10I. Clearly xA = Ax = A for all x ∈ A \ (0, 0, 0, 0), ([0, 5], [0, 5], [0, 5], [0, 5]). For [0, 0, 0, 0] M = {0} and ([0, 5], [0, 5], [0, 5], [0, 5]) M = (0, 0, 0, 0). Thus M is a Smarandache normal neutrosophic matrix interval subgroup of the S-neutrosophic matrix interval semigroup A.

Several other properties enjoyed by semigroups can be studied by any interested reader in case of neutrosophic interval matrix semigroups and pure neutrosophic interval matrix semigroups.

## 6.2 Neutrosophic Interval Polynomial Semigroups

In this section we for the first time introduce notions of pure neutrosophic interval polynomial semigroups and neutrosophic interval polynomial semigroups and discuss a few of the properties related with them. We will be using the notations described in chapter five of this book.

Throughout this book x will denote a variable or an indeterminate.

When we say a polynomial with interval coefficients we mean a polynomial which has coefficients as intervals.

We will first illustrate this before we go for the routine definition.

$p(x) = [0, 5I] x^7 + [0, \sqrt{3}I] x^2 + [0, 7I] x + [0, 8I]$ and
$q(x) = [0, 3/7I] x^8 + [0, I/8] x^5 + [0, 9I] x^3 + [0, 8I] x + [0, 3I]$

are examples of pure neutrosophic interval coefficient polynomials.

$f(x) = [0, \sqrt{3} + I] x^3 + [0, 5I] x^2 + [0, 2/9] x + [0, \sqrt{5} + 7I]$

and

$g(x) = [0, 8I] x^9 + [0, 9I + 8] x^8 + [0, \sqrt{7}] x^7 + [0, 12/7 + 3I] x^5$
$+ [0, 9I] x^2 + [0, 8] x^2 + [0, 8/9 + \sqrt{3} \, I]$

are neutrosophic interval coefficient polynomials.

Now we will proceed onto give the definition.

**DEFINITION 6.2.1:** *Let* $S = \left\{ \sum_{i=0}^{\infty} [0, a_i] x^i \,\middle|\, a_i \in Z_n I \text{ (or } Z^+I \cup \{0\} \text{ } Q^+I \cup \{0\}, R^+I \cup \{0\} \text{ )}, S \text{ under polynomial multiplication is a semigroup; } S \text{ is defined as the pure neutrosophic interval polynomial semigroup under multiplication.}\right.$

We will just illustrate how the product of two pure neutrosophic interval polynomials are done.

Suppose $p(x) = [0, 3I] x^2 + [0, 2I] x + [0, 7I]$ and $g(x) = [0, 12I] x^3 + [0, 3I]x + [0, 2I]$ be two pure neutrosophic interval polynomials in the variable x, then

$p(x) g(x) = ([0, 3I] x^2 + [0, 2I]x + [0, 7I]) ([0, 12I] x^3 + [0, 3I] x + [0, 2I])$

$= [0, 3I] x^2 \cdot [0, 12I] x^3 + [0, 3I] x^2 \cdot [0, 3I] x + [0, 3I] x^2 [0, 2I] + [0, 2I] x [0, 12I] x^3 + [0, 2I] x [0, 3I] x + [0, 2I] x [0, 2I] + [0, 7I] [0, 12I] x^3 + [0, 7I] [0, 3I]x + [0, 7I] [0, 2I]$

$= [0, 36I] x^5 + [0, 9I] x^3 + [0, 6I] x^2 + [0, 24I] x^4 + [0, 6I] x^2 + [0, 4I]x + [0, 84I] x^3 + [0, 21I] x + [0, 14I]$

$= [0, 36I] x^5 + [0, 24I] x^4 + [0, 93I] x^3 + [0, 12I] x^2 + [0, 25I] x + [0, 14I].$

Thus product of two pure neutrosophic interval polynomials are carried out in this manner. It is pertinent to mention here that $I^2$

= I is also used in every multiplication of pure neutrosophic intervals. [10]

We will now give examples of pure neutrosophic interval polynomial semigroups under multiplication.

*Example 6.2.1:* Let

$$S = \left\{ \sum_{i=0}^{\infty} [0, a_i] x^i \;\middle|\; a_i \in Z_n I \right\}$$

be a pure neutrosophic interval polynomial semigroup under multiplication. Clearly S is of infinite order.

*Example 6.2.2:* Let

$$P = \left\{ \sum_{i=0}^{\infty} [0, a_i] x^i \;\middle|\; a_i \in Z^+ I \cup \{0\} \right\}$$

be a neutrosophic interval polynomial semigroup under multiplication of infinite order.

*Example 6.2.3:* Let

$$T = \left\{ \sum_{i=0}^{\infty} [0, a_i] x^i \;\middle|\; a_i \in Q^+ I \cup \{0\} \right\}$$

be a pure neutrosophic interval polynomial semigroup under multiplication of infinite order.

*Example 6.2.4:* Let

$$W = \left\{ \sum_{i=0}^{\infty} [0, a_i] x^i \;\middle|\; a_i \in R^+ I \cup \{0\} \right\}$$

be a pure neutrosophic interval polynomial semigroup under multiplication of infinite order.

Now the natural question which arises in can we have pure neutrosophic polynomial interval semigroup of finite order. The answer of yes we give examples of them.

*Example 6.2.5:* Let

$$P = \left\{ \sum_{i=0}^{3} [0, a_i] x^i \,\middle|\, a_i \in Z_{12}I; \right.$$

$x^4 = 1$, $x^5 = x$, $x^6 = x^2$, $x^7 = x^3$, $x^8 = 1$ and so on} be a pure neutrosophic interval polynomial semigroup under multiplication. P is of finite order.

Thus we have infinite number of pure neutrosophic interval polynomial semigroups under multiplication. However all of them are built using only $Z_nI$; $n < \infty$. All other pure neutrosophic interval polynomial semigroups constructed using $Z^+I \cup \{0\}$ or $Q^+I \cup \{0\}$ or $R^+I \cup \{0\}$ are of infinite order.

*Example 6.2.6:* Let $P = \left\{ \sum_{i=0}^{8} [0, a_i] x^i \,\middle|\, x^9 = 1, x^{10}, = x, x^{11} = x^2, \right.$ ..., $x^{16} = x^7$ $x^{17} = x^8$, $x^{18} = 1$ and so on where $a_i \in Q^+I \cup \{0\}\}$ be a pure neutrosophic interval polynomial semigroup under multiplication. Clearly P is also of infinite order.

Now we can construct pure neutrosophic interval polynomial semigroups using the operation of interval polynomial addition.

*Example 6.2.7:* Let

$$S = \left\{ \sum_{i=0}^{7} [0, a_i] x^i \,\middle|\, a_i \in Q^+I \cup \{0\} \right\}$$

be a pure neutrosophic polynomial interval semigroup under addition. S is of infinite order.

*Example 6.2.8:* Let

$$P = \left\{ \sum_{i=0}^{7} [0, a_i] x^i \,\middle|\, a_i \in Z_{27}I \right\}$$

be a pure neutrosophic polynomial interval semigroup under addition of finite order.

*Example 6.2.9:* Let

$$W = \left\{ \sum_{i=0}^{3} [0, a_i] x^i \ \middle|\ a_i \in Z_{31} \right\}$$

be a pure neutrosophic interval polynomial semigroup under addition. W is of finite order.

*Example 6.2.10:* Let

$$B = \left\{ \sum_{i=0}^{2} [0, a_i] x^i \ \middle|\ a_i \in R^+ I \cup \{0\} \right\}$$

be a pure neutrosophic interval polynomial semigroup under addition, B is of infinite order.

Now we will proceed onto study the substructures such as pure neutrosophic interval polynomial subsemigroup and pure neutrosophic interval polynomial ideals.

*Example 6.2.11:* Let

$$P = \left\{ \sum_{i=0}^{\infty} [0, a_i] x^i \ \middle|\ a_i \in Z_{17} \right\}$$

be a pure neutrosophic interval polynomial semigroup under multiplication.
    Consider

$$T = \left\{ \sum_{i=0}^{\infty} [0, a_i] x^{2i} \ \middle|\ a_i \in Z_{17} \right\} \subseteq P,$$

T is a pure neutrosophic interval polynomial subsemigroup of P and is not an ideal of P.

*Example 6.2.12:* Let

$$W = \left\{ \sum_{i=0}^{\infty} [0, a_i] x^i \,\middle|\, a_i \in R^+ I \cup \{0\} \right\}$$

be a pure neutrosophic interval polynomial semigroup under multiplication. Clearly

$$B = \left\{ \sum_{i=0}^{\infty} [0, a_i] x^i \,\middle|\, a_i \in Q^+ I \cup \{0\} \right\} \subseteq W$$

is a pure neutrosophic polynomial interval subsemigroup of W and is not an ideal of W.

From these examples it is clear that we can have pure neutrosophic interval polynomial subsemigroups which are not ideals.

*Example 6.2.13:* Let

$$W = \left\{ \sum_{i=0}^{\infty} [0, a_i] x^i \,\middle|\, a_i \in Z^+ I \cup \{0\} \right\}$$

be a pure neutrosophic interval polynomial semigroup under addition. Consider

$$T = \left\{ \sum_{i=0}^{\infty} [0, a_i] x^{2i} \,\middle|\, a_i \in 2Z^+ I \cup \{0\} \right\} \subseteq W;$$

T is a pure neutrosophic interval polynomial subsemigroup of W under addition. However T is not an ideal of W.

*Example 6.2.14:* Let

$$M = \left\{ \sum_{i=0}^{\infty} [0, a_i] x^i \,\middle|\, a_i \in Z^+ I \cup \{0\} \right\}$$

be a pure neutrosophic interval polynomial semigroup under multiplication.
Take

$$N = \left\{ \sum_{i=0}^{\infty} [0, a_i] x^i \;\Big|\; a_i \in 2Z^+ I \cup \{0\} \right\} \subseteq M;$$

N is a pure neutrosophic interval polynomial subsemigroup as well as pure neutrosophic interval polynomial ideal of S.

*Example 6.2.15:* Let

$$P = \left\{ \sum_{i=0}^{\infty} [0, a_i] x^i \;\Big|\; a_i \in Z_{12} I \right\}$$

be a pure neutrosophic interval polynomial semigroup under multiplication.
Take

$$T = \left\{ \sum_{i=0}^{\infty} [0, a_i] x^i \;\Big|\; a_i \in \{0, 2I, 4I, 6I, 10I, 8I\} \subseteq Z_{12} I \right\} \subseteq P;$$

T is a pure neutrosophic interval polynomial subsemigroup of P as well as pure neutrosophic interval polynomial ideal of P. If the operation multiplication on P is replaced by addition we see certainly P has only pure neutrosophic interval polynomial subsemigroups but has no pure neutrosophic interval polynomial ideals in P.
 Let

$$M = \left\{ \sum_{i=0}^{n} [0, a_i] x^i \;\Big|\; a_i \in Z_m \right\}$$

be a pure neutrosophic polynomial interval semigroup under addition. M has no pure neutrosophic interval polynomial ideals.

*Example 6.2.16:* Let
$$W = \left\{ \sum_{i=0}^{8} [0, a_i] x^i \,\middle|\, a_i \in Z_{40}I \right\}$$
be a pure neutrosophic interval polynomial semigroup. Take
$$P = \left\{ \sum_{i=0}^{8} [0, a_i] x^i \,\middle|\, a_i \in \{0, 10I, 20I, 30I\} \right\} \subseteq W;$$

P is a pure neutrosophic interval polynomial ideal of W only if W is taken as a pure neutrosophic interval polynomial semigroup under multiplication. (Here $x^9 = 1$, $x^{10} = x$ and so on). Clearly W is of finite order. Infact W is a S-pure neutrosophic polynomial interval semigroup.

*Example 6.2.17:* Let
$$S = \left\{ \sum_{i=0}^{\infty} [0, a_i] x^i \,\middle|\, a_i \in Z_2 \right\}$$

be a pure neutrosophic interval polynomial semigroup under multiplication. Clearly S has non trivial pure neutrosophic interval polynomial subsemigroups.

Now having seen examples of substructures we proceed onto study more properties related with them.

*Example 6.2.18:* Let
$$P = \left\{ \sum_{i=0}^{\infty} [0, a_i] x^i \,\middle|\, a_i \in Q^+I \cup \{0\} \right\}$$

be a pure neutrosophic interval polynomial semigroup. The notion of S p-Sylow subsemigroup cannot be studied as P is of infinite order.

Likewise the notion of S-Cauchy interval element cannot be analysed in case of these infinite structures. Thus only when we have pure neutrosophic polynomial interval semigroup of finite order we can think of the above said notions.

In view of this we now proceed onto study these concepts.

***Example 6.2.19:*** Let P = {[0, t] $x_i$ | t ∈ $Z_6$I; 0 ≤ i ≤ 6; $x^7$ = I or 1} be a pure neutrosophic interval polynomial semigroup under multiplication.

P = {0, [0, I], [0, 2I], [0, 3I], [0, 4I], [0, 5I], [0, Ix], [0, 2Ix], [0, 3Ix], [0, 4Ix], [0, 5Ix], [0, $Ix^2$], [0, $2Ix^2$], [0, $3Ix^2$], [0, $4Ix^2$], [0, $5Ix^2$], [0, $Ix^3$], [0, $2Ix^3$], [0, $3Ix^3$], [0, $4Ix^3$], [0, $5Ix^3$], [0, $Ix^4$], [0, $2Ix^4$], [0, $3Ix^4$], [0, $4Ix^4$], [0, $5Ix^4$], [0, $Ix^5$], [0, $2Ix^5$], [0, $3Ix^5$], [0, $4Ix^5$], [0, $5Ix^5$], [0, $Ix^6$], [0, $2Ix^6$], [0, $3Ix^6$], [0, $4Ix^6$], [0, $5Ix^6$]} is of order 36.

Consider W = {0, [0, 2I], [0, 2Ix], [0, $2Ix^2$], [0, $2Ix^3$], [0, $2Ix^4$], [0, $2Ix^5$], [0, $2Ix^6$], [0, 4I], [0, 4Ix], [0, $4Ix^2$], [0, $4Ix^3$], [0, $4Ix^4$], [0, $4Ix^5$], [0, $4Ix^6$]} ⊆ P is a pure neutrosophic polynomial interval subsemigroup of P. W is also a pure neutrosophic interval polynomial ideal of P.

Consider B= {[0, 3I], [0, 3Ix], [0, $3Ix^2$], [0, $3Ix^3$], [0, $3Ix^4$], [0, $3Ix^5$], [0, $3Ix^6$], 0} ⊆ P, is both a pure neutrosophic interval polynomial subsemigroup as well as ideal. However o(B) $\not|$ o(P). If C = B \ {0} ⊆ P is taken. C is only a pure neutrosophic interval polynomial subsemigroup and is not an ideal of P. By change of conventition in notation we accept [0, aIx] = [0, aI]x.

***Example 6.2.20:*** Let V= {0, [0, I], [0, 2I], [0, 3I], [0, 4I], [0, 5I], [0, 6I], [0, 7I], [0, Ix], [0, 2Ix], [0, 3Ix], [0, 4Ix], [0, 5Ix], [0, 6Ix], [0, 7Ix], [0, $Ix^2$], [0, $2Ix^2$], [0, $3Ix^2$], [0, $4Ix^2$], [0, $5Ix^2$], [0, $6Ix^2$], [0, $7Ix^2$] | aI ∈ $Z_8$I; 0 ≤ a ≤ 7 and $x^3$ = I, be a pure neutrosophic interval polynomial semigroup of order 22 under polynomial multiplication.

Take W= {0, [0, 2I], [0, 4I], [0, 6I], [0, 2I] x, [0, 4I] x, [0, 6I] x, [0, 2I] $x^2$, [0, 4I] $x^2$, [0, 6I] $x^2$} ⊆ V is a pure neutrosophic interval polynomial subsemigroup of V Clearly W is also an ideal of V. We see o(W) $\not|$ o(V). Consider the set M= {I, 7I, Ix, 7Ix, $Ix^2$, $7Ix^2$} ⊆ V. M is only a pure neutrosophic interval polynomial subsemigroup. We see M is a S-pure neutrosophic interval polynomial subsemigroup as A = {I, 7I} is a group under multiplication with I as its identity.

**Example 6.2.21:** Let V = {0, [0, I], [0, 2I], [0, 3I], ..., [0, 16I], [0, I]x, [0, 2I]x, [0, 3I]x, ..., [0, 16I]x | $x^2$ = I and 0 ≤ a ≤ 17} be a pure neutrosophic interval polynomial semigroup built using $Z_{17}$ I. Clearly o(V) = 33. V has pure neutrosophic polynomial interval subsemigroups whose order does not divide order of V.

**Example 6.2.22:** Let S = {[0, a]x, [0, a], [0, a]$x^2$ | a ∈ $Z_5$I; $x^3$ = I} = { 0, [0, I] [0, 2I], [0, 3I], [0, 4I], [0, I]x, [0, 2I] x, [0, 3I] x, [0, 4I]x, [0, I]$x^2$, [0, 2I]$x^2$, [0, 3I]$x^2$, [0, 4I]$x^3$} be a pure neutrosophic interval polynomial semigroup. Clearly o(S) = 13. But S has proper pure neutrosophic interval polynomial subsemigroups.

**Example 6.2.23:** Let S = {[0, a] $x^i$ | i = 0, 1, 2, 3 | $x^4$ = I, a ∈ $Z_8$I} = {0, [0, I], [0, 2I], ..., [0, 7I] [0, I] x, [0, 2I] x , ..., [0, 7I] x, [0, I] $x^2$, [0, 2I] $x^2$, ..., [0, 7I] $x^2$, [0, I] $x^3$, [0, 2I] $x^3$, ... , [0, 7I] $x^3$} be a pure neutrosophic interval polynomial semigroup of order 29. Clearly S is a S- pure neutrosophic interval polynomial semigroup.

Consider T = {0, [0, 2I], [0, 4I], [0, 6I], [0, 2I] x, [0, 4I] x, [0, 6I] x, [0, 2I] $x^2$, [0, 4I] $x^2$, [0, 4I] $x^3$, [0, 2I] $x^3$, [0, 6I] $x^2$, [0, 6I] $x^3$} ⊆ S; T is also a pure neutrosophic interval polynomial subsemigroup of S. o(T) = 13, 13 ∤ 29. Further take Y = {0, [0, 4I], [0, 4I]x, [0, 4I] $x^2$, [0, 4I] $x^3$} ⊆ S; Y is also a pure neutrosophic interval polynomial subsemigroup of S and T and o(Y) = 5, 5 ∤ 13 , 5 ∤ 29.

Consider M = {[0, I], [0, 7I], [0, I]x, [0, 7I] x, [0, 7I] $x^2$, [0, I] $x^2$, [0, I] $x^3$, [0, 7I] $x^3$} ⊆ S. M is also a pure neutrosophic interval polynomial subsemigroup of order 8. We see T and Y are also pure neutrosophic polynomial interval ideals but M is not a pure neutrosophic polynomial interval ideal, only a subsemigroup. Consider N = {0, [0, I], ..., [0, 7I], [0, I] $x^2$, ..., [0, 7I] $x^2$} ⊆ S, N is a pure neutrosophic interval polynomial subsemigroup of S. o(N) = 15, 15 ∤ 29. Clearly N is also not a pure neutrosophic interval polynomial ideal of S.

**Example 6.2.24:** Let V = {[0, $a_i$] $x^i$ | $x^2$ = I, i = 0, 1, $a_i$ ∈ $Z_{13}$ I} be a pure neutrosophic interval polynomial semigroup.

V = {0, [0, I], ..., [0, 12I], [0, I] x, [0, 2I] x, ..., [0, 12I] x} and o (V) = 25. Consider T = {0, [0, I], [0, 12I], [0, I]x, [0, 12I] x} $\subseteq$ V is a pure neutrosophic interval polynomial subsemigroup of V. Further o(T) = 5 and 5 | 25. Also T is a S- pure neutrosophic interval polynomial subsemigroup as P = {[0, I], [0, 12 I], [0, I] x, [0, 12I] x} $\subseteq$ T is a group under product. Thus V is a S pure neutrosophic polynomial interval semigroup.

***Example 6.2.25:*** Let M = {[0, $a_i$] $x^i$ | $x^3$ = I, $a_i \in Z_{11}$ I} = { 0, [0, I], [0, 2I], ..., [0, 10 I] , [0, I] x, [0, I] x, [0, 2I] x, ..., [0, 10I] x, [0, I] $x^2$, [0, 2I] $x^2$, ..., [0, 10I] $x^2$} be a pure neutrosophic interval polynomial semigroup. o(M) = 31. M is S - pure neutrosophic interval polynomial semigroup, as A = { [0, I];[ 0, 10 I]} $\subseteq$ M is a subgroup.
Also B = {[0, I], [0, 10I], [0, I] x, [0, 10 I] x , [0, 10 I] $x^2$, [0, 1 I] $x^2$} $\subseteq$ M is a subgroup. Further T = {0, [0, I], [0, 10 I]} $\subseteq$ M is a S-pure neutrosophic interval polynomial subsemigroup of M. Likewise P = { [0, I], [0, 10 I], [0, I] x, [0, 10 I] x, [0, I] $x^2$, [0, 10 I] $x^2$, 0} is a S-pure neutrosophic interval polynomial subsemigroup of M. But both T and P are not ideals of M. Clearly M has no pure neutrosophic interval polynomial ideals.

***Example 6.2.26:*** Let S = {[0, $a_i$] $x^i$ | $x^3$ = I ; $a_i \in Z_{12}$ I} be a pure neutrosophic interval polynomial semigroup of order 34. Clearly S is a S pure neutrosophic interval polynomial semigroup. S has nontrivial pure neutrosophic interval polynomial ideal. For take T = {0, [0, 3I], [0, 6I], [0, 9I], [0, 3I]x, [0, 6I] x, [0, 9 I] $x^2$, [0, 3I] $x^2$, [0, 6I] $x^2$, [0, 9I] x} $\subseteq$ S is a pure neutrosophic interval polynomial subsemigroup as well ideal of S. Take K = {0, [0, 4I] [0, 8 I], [0, 4 I] x, [0, 8 I] x, [0, 4 I] $x^2$, [0, 8I] $x^2$} $\subseteq$ S is a pure neutrosophic interval polynomial subsemigroup as well as ideal of S. W = {0, [0, 6 I], [0, 6 I] x, [0, 6 I] $x^2$} is a pure neutrosophic interval polynomial subsemigroup as well as ideal of S.

We see none of these substructures order divide the order of S.

***Example 6.2.27:*** Let M = {[0, $a_i$] $x^i$ | $a_i \in Z_3$I; $x^4$ = I} be a pure neutrosophic interval polynomial semigroup. M = { 0, [0, I], [0,

2 I], [0, I]x, [0, 2 I]x, [0, I]$x^2$, [0, 2 I] $x^2$, [0, I]$x^3$, [0, 2 I]$x^3$} is a pure neutrosophic interval polynomial semigroup of order 9. P = {0, [0, I], [0, 2 I]} $\subseteq$ M; is a pure neutrosophic interval polynomial subsemigroup of M.

Clearly P is not an ideal of M. Also T = {0, [0, I]x, [0, I], [0, I]$x^2$, [0, I]$x^3$} $\subseteq$ M is again a S – pure neutrosophic interval polynomial subsemigroup of M which is not an ideal of M. M has no proper ideals.

***Example 6.2.28:*** Let P = {[0, $a_i$] $x^i$ | $a_i \in Z_2$ I, $x^5$ = I} = {0, [0, I], [0, I]x, [0, I ]$x^2$, [0, I]$x^3$, [0, I]$x^4$} be a pure neutrosophic interval polynomial semigroup of order six. Suppose we take $x^6$ = I instead of $x^5$ = I we get a pure neutrosophic interval polynomial semigroup M of order seven. If we take $x^4$ = I then we get a pure neutrosophic interval polynomial semigroup W of order five. But this W has nontrivial pure neutrosophic polynomial interval subsemigroup of order three given by {0, [0, I], [0, I] $x^2$}; 3 $\nmid$ 5.

Now one natural question is how many prime order pure neutrosophic interval semigroups of this type we can have.

The answer is infinite, for when even $x^n$ = I and n odd we get either a prime order semigroup or an odd order semigroup.

Similarly when $x^m$ = I with m odd we get an even order semigroup. This is true when we take the basic set on which the pure neutrosophic polynomial interval semigroup is built using $Z_2$I

***Example 6.2.29:*** Let M = {[0, $a_i$] $x^i$ | $x^2$ = I, $a_i \in Z_3$ I} be a pure neutrosophic interval polynomial semigroup of this special type under multiplication. M = {0, [0, I], [0, 2 I], [0, I]x, [0, 2 I]x}, order of M is 5. If we take $x^3$ = I instead of $x^2$ = I in M. we get a pure neutrosophic interval polynomial semigroup say T = {0, [0, I], [0, 2 I], [0, I]x, [0, 2 I]$x^2$, [0, 2 I]x, [0, I]$x^2$} is of order seven.

If we take $x^4$ = I then the pure neutrosophic interval polynomial semigroup K = {0, [0, I], [0, 2 I], [0, I]x, [0, 2 I]x, [0, I]$x^2$, [0, 2 I]$x^2$, [0, I]$x^3$, [0, 2 I]$x^3$} is of order 9. Thus we see incase of pure neutrosophic interval polynomial semigroup built using $Z_3I$ is of order 5, 7, 9, 11, 13 and so on and the order is $(2n + 1)$ ; $n = 2, 3, ..., \infty$.

By this method we get very many special pure neutrosophic interval polynomial semigroups of varying order which is independent of the order of the basic set $(Z_p I)$ on which it is built.

***Example 6.2.30:*** Let K = {$[0, a_i] x^i \mid x^3 = I, a_i \in Z_4 I$ } = {0, [0, I], [0, 2 I], [0, 3 I], [0, I]x, [0, 2 I]x, [0, 3 I]x, [0, I]$x^2$, [0, 2 I]$x^2$, [0, 3 I]$x^2$, [0, I]$x^3$, [0, 2 I]$x^3$, [0, 3 I] $x^3$} be a pure neutrosophic interval polynomial semigroup of order 13. If $x^3 = I$ we see the order of the semigroup would be, 10. If $x^2 = I$ then the order of the semigroup is 7 and if $x^8 = I$ then order of the semigroup would be 25 and so on. Thus using $Z_4$ I we get pure neutrosophic interval polynomial semigroups of order 7, 10, 13, 16, 19, 22, 25, 28 and so on.

Thus if $Z_n$ I is used we will get pure neutrosophic interval polynomial semigroups of order n+(n–1), n+2 (n–1), n+3 (n–1), and so on.

All properties like S-Lagrange semigroup, S-p-Sylow semigroup, S-Cauchy element etc can be studied incase of these special polynomial pure neutrosophic interval semigroups also.

Now we proceed on to define neutrosophic interval polynomial semigroups.

**DEFINITION 6.2.2:** *Let*

$$V = \left\{ \sum_{i=0}^{\infty} [0, a_i] x^i \bigg| a_i = x_i + y_i I \in N(Z_n) \right\}$$

*(or $N(Z^+ \cup \{0\})$ or $N(Q^+ \cup \{0\})$ or $N(R^+ \cup \{0\})$); V under interval polynomial multiplication is a semigroup known as neutrosophic interval polynomial semigroup.*

We will first illustrate this situation by some examples.

*Example 6.2.31:* Let

$$S = \left\{ \sum_{i=0}^{\infty} [0, a_i] x^i \;\middle|\; a_i \in N(Z_5) \right\}$$

be a neutrosophic interval polynomial semigroup under multiplication. Clearly S is of infinite order.

*Example 6.2.32:* Let

$$P = \left\{ \sum_{i=0}^{\infty} [0, a_i] x^{2i} \;\middle|\; a_i \in N(Z^+ \cup \{0\}) \right\}$$

be a neutrosophic interval polynomial semigroup under multiplication. Clearly P is of infinite order.

*Example 6.2.33:* Let

$$W = \left\{ \sum_{i=0}^{\infty} [0, a_i] x^i \;\middle|\; a_i \in N(R^+ \cup \{0\}) \right\}$$

be a neutrosophic interval polynomial semigroup under multiplication. Clearly W is also of infinite order.

*Example 6.2.34:* Let

$$S = \left\{ \sum_{i=0}^{7} [0, a_i] x^i \;\middle|\; a_i \in N(Z_{40}) \right\}$$

be a neutrosophic interval polynomial semigroup under addition of finite order. Clearly S is not closed under multiplication unless some condition on x is put.

*Example 6.2.35:* Let

$$S = \left\{ \sum_{i=0}^{\infty} [0, a_i] x^i \;\middle|\; a_i \in N(Z^+ \cup \{0\}) \right\}$$

be a neutrosophic polynomial interval semigroup under addition of infinite order.

*Example 6.2.36:* Let

$$T = \left\{ \sum_{i=0}^{4} [0, a_i] x^i \,\middle|\, a_i \in N(R^+ \cup \{0\}) \right\}$$

be a neutrosophic polynomial interval semigroup under addition. T is of infinite order.

Now we have seen finite and infinite neutrosophic polynomial interval semigroups we will proceed to illustrate their substructures.

*Example 6.2.37:* Let

$$T = \left\{ \sum_{i=0}^{9} [0, a_i] x^i \,\middle|\, x^{10} = 1 \quad a_i \in N(Z_7) \right\}$$

be a neutrosophic interval polynomial semigroup under multiplication for $x^{11} = x$, $x^{12} = x^2$ and so on.

*Example 6.2.38:* Let

$$M = \left\{ \sum_{i=0}^{3} [0, a_i] x^i \,\middle|\, a_i \in N(Z^+ \cup \{0\}), x^4 = 1 \right\}$$

be a neutrosophic interval polynomial semigroup under multiplication (or addition, or used in the mutually exclusive sense) is of infinite order.

*Example 6.2.39:* Let

$$W = \left\{ \sum_{i=0}^{3} [0, a_i] x^i \,\middle|\, a_i \in N(Z_3), x^4 = 1 \right\}$$

be a neutrosophic polynomial interval semigroup under addition (or under multiplication); W is a finite order semigroup.

Now we proceed to give examples of substructures.

*Example 6.2.40:* Let

$$K = \left\{ \sum_{i=0}^{\infty} [0, a_i] x^i \,\middle|\, a_i \in N(Z^+ \cup \{0\}) \right\}$$

be a neutrosophic polynomial interval semigroup under multiplication. Consider

$$T = \left\{ \sum_{i=0}^{\infty} [0, a_i] x^{2i} \,\middle|\, a_i \in N(Z^+ \cup \{0\}) \right\} \subseteq K;$$

T is a neutrosophic polynomial interval subsemigroup of K.

*Example 6.2.41:* Let

$$P = \left\{ \sum_{i=0}^{\infty} [0, a_i] x^i \,\middle|\, a_i \in N(Z^+ \cup \{0\}) \right\}$$

be a neutrosophic polynomial interval semigroup under multiplication. Consider

$$T = \left\{ \sum_{i=0}^{\infty} [0, a_i] x^i \,\middle|\, a_i \in N(3Z^+ \cup \{0\}) \right\} \subseteq P;$$

T is a neutrosophic polynomial interval subsemigroup of P.

*Example 6.2.42:* Let

$$B = \left\{ \sum_{i=0}^{\infty} [0, a_i] x^i \,\middle|\, a_i \in N(Z_{12}) \right\}$$

be a neutrosophic interval polynomial semigroup of infinite order under multiplication. Consider

$$C = \left\{ \sum_{i=0}^{\infty} [0, a_i] x^i \,\bigg|\, a_i \in \{0, 2I, 4I, 6I, 8I, 10I\} \subseteq (Z_{12}) \right\},$$

C is a neutrosophic interval polynomial subsemigroup of infinite order under multiplication.

*Example 6.2.43:* Let

$$T = \left\{ \sum_{i=0}^{8} [0, a_i] x^i \,\bigg|\, a_i \in N(Z_{11}) \right\}$$

be a neutrosophic polynomial interval semigroup under addition of finite order

$$P = \left\{ \sum_{i=0}^{4} [0, a_i] x^i \,\bigg|\, a_i \in N(Z_{11}) \right\}$$

is a neutrosophic polynomial interval subsemigroup of T of finite order.
  Clearly T has no ideals.

*Example 6.2.44:* Let

$$W = \left\{ \sum_{i=0}^{19} [0, a_i] x^i \,\bigg|\, a_i \in N(Q^+ \cup \{0\}) \right\}$$

be a neutrosophic interval polynomial semigroup under addition of infinite order. Take

$$M = \left\{ \sum_{i=0}^{12} [0, a_i] x^i \,\bigg|\, a_i \in N(Z^+ \cup \{0\}) \right\} \subseteq W,$$

M is a neutrosophic interval polynomial subsemigroup of W. Clearly M is not an ideal of W.

*Example 6.2.45:* Let

$$S = \left\{ \sum_{i=0}^{9} [0, a_i]x^i \,\middle|\, x^{10} = 1; a_i \in N(Z_{11}) \right\}$$

be a neutrosophic interval polynomial semigroup under multiplication. Take

$$T = \left\{ \sum_{i=0}^{9} [0, a_i]x^i \,\middle|\, a_i \in Z_{11}I \right\} \subseteq S,$$

is a neutrosophic interval polynomial subsemigroup of S and T is also a neutrosophic interval polynomial ideal of S. S is of finite order.

*Example 6.2.46:* Let

$$M = \left\{ \sum_{i=0}^{12} [0, a_i]x^i \,\middle|\, a_i \in N(Z_3) \right\}$$

be a neutrosophic interval polynomial semigroup under addition. Clearly M is of finite order M has neutrosophic interval polynomial subsemigroup but has no ideals.

*Example 6.2.47:* Let

$$T = \left\{ \sum_{i=0}^{12} [0, a_i]x^i \,\middle|\, a_i \in N(Z_3) \text{ where } x^{13} = 1 \right\}$$

be a neutrosophic interval polynomial semigroup under multiplication.

$$W = \left\{ \sum_{i=0}^{12}[0,a_i]x^i \;\middle|\; a_i \in Z_3 \text{ with } x^{13}=1 \right\} \subseteq T;$$

W is a interval polynomial subsemigroup but is not neutrosophic. Thus W is a pseudo neutrosophic interval polynomial subsemigroup of T. Further

$$A = \left\{ \sum_{i=0}^{12}[0,a_i]x^i \;\middle|\; a_i \in Z_3\, I;\, x^{13}=1 \right\} \subseteq T$$

is a pure neutrosophic interval polynomial subsemigroup of T. A is an ideal of T but W is not an ideal of T.

*Example 6.2.48:* Let

$$M = \left\{ \sum_{i=0}^{\infty}[0,a_i]x^i \;\middle|\; a_i \in N(Q^+ \cup \{0\}) \right\}$$

be a neutrosophic interval polynomial semigroup under multiplication. Take

$$N = \left\{ \sum_{i=0}^{\infty}[0,a_i]x^i \;\middle|\; a_i \in Q^+I \cup \{0\} \right\} \subseteq M,$$

N is a neutrosophic interval polynomial subsemigroup. Also N is an ideal of M.

We give the following theorem the proof of which is left as an exercise for the reader.

**THEOREM 6.2.1:** *Let* $S = \left\{ \sum_{i=0}^{\infty}[0,a_i]x^i \;\middle|\; a_i \in N(Z_n) \right\}$, *(or N $(Q^+ \cup \{0\})$ or (or N $(R^+ \cup \{0\})$, (N $(Z^+ \cup \{0\})$ be a neutrosophic interval polynomial semigroup under multiplication. S has a nontrivial neutrosophic interval polynomial ideal.*

We give an hint to the proof.
Take

$$W = \left\{ \sum_{i=0}^{\infty} [0, a_i] x^i \,\middle|\, a_i \in Z_n I \right.$$

(or $(Q^+ I \cup \{0\})$ or $(R^+ I \cup \{0\})$, or $(Z^+ I \cup \{0\}))\} \subseteq S$; W is an ideal of S.

**THEOREM 6.2.2:** *Let $S = \left\{ \sum_{i=0}^{\infty} [0, a_i] x^i \,\middle|\, a_i \in N(Z_n) \text{ (or } N(Q^+ \cup \{0\}) \text{ or (or } N(R^+ \cup \{0\}), (N(Z^+ \cup \{0\})\}$ be a neutrosophic interval polynomial semigroup under addition. S has no proper ideal but has infinite number of neutrosophic interval polynomial subsemigroup.*

The proof is left as an exercise to the reader.

*Example 6.2.49:* Let

$$V = \left\{ \sum_{i=0}^{\infty} [0, a_i] x^i \,\middle|\, a_i \in N(R^+ \cup \{0\}) \right\}$$

be a neutrosophic interval polynomial semigroup under multiplication. Take

$$W = \left\{ \sum_{i=0}^{\infty} [0, a_i] x^i \,\middle|\, a_i \in N(R^+ I \cup \{0\}) \right\} \subseteq V;$$

W is a neutrosophic interval polynomial subsemigroup of V as well as ideal of V.

*Example 6.2.50:* Let

$$M = \left\{ \sum_{i=0}^{\infty} [0, a_i] x^i \,\middle|\, a_i \in N(Z_{12}) \right\}$$

be a neutrosophic interval polynomial semigroup under multiplication. Consider

$$C = \left\{ \sum_{i=0}^{\infty} [0, a_i] x^i \;\middle|\; a_i \in Z_{12} I \right\} \subseteq M,$$

C is a neutrosophic interval polynomial ideal of M.

*Example 6.2.51:* Let

$$T = \left\{ \sum_{i=0}^{\infty} [0, a_i] x^i \;\middle|\; a_i \in N(Z^+ \cup \{0\}) \right\}$$

be a neutrosophic interval polynomial semigroup under addition. Choose

$$B = \left\{ \sum_{i=0}^{\infty} [0, a_i] x^i \;\middle|\; a_i \in Z^+ I \cup \{0\} \right\} \subseteq T,$$

B is only a neutrosophic interval polynomial subsemigroup of T. Infact

$$S_n = \left\{ \sum_{i=0}^{\infty} [0, a_i] x^i \;\middle|\; a_i \in N(nZ^+ \cup \{0\}) \subseteq N(Z^+ \cup \{0\}) \right\} \subseteq T$$

for $n = 2, 3, \ldots, \infty$ are neutrosophic interval polynomial subsemigroups of T and none of them are ideals of T.

*Example 6.2.52:* Let

$$T = \left\{ \sum_{i=0}^{\infty} [0, a_i] x^i \;\middle|\; a_i \in N(Z_n) \right\}$$

be a neutrosophic interval polynomial semigroup under addition. T has several neutrosophic interval polynomial subsemigroups but has no proper ideals. Take for example

$$P_m = \left\{ \sum_{i=0}^{m} [0, a_i] x^i \;\middle|\; a_i \in N(Z_n) \right\} \subseteq T;$$

$P_m$ is a neutrosophic interval polynomial subsemigroup of T for m = 2, 3, ... but is not an ideal of T.

Now we proceed on to define special type of neutrosophic interval polynomial semigroup under multiplication.

**DEFINITION 6.2.3:** *Let $S = \{[0, a_i] x^i \mid i = 0, 1, 2, ..., n, x^n = 1, n < \infty, a_i \in N(Z_m)$ (or $N(Q^+ \cup \{0\})$ or $N(Z^+ \cup \{0\})$ or $N(R^+ \cup \{0\}))\}$. S under multiplication is a semigroup known as the special neutrosophic interval polynomial semigroup. If $a_i \in N(Z_m)$, S is of finite order otherwise S is of infinite order.*

We are more interested in the study of these structures over $N(Z_m)$.

We will first illustrate this situation by some examples.

*Example 6.2.53:* Let $S = \{[0, a_i] x^i \mid x^2 = 1, a_i \in N(Z_3)\}$ be a special neutrosophic interval polynomial semigroup of finite order. S = {0, [0, 1]x, [0, I] [0, 2], [0, I], [0, 2I], [0, 1 + I], [0, 2 + I], [0, 2 + 2I], [0, 2 + 2I], [0, 1]x, [0, 2]x, [0, 2I]x, [0, 1 + I]x, [0, 2 + I]x, [0, 1 + 2I]x, [0, 2 I + 2]x} and order of S is 17. We see S is a S- special neutrosophic interval polynomial semigroup. S has neutrosophic interval polynomial subsemigroup T = {0, [0, I], [0, 1], [0, 2], [0, 2I], [0, 1 + I], [0, 1 + 2I], [0, 2 + I], [0, 2 + 2I]} ⊆ S. Also W = {0, [0, I], [0, 2 I], [0, I]x, [0, 2I]x} ⊆ S is a special neutrosophic interval polynomial subsemigroup. Though the order of S is 17 still S has proper special neutrosophic interval polynomial subsemigroups.

*Example 6.2.54:* Let $V = \{[0, a] x^i \mid x^3 = 1, a \in N(Z_2)\}$ be a special neutrosophic interval polynomial semigroup under multiplication. V = {0, [0, 1], [0, 1]x, [0, I]x, [0, I], [0, 1 + I], [0, 1 + I]x, [0, 1]$x^2$ [0, I]$x^2$, [0, I + 1]$x^2$}; and order of V is 10. S = {0, [0, 1], [1, I], [0, 1 + I]} ⊆ V is a neutrosophic interval subsemigroup of V. T = {0, [0, 1]x [0, 1]$x^2$ [0, 1]} $\subseteq$ V is again a special neutrosophic interval polynomial subsemigroup of V.

M = {0, [0, I], [0, I]x, [0, I]$x^2$} $\subseteq$ S is a special pure neutrosophic interval polynomial subsemigroup of V. N = {0, [0, I], [0, 1] [0, I]x, [0, I]$x^2$, [0, 1]x, [0, 1]$x^2$} $\subseteq$ V is also a special neutrosophic interval polynomial subsemigroup of V.

B = {0, [0, 1 + I], [0, 1 + I]x, [0, 1 + I]$x^2$} $\subseteq$ V is also a special pure neutrosophic interval polynomial subsemigroup of V.

C = {0, [0, 1 + I], [0, 1 + I]x, [0, 1 + I]$x^2$, [0, 1]} $\subseteq$ V is also a special neutrosophic interval polynomial subsemigroup of V and o(C) | o(V).

Take A = {0, [0, 1 + I], [0, 1 + I]x, [0, 1 + I]$x^2$, [0, I] } $\subseteq$ V is a special neutrosophic interval polynomial subsemigroup of V and o(A) | o(V). P = {0, [0, 1], [0, I], [0, I]x, [0, I]$x^2$} $\subseteq$ V is a special neutrosophic interval polynomial subsemigroup of V and o(P) | o(V).

*Example 6.2.55:* Let V = {[0, a] $x^i$ | $x^2$ = 1, a $\in$ N($Z_4$)} be a special neutrosophic interval polynomial semigroup under multiplication.

V = {0, [0, 1], [0, 2], [0, 3], [0, 1]x, [0, 2]x, [0, 3]x, [0, I], [0, 2I], [0, 3I], [0, I]x, [0, 2I]x, [0, 3I]x, [0, 1 + I] [0, 1 + 2I], [0, 1 + 3I], [0, 2 + I] [0, 2 + 2I], [0, 2 + 3I], [0, 3 + I] [0, 3 + 2I], [0, 3 + 3I], [0, 1 + I]x [0, 1 + 2I]x [0, 1 + 3I]x, [0, 2 + I]x, [0, 2 + 2I]x, [0, 2 + 3I]x, [0, 3 + I]x [0, 3 + 2I]x, [0, 3 + 3I]x}, V is of order 31.

Clearly o(V) is a prime. But V has several special neutrosophic interval polynomial subsemigroup. $T_1$ = {0, [0, 1], [0, 2], [0, 3], [0, 1]x [0, 2]x, [0, 3]x} $\subseteq$ V is a special neutrosophic polynomial subsemigroup of V, which is not an ideal of V.

$T_2$ = {0, [0, 2], [0, 2]x, [0, 2I], [0, 2I]x} $\subseteq$ V is again a special neutrosophic interval polynomial subsemigroup of V.

$T_3$ = {0, [0, 2], [0, 2]x, [0, 2I], [0, 2I]x, [0, 2 + 2I] [0, 2 + 2I]x} $\subseteq$ V is again a special neutrosophic interval polynomial subsemigroup of V.

Infact $T_3$ is also a special neutrosophic interval polynomial ideal of V.

V is a S special neutrosophic interval polynomial semigroup for X = {[0, 1], [0, 3]} ⊆ V is a group under multiplication.

***Example 6.2.56:*** Let $V = \{[0, a] x^i \mid x^2 = 1, a \in N(Z_5)\}$ be a special neutrosophic interval polynomial semigroup under multiplication. V = {0, [0, 1], [0, 2], [0, 3], [0, 4], [0, 1]x, [0, 2]x, [0, 3]x, [0, 4]x, [0, I], [0, 2 I], [0, 3 I], [0, 4I], [0, I]x, [0, 2I]x, [0, 3I]x, [0, 4 I]x, [0, 1 + I], [0, 1 + 2 I], [0, 1 + 3 I], [0, 1 + 4 I], [0, 1 + I]x, [0, 1 + 2 I]x, [0, 1 + 3 I]x, [0, 1 + 4 I]x, [0, 2 + I], [0, 2 + 2 I], [0, 2 + 3 I], [0, 2 + 4 I], [0, 2 + I]x [0, 2 + 2 I]x, [0, 2 + 3 I]x, [0, 2 + 4 I]x, [0, 3 + I], [0, 3 + 2 I], [0, 3 + 3 I], [0, 3 + 4 I], [0, 3 + I]x, [0, 3 + 2 I]x, [0, 3 + 3 I]x, [0, 3 + 4 I]x, [0, 4 + I], [0, 4 + 2 I], [0, 4 + 3 I], [0, 4 + 4 I], [0, 4 + I]x, [0, 4 + 2 I]x, [0, 4 + 3 I]x, [0, 4 + 4 I]x}, that is order of V is 49. Now we will give some of its substructures.

Take $K_1$ = {0, [0, 1], [0, 2], [0, 3], [0, 4]} ⊆ V is a pseudo special neutrosophic interval subsemigroup of order 5. Clearly $K_1$ is not an ideal of V. Consider $K_2$ = {0, [0, I], [0, 2 I], [0, 3 I], [0, 4 I]} ⊆ V, $K_2$ is a special pure neutrosophic interval subsemigroup of V and is not an ideal of V. $K_3$ = {0, [0, 1], [0, 2], [0, 3], [0, 4] [0, I], [0, 2 I], [0, 3 I], [0, 4 I]} ⊆ V is a neutrosophic interval subsemigroup of V and is not an ideal of V. $K_4$ = {0, [0, 1] [0, 2], [0, 3], [0, 4], [0, 1]x, [0, 2]x, [0, 3]x, [0, 4]x} ⊆ V is a special neutrosophic interval polynomial subsemigroup of V and is not an ideal of V.

$K_5$ = {0, [0, 1 + I], [0, 2 + 2 I], [0, 3 + 3 I], [0, 4 + 4 I], [0, 1 + I]x, [0, 2 + 2 I]x, [0, 3 + 3 I]x, [0, 4 + 4 I]x} ⊆ V is a special neutrosophic interval polynomial subsemigroup of V which is also not an ideal of V.

We have the following interesting theorem the proof of which is left as an exercise for the reader.

**THEOREM 6.2.3:** *Let $V = \{[0, a]x^i \mid x^n = 1, n < \infty, a \in N(Z_m)\}$ be a special neutrosophic polynomial interval semigroup under multiplication.*
*Then*
*(a) V has special pure neutrosophic polynomial interval subsemigroups*

*(b) V has interval subsemigroups*
*(c) V has neutrosophic interval subsemigroups*
*(d) V is a S semigroup.*

***Example 6.2.57:*** Let $V = \{[0, a]x^i \mid 0 \leq i \leq \infty$ and $a \in N(Z^+ \cup \{0\})\}$ be the special neutrosophic interval polynomial semigroup under multiplication of infinite order. Clearly V is not a S-semigroup.

In view of this we have the following theorem.

**THEOREM 6.2.4:** *Let $V = \{[0, a] x^i \mid 0 \leq i \leq \infty, a_i \in N(Z^+ \cup \{0\})$, (or $N(Q^+ \cup \{0\})$ or $N(R^+ \cup \{0\}))\}$ be a special neutrosophic interval polynomial semigroup under multiplication of infinite order. Then*
  *(a) V has neutrosophic interval subsemigroup of infinite order.*
  *(b) V has interval subsemigroup of infinite order.*
  *(c) V has no S –p –Sylow subgroup.*
  *(d) V has no S – Lagrange subgroup.*
  *(e) V has no S- Cauchy elements.*

Now we see in case of interval semigroups or interval polynomial semigroups when we use $N(Z^+ \cup \{0\})$ or $N(Q^+ \cup \{0\})$ or $N (R^+ \cup \{0\})$ are of infinite order. However some of them are S semigroups.
  Several interesting properties related with them can be studied.

## 6.3 Fuzzy Interval Semigroups

Now we briefly describe fuzzy interval semigroups and the operations on them. Through out this book the set of special fuzzy intervals would be [0, a]; $0 \leq a \leq 1$. We take the collection of fuzzy intervals and define operations on them.

**DEFINITION 6.3.1:** *Let S = {[0, a] | 0 ≤ a ≤ 1} be the collection of fuzzy intervals. For any x = [0, b] and y = [0, a] define min {x, y} = min {[0, b], [0, a]} = {[0, min {a, b}]}.*

*Clearly S with min. operation on fuzzy intervals is a semigroup known as the fuzzy interval semigroup. Infact S is a semilattice with min operation.*

*If we replace the min. operation on S by the max operation still we see S with max operation is a semigroup called the fuzzy interval semigroup with max operation.*

*Example 6.3.1:* Let S = {[0, $1/2^n$] | n = 0, 1, 2, …, ∞}; S is a fuzzy interval semigroup under min. operation.

*Example 6.3.2:* Let T = {[0, $(7/10)^n$] | n = 0, 1, 2, …, ∞}; S is a fuzzy interval semigroup under max operation.

*Example 6.3.3:* Let W = {[0, 1/n] | 1 ≤ n ≤ 29}, W is a fuzzy interval semigroup under min. operation.

It is interesting and important to note that the fuzzy interval semigroups in examples 6.3.1 and 6.3.2 are of infinite order where as the fuzzy interval semigroup given in example 6.3.3 is of finite order.

We can define fuzzy interval sub semigroups.

**DEFINITION 6.3.2:** *Let S = {[0, a] | 0 ≤ a ≤ 1} be a fuzzy interval semigroup with min. operation Suppose W = {[0, b] | 0 ≤ b ≤ 1} ⊆ S, and if W itself is a fuzzy interval semigroup with min-operation then we call W to be a fuzzy interval subsemigroup with min operation.*

We will illustrate this situation by some examples.

*Example 6.3.4:* Let S = [0, a] | 0 ≤ a ≤ 1} be a fuzzy interval semigroup under min operation. T = {[0, $(3/7)^n$] | n = 0, 1, 2, …, ∞} ⊆ S under min operation is a fuzzy interval subsemigroup of S.

*Example 6.3.5:* Let
$$W = \left\{ \left[0, \frac{1}{2^n}\right] \mid n = 0,1,2,...,\infty \right\}$$
be a fuzzy interval semigroup under min. operation.
Take
$$M = \left\{ \left[0, \frac{1}{8^n}\right] \mid n = 0,1,2,...,\infty \right\} \subseteq W;$$
M is a fuzzy interval subsemigroup of W.

*Example 6.3.6:* Let
$$T = \left\{ \left[0, \frac{1}{3^n}\right] \mid n = 0,1,2,...,\infty \right\}$$
be a fuzzy interval semigroup under min. operation.
$$N = \left\{ \left[0, \frac{1}{3^n}\right] \mid n = 0,1,2,...,15 \right\} \subseteq T$$
is a fuzzy interval subsemigroup of T.

Now one can define ideals of fuzzy interval semigroups identical with usual semigroup as follows.

If $S = \{[0, a] \mid 0 \leq a \leq 1\}$ is a fuzzy interval semigroup. Let $T \subseteq S$ be a fuzzy interval subsemigroup of S. We say T a fuzzy interval ideal of S if for all $[0, t] \in T$ and for all $[0, s] \in S$ we have min $\{[0, s], [0, t]\} = \{[0, \min \{s, t\}]$ is in T.

We will give examples of fuzzy interval ideal.

*Example 6.3.7:* Let
$$T = \left\{ \left[0, \frac{1}{5^n}\right] \mid n = 0,1,2,..., 40 \right\}$$
be a fuzzy interval semigroup with min. operation.

$$I = \left\{ \left[0, \frac{1}{5^m}\right] \middle| m = 20, 21, \ldots, 40 \right\} \subseteq T,$$

is a fuzzy interval subsemigroup. It is easily verified I is an ideal of T.

*Example 6.3.8:* Let
$$S = \left\{ \left[0, \frac{1}{7^n}\right] \middle| n = 0, 1, 2, \ldots, 60 \right\}$$

be a fuzzy interval semigroup under min. operation. Suppose

$$W = \left\{ \left[0, \frac{1}{7^m}\right] \middle| m = 25, 26, \ldots, 60 \right\} \subseteq S,$$

W is a fuzzy interval semigroup under min. operation. W under min. operation is fuzzy interval ideal of S under min. operation.

W under min. operation is fuzzy interval ideal of S under min. operation.

*Example 6.3.9:* Let
$$S = \left\{ \left[0, \frac{1}{8^m}\right] \middle| m = 0, 1, 2, \ldots, \infty \right\}$$

be a fuzzy interval semigroup under max operation.

$$W = \left\{ \left[0, \frac{1}{8^n}\right] \middle| n = 0, 1, 2, \ldots, 20 \right\} \subseteq S,$$

is a fuzzy interval subsemigroup under max operation. W is a fuzzy interval ideal of S.

*Example 6.3.10:* Let
$$M = \left\{ \left[0, \left(\frac{4}{9}\right)^n\right] \middle| n = 0, 1, 2, \ldots, \infty \right\}$$

be a fuzzy interval semigroup under max operation.

$$T = \left\{ \left[ 0, \left( 4/9 \right)^m \right] \,\middle|\, m = 0, 1, 2, \ldots, 120 \right\} \subseteq M$$

is a fuzzy interval subsemigroup under max operation. T is a fuzzy interval ideal under max operation.

Now having seen fuzzy interval semigroups.

We now proceed onto define the concept of fuzzy interval matrix semigroup and fuzzy interval polynomial semigroup.
We will give examples of these structures.

***Example 6.3.11:*** Let $V = \{([0, a_1], [0, a_2], \ldots, [0, a_n]) \mid 0 \leq a_i \leq 1, 1 \leq i \leq n\}$, V is a fuzzy interval row matrix. V under max. operation (or min. operation) or used in the mutually exclusive sense (V, min) is defined as a fuzzy interval row matrix semigroup.
We will describe how the operation is carried out; suppose $X = ([0, a_1], [0, a_2], \ldots, [0, a_n])$ and $Y = ([0, a_1], [0, a_2], \ldots, [0, a_n])$ are in V then

$$\begin{aligned}
\min \{X, Y\} &= \min \{([0, a_1], [0, a_2], \ldots, [0, a_n])\,([0, b_1], [0, b_2], \ldots, [0, b_n])\} \\
&= (\min \{[0, a_1], [0, b_1]\}, \min. \{[0, a_2], [0, b_2]\}, \ldots, \min \{[0, a_n], [0, b_n]\}) \\
&= ([0, \min \{a_1, b_1\}], [0, \min \{a_2, b_2\}], \ldots, [0, \min \{a_n, b_n\}]).
\end{aligned}$$

***Example 6.3.12:*** Let $W = \{([0, a_1], [0, a_2], \ldots, [0, a_m]) \mid 0 \leq a_i \leq 1, i = 1, 2, \ldots, m\}$ be a fuzzy interval row matrix semigroup; W under the operation max is a semigroup.
We define for any $X = ([0, a_1], [0, a_2], \ldots, [0, a_m])$ and $Y = \{([0, b_1], [0, b_2], \ldots, [0, b_m])\}$ in W, max. $\{X, Y\}$ = max. $\{([0, a_1], [0, a_2], \ldots, [0, a_m]), ([0, b_1], [0, b_2], \ldots, [0, b_m])\}$ = (max. $\{[0, a_1], [0, b_1]\}$, max. $\{[0, a_2], [0, b_2]\}$, …, max. $\{[0, a_m], [0,$

$b_m\}]) = ([0, \max.\{a_1, b_1\}], [0, \max\{a_2, b_2\}], \ldots, [0, \max\{a_m, b_m\}])$ is in W.

We illustrate it by an example if $X = ([0, 1], [0, \frac{1}{2}], [0, 1/3], [0, 1/7]$ and $Y = ([0, 1/3], [0, 1], [0, 2/5], [0, 5/9])$, max $\{X, Y\}$ $= \{[0, 1], [0, 1], [0, 2/5], [0, 5/9])$.

All properties studied in case of interval row matrix semigroup can be derived in case of fuzzy interval row matrix semigroup with appropriate modifications.

*Example 6.3.13:* Let

$$X = \left\{ \begin{bmatrix} [0, a_1] \\ [0, a_2] \\ \vdots \\ [0, a_n] \end{bmatrix} \middle| \begin{matrix} 0 \leq a_i \leq 1 \\ 1 \leq i \leq n \end{matrix} \right\}$$

be a collection of all column interval matrices with intervals constructed using the fuzzy set [0, 1]. We define for any

$$x = \begin{bmatrix} [0, a_1] \\ [0, a_2] \\ \vdots \\ [0, a_n] \end{bmatrix}$$

and

$$y = \begin{bmatrix} [0, b_1] \\ [0, b_2] \\ \vdots \\ [0, b_n] \end{bmatrix}$$

in X define max $\{x, y\}$ as

$$\max \left\{ \begin{bmatrix} [0, a_1] \\ [0, a_2] \\ \vdots \\ [0, a_n] \end{bmatrix}, \begin{bmatrix} [0, b_1] \\ [0, b_2] \\ \vdots \\ [0, b_n] \end{bmatrix} \right\}$$

$$= \begin{bmatrix} \max[0,a_1],[0,b_1] \\ \max[0,a_2],[0,b_2] \\ \vdots \\ \max[0,a_n],[0,b_n] \end{bmatrix} = \begin{bmatrix} [0,\max\{a_1,b_1\}] \\ [0,\max\{a_2,b_2\}] \\ \vdots \\ [0,\max\{a_n,b_n\}] \end{bmatrix};$$

$$\max\{x, y\} \in X.$$

Thus X is a fuzzy interval column matrix under max. operation.

*Example 6.3.14:* Let

$$Y = \left\{ \begin{bmatrix} [0,a_1] \\ [0,a_2] \\ \vdots \\ [0,a_m] \end{bmatrix} \,\middle|\, \begin{array}{c} 0 \le a_i \le 1 \\ i = 1, 2, ..., m \end{array} \right\}$$

be a fuzzy interval column matrix under min. operation.
If

$$x = \begin{bmatrix} [0,a_1] \\ [0,a_2] \\ \vdots \\ [0,a_m] \end{bmatrix} \text{ and } y = \begin{bmatrix} [0,b_1] \\ [0,b_2] \\ \vdots \\ [0,b_m] \end{bmatrix}$$

are in Y the min $\{x, y\}$ =

$$\min \left\{ \begin{bmatrix} [0,a_1] \\ [0,a_2] \\ \vdots \\ [0,a_m] \end{bmatrix}, \begin{bmatrix} [0,b_1] \\ [0,b_2] \\ \vdots \\ [0,b_m] \end{bmatrix} \right\}$$

$$= \begin{bmatrix} \min[0,a_1],[0,b_1] \\ \min[0,a_2],[0,b_2] \\ \vdots \\ \min[0,a_m],[0,b_m] \end{bmatrix} = \begin{bmatrix} [0,\min\{a_1,b_1\}] \\ [0,\min\{a_2,b_2\}] \\ \vdots \\ [0,\min\{a_m,b_m\}] \end{bmatrix}$$

is in Y.

Substructures can be defined for these fuzzy interval column (row) matrix semigroups.

*Example 6.3.15:* Let M = {set of all m × n fuzzy intervals matrices with fuzzy intervals of the form $[0, a_i]$; $0 \leq a_i \leq 1$, m ≠ n}, M under max or min operation is a fuzzy interval matrix semigroup.

Properties related with these fuzzy interval matrix semigroups can be analysed as in case of usual matrices.

If m = n we have the fuzzy interval square matrix semigroup. In this case apart from max or min we can define the operation of max, min. Under all these three operation fuzzy interval square matrices form a semigroup. Finally we can also define special fuzzy interval matrix semigroups and special fuzzy interval semigroups.

**DEFINITION 6.3.3:** *Let (S, .) be a interval semigroup. A map $\mu : S \rightarrow [0, 1]$ is called a special fuzzy interval semigroup if $\mu(x.y) = \min \{\mu(x), \mu(y)\}$ for all $x, y \in S$.*

We will illustrate this situation by some examples.

*Example 6.3.16:* Let $S = \{[0, x] \mid x \in Z^+ \cup \{0\}\}$ be an interval semigroup under interval multiplication. That is $[0, x], [0, y] \in S$ then $[0, x], [0, y] = [0, xy]$.

Define $\eta : S \rightarrow [0, 1]$ by $\eta([0, x]) = \begin{cases} \dfrac{1}{x} \text{ if } x \neq 0 \\ 1 \text{ if } x = 0 \end{cases}$

$\eta$ is a special fuzzy interval semigroup.

*Example 6.3.17:* Let $S = \{[0, x] \mid x \in Z_{45}\}$ be an interval semigroup.

Define $\mu : S \to [0, 1]$ by $\mu([0, x]) = \begin{cases} \dfrac{1}{x} \text{ if } x \neq 0 \\ 1 \text{ if } x = 0 \end{cases}$

$\mu$ is a special fuzzy interval semigroup.

Now for any interval row matrix semigroup $S = \{([0, a_1], [0, a_2], \ldots, [0, a_n]) \mid a_i \in Z_n \text{ or } Z^+ \cup \{0\} \text{ or } Q^+ \cup \{0\} \text{ or } R^+ \cup \{0\}\}$
$\eta_s : S \to [0, 1]$ is such that each interval $[0, a_i]$ in the row interval matrix is mapped on to a fuzzy interval $[0, b_i]$, $0 \leq b_i \leq 1$ so that the resultant row interval matrix is a fuzzy row interval matrix.

Then $\eta_s$ is a fuzzy row interval semigroup of the interval row matrix semigroup S.

Similarly we define the special fuzzy column interval matrix semigroup and special fuzzy m × n interval matrix semigroup.

We will illustrate this situation by some examples.

*Example 6.3.18:* Let

$$S = \left\{ \begin{bmatrix} [0, a_1] & [0, a_2] \\ [0, a_3] & [0, a_4] \end{bmatrix} \right\}$$

where $a_i \in Z^+ \cup \{0\}\}$ be a interval square matrix semigroup.
Define $\eta : S \to [0, 1]$ as follows:

$$\eta \left( \begin{bmatrix} [0, a_1] & [0, a_2] \\ [0, a_3] & [0, a_4] \end{bmatrix} \right) = \begin{bmatrix} \left[0, \dfrac{1}{a_1}\right] & \left[0, \dfrac{1}{a_2}\right] \\ \left[0, \dfrac{1}{a_3}\right] & \left[0, \dfrac{1}{a_4}\right] \end{bmatrix}$$

where $a_i \neq 0$ if $a_i = 0$ we replace the interval by $[0, 1]$, $\eta$ is a special fuzzy interval square matrix semigroup.

*Example 6.3.19:* Let

$$V = \left\{ \begin{bmatrix} [0,a_1] & [0,a_2] \\ [0,a_3] & [0,a_4] \\ [0,a_5] & [0,a_6] \\ [0,a_7] & [0,a_8] \end{bmatrix} \middle| a_i \in Q^+ \cup \{0\}; 1 \leq i \leq 8 \right\}$$

be a interval $4 \times 2$ matrix semigroup. Define $\eta : V \to [0, 1]$ as follows:

(We wish to state that be it any interval matrix the interval $[0, a_i]$ is always mapped on to the some fuzzy interval wherever $[0, a_i]$ is present in which ever matrix in S, S the collection of interval matrices.)

$$\eta \left( \begin{bmatrix} [0,a_1] & [0,a_2] \\ [0,a_3] & [0,a_4] \\ [0,a_5] & [0,a_6] \\ [0,a_7] & [0,a_8] \end{bmatrix} \right) = \left\{ \begin{bmatrix} [0,b_1] & [0,b_2] \\ [0,b_3] & [0,b_4] \\ [0,b_5] & [0,b_6] \\ [0,b_7] & [0,b_8] \end{bmatrix} \right.$$

$0 \leq b_i \leq 1$ if $a_i \neq 0$ and 1 if $a_i = 0$}.

Thus this map can be realized as a map from the set of intervals built using $Z_n$ or $Z^+ \cup \{0\}$ or $Q^+ \cup \{0\}$ or $R^+ \cup \{0\}$ into the set of fuzzy interval built using $[0, 1]$.

Thus the interested reader can study this situation.

Thus we get two types, viz. fuzzy interval semigroups and special fuzzy interval semigroups.

Now the same method is used for these two types of fuzzy polynomial interval semigroups.

Thus if

$$S = \left\{ \sum_{i=0}^{\infty} [0, a_i] x^i \middle| a_i \in [0,1]; 0 \leq a_i \leq 1 \right\}$$

S under min (or max) operation is a fuzzy interval polynomial semigroup.

*Example 6.3.20:* Let $S = \{0, [0, 1]x, [0, .5]x^2, [0, .3]x^3, [0, 1]x^7, [0, 0.7]x^5\}$ is not a interval fuzzy polynomial semigroup.

*Example 6.3.21:* Let $S = \{[0, a_i]x^i \mid 0 \leq a_i \leq 1; x^8 = 1\}$ be a fuzzy interval polynomial semigroup. Clearly S is of infinite order. This semigroup can be formed by min or max or product operation.
 For if
$$\min \{[0, 0.5]x^8, [0, .7]x^3\} = [0, \min \{0.5, 0.7\} \min \{x^8, x^3\}]$$
$$= [0.0.5] x^3.$$

Now we give an illustration of special interval fuzzy semigroup.

*Example 6.3.22:* Let $V = \left\{ \sum_{i=0}^{\infty} [0, a_i]x^i \mid a_i \in Z^+ \right\}$ be a interval polynomial semigroup.
Define $\eta : V \to [0, 1]$ as follows

$$\eta\left( \sum_{i=0}^{\infty} [0, a_i]x^i \right) = \sum_{i=0}^{\infty} \left[0, \frac{1}{a_i}\right] x^i \; ;$$

$\eta(V)$ or $(\eta, V)$ is a special fuzzy interval polynomial semigroup.

Interested reader can study this concept as in case of interval semigroups of all types.

**Chapter Seven**

# APPLICATIONS OF INTERVAL SEMIGROUPS

In this chapter we just indicate the applications of these various interval semigroups.

We can think of the challenging problem of replacing the alphabets in finite automation by intervals. This has advantage for each bit can be visualized as an interval.

Hence a study of computer applications to interval semigroups would be invaluable.

These interval semigroups can also be used in finite element analysis and fuzzy finite element analysis. As in free semigroups we use only justraposed symbols; these symbols can be replaced by intervals.

Interval semigroups can be used in biology as semigroups are used to describe certain aspects in the crossing or organisms in genetics.

Fuzzy interval semigroups would be much useful in sociology.

When min or max operations are used instead of product or addition these interval semigorups can be realized as interval semilattices and this structures can find a lot of applications in combinatoric technique to problems.

Numerical interval semigroups has applications to algebraic error correcting codes.

Since all these infinite interval semigroups built using $Z^+ \cup \{0\}$ or $Q^+ \cup \{0\}$ or $R^+ \cup \{0\}$ are orderable these structures can be used in the algebraic modeling of pattern design. Introduction of representation theory of finite symmetric interval semigroups and special symmetric interval semigroups would give an innovative theory with lot of applications, as both these structure contain in them the symmetric interval group as a substructure, that is they are Smarandache semigroups.

As these new structures are just introduced they will in due course of time find several other valuable applications.

**Chapter Eight**

# SUGGESTED PROBLEMS

In this chapter we suggest 241 number of problems for the reader to solve and some of them are innovative.

1. What is the order of the interval semigroup $V = \{[0, a_i] \mid a_i \in Z_{12}\}$. Find at least 2 interval subsemigroups.

2. Obtain some interesting properties about interval semigroups.

3. Can $V = \{[0, a_i] \mid a_i \in R^+ \cup \{0\}\}$, the interval semigroup have ideals? Justify your claim.

4. Can $V = \{[0, a] \mid a \in Q^+ \cup \{0\}\}$, the interval semigroup have ideals? Justify.

5. Find ideals of $V = \{[0, a] \mid a \in Z^+ \cup \{0\}\}$.

6. Give some interesting properties about row matrix interval semigroup.

7. Find the order of $V = \{([0, a_1], [0, a_2], \ldots, [0, a_{10}]) \mid a_i \in Z_7; 1 \le i \le 10\}$, the row matrix interval semigroup.

8. Let $V = \{([0, a_1], [0, a_2], [0, a_3]) \mid a_i \in Z^+ \cup \{0\}\}$ be the row matrix interval semigroup. Find two ideals of V.

9. Let $V = \{([0, a_1], [0, a_2], [0, a_3], [0, a_4], [0, a_5]) \mid a_i \in Z_{40}; 1 \le i \le 5\}$ be the row matrix interval semigroup. What is the order of V? Does V have row matrix interval subsemigroups which does not divide the order of V? Find atleast three row matrix interval subsemigroups.

10. Let $V = \{([0,a], [0,a], [0,a], [0,a]) \mid a \in Z_{19}\}$ be a interval row matrix semigroup. Can V have interval row matrix subsemigroups?

11. Let $V = \{([0,a_1], [0,a_2], [0,a_3], [0,a_4]) \mid a_i \in Z_{19}\}$ be a interval row matrix semigroup. Find atleast five interval row matrix subsemigroups of V. Can V have interval row matrix ideals?

12. Prove in any interval row matrix semigroup of finite order the order of row matrix interval subsemigroup need not divide the order of V.

13. Let $V = \left\{ \begin{bmatrix} [0, a_1] \\ [0, a_2] \\ \vdots \\ [0, a_9] \end{bmatrix} \middle| a_i \in Z_{12}; 1 \le i \le 9 \right\}$ be a column interval matrix semigroup. Does V have ideals? Can V have interval subsemigroups? Justify.

14. Obtain some interesting properties about column interval semigroup.

15. Can ever a column interval semigroup have ideals? Prove your claim. Let $V = \{$all $3 \times 5$ interval matrices with intervals of the form $[0, a_i]$ where $a_i \in Z^+ \cup \{0\}\}$ be a $3 \times 5$ matrix interval semigroup under addition. Can V have ideals? Find at least five distinct matrix interval subsemigroups of V.

16. Let V = {All 7 × 2 interval matrices with intervals of the form [0, $a_i$] where $a_i \in Z_{16}$} be a 7 × 2 matrix interval semigroup. Can V have proper ideals? Find matrix interval subsemigroups of V.

17. Let $V = \left\{ \begin{bmatrix} [0,a] & [0,a] & [0,a] \\ [0,a] & [0,a] & [0,a] \end{bmatrix} \middle| a_i \in Z_{11} \right\}$ be a matrix interval semigroup under addition. Can V have matrix interval subsemigroups? Can V have non trivial matrix interval ideals? Justify your answer.

18. Obtain some interesting results about m × n (m≠n) interval matrix semigroups built using $Z_p$, p a prime.

19. Let V = {5 × 5 interval matrices with entries from $Z_{90}$} be the interval matrix semigroup under multiplication.
    a. Find 3 ideals of V.
    b. What is the order of V?
    c. Does V have interval matrix subsemigroups whose order does not divide the order of V?

20. Give an example of a interval matrix semigroup which has no matrix interval ideals?

21. Does there exist a matrix interval semigroup which has no interval matrix subsemigroups?

22. Obtain some interesting properties about interval matrix semigroups.

23. Does there exist an interval matrix semigroup of order 20?

24. What is the order of the 2 × 2 matrix interval semigroup built using intervals in $Z_5$?

25. What is the order of the 7 × 1 column interval matrix semigroup built using $Z_{10}$?

26. Find the order of the 6 × 3 interval matrix semigroup built using $Z_{17}$.

27. What is the order of V = 
$\left\{ \begin{bmatrix} [0,a_1] & [0,a_2] & [0,a_3] & [0,a_4] & [0,a_5] \\ [0,a_6] & [0,a_7] & [0,a_8] & [0,a_9] & [0,a_{10}] \\ [0,a_{11}] & [0,a_{12}] & [0,a_{13}] & [0,a_{14}] & [0,a_{15}] \end{bmatrix} \middle| \begin{array}{l} a_i \in Z_2; \\ 1 \le i \le 14 \end{array} \right\}$? Can V have ideals? Justify.

28. Find ideals in V = {3 × 3 interval matrices with entries from $Z^+ \cup \{0\}$}; V a matrix interval semigroup under multiplication.

29. Does a matrix interval subsemigroup always have an interval matrix ideal? Prove your claim.

30. Does there exists a matrix interval semigroup in which every matrix interval subsemigroup in an ideal?

31. Does there exist a matrix interval semigroup in which no matrix interval subsemigroup is an ideal?

32. Give an example of a simple matrix interval semigroup.

33. Give an example of a doubly simple matrix interval semigroup.

34. Give an example of a S-matrix interval semigroup.

35. Give an example of a matrix interval semigroup which is not a S-matrix interval semigroup.

36. Give an example of a S-Lagrange matrix interval semigroup.

37. Is every matrix interval semigroup a S-Lagrange matrix interval semigroup? Justify your claim.

38. Give an example of a S-matrix interval hyper subsemigroup.

39. Does there exists a matrix interval semigroup which has no proper S-matrix interval subsemigroup.

40. Define S-coset in case of matrix interval semigroup and illustrate the situation by some examples.

41. Give an example of S-Cauchy matrix interval semigroup.

42. Is every matrix interval semigroup a S-Cauchy interval matrix semigroup?

43. Let $V = \left\{ \begin{bmatrix} [0,a_1] & [0,a_2] \\ [0,a_3] & [0,a_4] \end{bmatrix} \mid a_i \in Z_{10}; 1 \leq i \leq 4 \right\}$ be a matrix interval semigroup.

    a. Is V a S-Lagrange matrix interval semigroup?
    b. Does V contain a S-p-Sylow subgroup?
    c. Can V have S-Cauchy elements?
    d. What is the order of V?
    e. Can V have S-ideals?
    f. Can V have S interval matrix hyper subsemigroups?
    g. Can V have matrix interval subsemigroups W so that o(W) / o(V) ?

44. Give an example of matrix interval semigroup which has no S-Cauchy elements?

45. Let $V = \left\{ \begin{bmatrix} [0,a_1] & [0,a_2] & [0,a_3] \\ [0,a_4] & [0,a_5] & [0,a_6] \end{bmatrix} \mid a_i \in Z_7; 1 \leq i \leq 6 \right\}$ be a matrix interval semigroup under addition modulo 7.
    a. What is the order of V?
    b. Is V a S-matrix interval semigroup?
    c. Does V have S-Lagrange subgroups?
    d. Can V have S-cauchy elements?
    e. Can V have ideals?
    f. Does V have matrix interval subsemigroups W so that o(W) / o(V)?

46. Does there exist an interval matrix semigroup of order 42?

47. Obtain some interesting properties about polynomial interval semigroups of finite order under addition built using $Z_n$ ($n < \infty$).

48. Find examples of polynomial interval semigroups which are not S-polynomial interval semigroups.

49. Find examples of polynomial interval semigroups which are S-polynomial interval semigroups.

50. Let $S = \left\{ \sum_{i=0}^{4} [0, a_i] x^i \mid a_i \in Z_5 \right\}$ be a polynomial interval semigroup under addition.
    a. Find the order S.
    b. Is S a S-polynomial interval semigroup?
    c. Can S have elements of finite order?
    d. Does S have S-Cauchy element?
    e. Find at least 2 polynomial interval subsemigroup of S.
    f. Can S have polynomial interval ideals?

51. Let $S = \left\{ \sum_{i=0}^{9} [0, a_i] x^i \mid a_i \in Z_7 \right\}$ be a polynomial interval semigroup under addition.
    a. Is S a S-polynomial interval semigroup?
    b. Can S have polynomial interval ideals?
    c. Does S have a polynomial interval subsemigroup which divides the order of S?
    d. Does S have a polynomial interval subsemigroup which does not divide the order of S?

52. Let $G = \left\{ \sum_{i=0}^{\infty} [0, a_i] x^i \mid a_i \in Q^+ I \cup \{0\} \right\}$ be a polynomial interval semigroup.
    a. Is G a S-polynomial interval semigroup?
    b. Find proper polynomial interval subsemigroups of G.

53. Let $S = \left\{ \sum_{i=0}^{\infty} [0, a_i] x^i \mid a_i \in Z^+ \cup \{0\} \right\}$ be a polynomial interval semigroup under addition. Can S have non zero polynomial interval ideals?

54. Let $S = \left\{ \sum_{i=0}^{\infty} [0, a_i] x^i \mid a_i \in Z^+ \cup \{0\} \right\}$ be a interval polynomial semigroup under multiplication. S has non trivial interval polynomial ideals. For take $T = \left\{ \sum_{i=0}^{\infty} [0, a_i] x^i \mid a_i \in 3Z^+ \cup \{0\} \right\} \subseteq S$; T is a interval polynomial ideal of S, prove. Can S have interval polynomial subsemigroups that are not interval polynomial ideals? Justify your claim.

55. Prove $S = \left\{ \sum_{i=0}^{\infty} [0, a_i] x^i \mid a_i \in R^+ \cup \{0\} \right\}$ has no interval polynomial ideals.

56. Can the polynomial interval semigroup $S = \left\{ \sum_{i=0}^{\infty} [0, a_i] x^i \mid a_i \in Z \right\}$ have polynomial interval ideals? Justify your claim.

57. Prove some interesting results about polynomial interval semigroup.

58. Give a class of polynomial interval semigroups which has no polynomial interval ideals but only polynomial interval subsemigroups.

59. Does there exists polynomial interval semigroups in which every polynomial interval subsemigroup is a polynomial interval ideal?

60. Give an example of a polynomial interval semigroup which has S-hyper polynomial interval subsemigroups.

61. Give examples of polynomial interval semigroup which has no S-hyper polynomial interval subsemigroups.

62. Give an example of a polynomial interval semigroup in which every S-polynomial interval subsemigroup is a S-polynomial interval hyper subsemigroup.

63. Does a polynomial interval semigroup S have S-Cauchy element? Justify.

64. Give an example of a polynomial interval semigroup S which has S-Lagrange polynomial interval semigroup.

65. Does there exist a polynomial interval semigroup S which is a p-Sylow interval subgroup?

66. Can $S = \left\{ \sum_{i=0}^{9} [0, a_i] x^i \mid a_i \in Z_{10} \text{ in which } x^{10} = 1 = x^0 \text{ and } x^{11} = x \text{ and so on} \right\}$ be a polynomial interval semigroup? Can S have p-Sylow subgroups? Justify your claim.

67. Let $S = \left\{ \sum_{i=0}^{3} [0, a_i] x^i \mid a_i \in Z_3 ; x^4 = x^0 = 1 \text{ and } x^5 = x \right\}$ be a polynomial interval semigroup.
    a. What is the order of S?
    b. Does S have p-Sylow subgroups?
    c. Is S a S-Langrange polynomial interval semigroup?
    d. Can S have polynomial interval ideals?
    e. Can S have polynomial interval subsemigroups which are not polynomial interval ideals?

68. Obtain some interesting properties about interval polynomial monoids.

69. Obtain some interesting results about the interval symmetric semigroups.

70. Find atleast two S-interval subsemigroups of S(X) where $o(S(X)) = 4^4$.

71. Find at least 5-S-Cauchy elements of S(X) where $o(S(X)) = 6^6$.

72. Is S(X) where $o(S(X) = 7^7$ a S-weakly interval symmetric semigroup?

73. Find the S-interval symmetric hyper subsemigroup of S(X) where X has 8 distinct intervals.

74. Can the concept of S-coset be defined for interval symmetric semigroup S(X)?

75. Obtain some interesting properties about S(X).

76. Can S(X) ever be simple?

77. The S-semigroup S(X) is a S-p-Sylow semigroup prove.

78. Find S-Cauchy elements of S(X).

79. Can the concept of S-double coset by defined for S(X)?

80. Does S(X) have S-normal interval subgroups? Justify your answer.

81. Find some interesting properties associated with S(⟨X⟩).

82. Let X be the interval set containing 19 intervals.

    a. Find the order of S(X).
    b. How many S-p-Sylow subgroups does S(X) have?
    c. Find S-interval symmetric subsemigroups of S(X).
    d. Does S(X) have S-interval symmetric hyper semigroups?
    e. Does S(X) have a interval symmetric subsemigroup which divides the order of S (X)?
    f. Can S(X) have a S-interval symmetric subsemigroup whose order does not divide the order of S(X)?

83. Find some interesting properties relating to S(X) and S(⟨X⟩).

84. Can $S(\langle X \rangle)$ have S-interval special symmetric hyper subsemigroup? Justify your answer.

85. Can $S(\langle X \rangle)$ have S-p-Sylow special symmetric subgroup if the cardinality of the interval set X is a composite number?

86. Let $S(X)$ and $S(\langle X \rangle)$ be interval symmetric semigroup and special interval semigroup related with the interval set $X = \{[a_1, b_1], [a_2, b_2], [a_3, b_3]\}$.
    a. Obtain all results for $S(X)$ and $S(\langle X \rangle)$ and compare them.
    b. Are $S(X)$ and $S(\langle X \rangle)$ S-Lagrange weakly interval symmetric semigroups?
    c. Does $S(X)$ have S-p-Sylow subgroups?

87. Obtain some interesting properties about pure neutrosophic interval semigroups built using $Q^+ \cup \{0\}$.

88. Determine the properties enjoyed by the pure neutrosophic interval semigroup built using $Z_n$, n a composite number.

89. Can ever a pure neutrosophic interval semigroup built on $Z^+ \cup \{0\}$ be a S-pure neutrosophic interval semigroup? Justify your claim.

90. Give an example of a S-pure neutrosophic interval semigroup.

91. What is the order of the neutrosophic interval semigroup $S = \{[0, a + bI] \mid a, b \in Z_2\}$? Can S be a S-neutrosophic interval semigroup?

92. Can $S = \{[0, a + bI] \mid a, b \in R^+ \cup \{0\}\}$ be a S-neutrosophic interval semigroup?

93. Let $S = \{[0, a + bI] \mid a, b \in Q^+ \cup \{0\}\}$ be a neutrosophic interval semigroup. Can S have pseudo neutrosophic interval subsemigroups?

94. Let $S = \{[0, a + bI] \mid a, b \in Z_n\}$ be a neutrosophic interval semigroup.

a. What is the order of S?
  b. Is S a S-Lagrange neutrosophic interval semigroup?
  c. Find some neutrosophic interval ideals of S.
  d. Does S have S-Cauchy elements?
  e. Find some proper S-neutrosophic interval subsemigroups of S.
  f. Is S a S-neutrosophic interval hyper subsemigroup?

95. Let $W = \{[0, a + bI] \mid a, b \in Z_{24}\}$ be a neutrosophic interval semigroup. Find all S-neutrosophic interval subsemigroups of W.

96. Find some interesting properties about neutrosophic interval semigroups built using $Z_p$, p a prime.

97. Let $S = \{[0, a + bI] \mid a, b \in Z_7\}$ be a neutrosophic interval semigroup;
    a. Can S have zero divisors?
    b. Does S have S-Cauchy elements?
    c. Is S a S-Lagrange neutrosophic interval semigroup?
    d. Is S a S-weakly Lagrange neutrosophic interval semigroup?
    e. Does S have S-p-Sylow neutrosophic interval subgroups?
    f. Can S have idempotents?
    g. Does S contain units?
    h. Can S have non trivial nilpotents?

98. Let $S = \{[0, a + bI] \mid a, b \in Z_{10}\}$ be a neutrosophic interval semigroup.
    a. Find all units in S.
    b. Find all zero divisors in S.
    c. Does S contain idempotents?
    d. Is S a S-Lagrange neutrosophic interval semigroup?

99. Let $A = \{[0, a + bI] \mid a, b \in Z_{12}\}$ be a neutrosophic interval semigroup.
    a. Can A have ideals? Does A contain S-ideals?
    b. Obtain some neutrosophic interval subsemigroups.
    c. Is A a S-neutrosophic interval semigroup?

       d. Does A contain S-neutrosophic interval hyper subsemigroups?
       e. Does A contain S-cauchy elements?

100. Characterize the zero divisors in $S = \{[0, a + bI] \mid a, b \in Z_n;$ n a composite number$\}$; a neutrosophic interval semigroup.

101. Can $A = \{[0, a + bI] \mid a, b \in Z_p;$ p a prime$\}$, the neutrosophic interval semigroup have zero divisors?

102. Find all units in A described in problem (101).

103. Can A in problem (101) have neutrosophic interval ideals?

104. Obtain some interesting properties about neutrosophic row matrix interval semigroups under addition.

105. Can pure neutrosophic column matrix interval semigroup S have non trivial pure neutrosophic column matrix interval ideals in S?

106. Let $V = \left\{ \begin{bmatrix} [0, a_1] \\ [0, a_2] \\ \vdots \\ [0, a_6] \end{bmatrix} \mid a_i \in R^+ \cup \{0\}, 1 \leq i \leq 6 \right\}$ be a pure neutrosophic interval matrix semigroup. Can V have non trivial pure neutrosophic interval matrix ideals? Justify your claim.

107. Obtain some important properties enjoyed by neutrosophic $5 \times 3$ interval matrix semigroups built using $N(Z^+ \cup \{0\})$.

108. Can a neutrosophic $3 \times 2$ interval matrix semigroup have proper neutrosophic $3 \times 2$ interval matrix ideals?

109. Give an example of a pure neutrosophic Smarandache normal subgroup of a pure neutrosophic $m \times m$ matrix interval semigroup.

110. Give an example of a pure neutrosophic 3 × 3 matrix interval semigroup which has no S-pure neutrosophic normal subgroups.

111. Define some new properties on pure neutrosophic interval matrix semigroups.

112. Define Smarandache inverse of pure neutrosophic interval matrix semigroups.

113. Give an example of a Smarandache non Lagrange neutrosophic matrix interval semigroup.

114. Give an example of a Smarandache Lagrange neutrosophic matrix interval semigroup.

115. Give an example of a S-weakly Lagrange pure neutrosophic interval matrix semigroup.

116. Give an example of a finite Smarandache Cauchy neutrosophic interval matrix semigroup.

117. What is the order of A = {All 2 × 2 interval matrices with special intervals of the form $[0, a_i]$, $a_i = x_i + y_i I$; $x_i, y_i \in Z_3$}, the neutrosophic interval matrix semigroup?
    a. Is A a S-neutrosophic interval matrix semigroup?
    b. Can A have S-neutrosophic interval matrix ideals?

118. Find the order of S where $S = \left\{ \begin{bmatrix} [0, a_1] \\ [0, a_2] \\ \vdots \\ [0, a_7] \end{bmatrix} \middle| a_i = x_i + y_i I; \text{where } x_i, y_i \in Z_8; 1 \leq i \leq 8 \right\}$ is the neutrosophic interval matrix semigroup.

119. Let $V = \left\{ \begin{bmatrix} [0, a_1] \\ [0, a_2] \\ \vdots \\ [0, a_9] \end{bmatrix} \middle| a_i \in Z_{15}I; 1 \leq i \leq 9 \right\}$ be a pure neutrosophic interval matrix semigroup under addition. What is the order of V? Does the order of every pure neutrosophic interval matrix semigroup divide the order of V?

120. Does a pure neutrosophic matrix interval semigroup of order 31 exist?

121. Let $P = \left\{ \begin{bmatrix} [0, a] \\ [0, a] \\ \vdots \\ [0, a] \end{bmatrix} \middle| a \in Z_{43}I \right\}$ be a pure neutrosophic 15 × 1 column interval matrix semigroup. What is the order of P? Can P have proper pure neutrosophic 15 × 1 column interval matrix subsemigroups?

122. Obtain some interesting properties about those pure neutrosophic interval matrix semigroups which have proper ideals in them.

123. Give a class of pure neutrosophic interval matrix semigroups which have no proper subsemigroups.

124. Give a class of pure neutrosophic interval matrix semigroups which have no proper ideals but has proper subsemigroups.

125. What is the order of the neutrosophic interval matrix semigroup $S = \left\{ \begin{bmatrix} [0, a_1] \\ [0, a_2] \\ \vdots \\ [0, a_6] \end{bmatrix} \middle| \begin{array}{l} a_i = x_i + y_i I; \\ x_i, y_i \in Z_7; \\ 1 \leq i \leq 6 \end{array} \right\}$. Can S have

proper neutrosophic interval matrix subsemigroups? Does the order of the subsemigroups divide the order of S?

126. Find the order of the pure neutrosophic interval matrix semigroup $A = \left\{ \begin{bmatrix} [0,a_1] \\ [0,a_2] \\ \vdots \\ [0,a_9] \end{bmatrix} \middle| \begin{array}{l} a_i = x_i + y_i I; \\ x_i, y_i \in Z_4; \\ 1 \leq i \leq 9 \end{array} \right\}$. Is A S-pure neutrosophic interval matrix semigroup? Justify.

127. Give an example of a S-Langrange pure neutrosophic column interval matrix semigroup.

128. Give an example of a pure neutrosophic column interval matrix semigroup which has S-Cauchy elements.

129. Let $P = \left\{ \begin{bmatrix} [0,a_1] \\ [0,a_2] \\ \vdots \\ [0,a_9] \end{bmatrix} \middle| \begin{array}{l} a_i = x_i + y_i I; \\ x_i, y_i \in Z_{12}; \\ 1 \leq i \leq 9 \end{array} \right\}$ be a neutrosophic column matrix interval semigroup.
 a. What is the order of P?
 b. Is P a S-Sylow semigroup?
 c. Can P be a S-neutrosophic column matrix interval semigroup?
 d. Does P contain non trivial ideals?
 e. Is P a S-Lagrange semigroup? Justify your answers.

130. Let $A = \left\{ \begin{bmatrix} [0,a_1] \\ [0,a_2] \\ \vdots \\ [0,a_7] \end{bmatrix} \middle| a_i \in N(Z_{18}); 1 \leq i \leq 7 \right\}$ be a neutrosophic column interval matrix semigroup.
 a. What is the order of A?
 b. Can A be S-Lagrange semigroup?

c. Prove A cannot have ideals!
   d. Can A have S-Cauchy elements?

131. Enumerate the special properties enjoyed by the neutrosophic column interval matrix semigroups of finite order.

132. Let V = {5 × 2 interval matrices with special intervals of the form [0, $a_i$]; $a_i \in Z_7I$} be a pure neutrosophic interval 5 × 2 matrix semigroup.
   a. What is the order of V?
   b. Can V be a S-semigroup? Justify.
   c. Does V have ideals?

133. Let V = {all 2 × 7 interval matrices with intervals of the form [0, $a_i$] where $a_i \in Z^+I \cup \{0\}$} be a pure neutrosophic 2 × 7 interval matrix semigroup matrix under addition.
   a. Can V be a S-pure neutrosophic matrix interval semigroup? Justify.
   b. Prove V cannot have proper ideals.
   c. Find atleast 3 pure neutrosophic interval matrix subsemigroups.
   d. What is the special property enjoyed by this semigroup regarding special elements?

134. Let $G = \left\{ \begin{bmatrix} [0, a_1] \\ [0, a_2] \\ \vdots \\ [0, a_{12}] \end{bmatrix} \middle| a_i \in R^+I \cup \{0\}; 1 \leq i \leq 12 \right\}$ be a pure neutrosophic interval column matrix semigroup.
   a. Does G have ideals?
   b. Find in G pure neutrosophic interval column matrix subsemigroup.
   c. Can G be S-pure neutrosophic interval column matrix semigroup?

135. Let S = {5 × 5 neutrosophic matrices with intervals of the form [0, $a_i$] | $a_i$ ∈ $Z_2$ I} be a pure neutrosophic interval matrix semigroup under multiplication.
    a. Find the order of S.
    b. Find ideals in S.
    c. Is S a S-Lagrange semigroup?
    d. Prove S has zero divisors.
    e. Does S contain idempotents?
    f. Can S have S-Cauchy elements?
    g. Is S a S-semigroup?
    h. What is the structure enjoyed by the 5 × 5 interval diagonal matrices?

136. Let T = {All 7 × 5 neutrosophic interval matrices with intervals of the form [0, $a_i$] where $a_i$ = $x_i$ + $y_i$ I with $x_i$, $y_i$ ∈ $Z_3$} be a neutrosophic interval matrix semigroup under addition.
    a. Find the order of T.
    b. Is T a S – semigroup?
    c. Find atleast 3 subsemigroups of T.
    d. Prove T has pure neutrosophic interval matrix subsemigroups.

137. Let W = {All 3 × 3 neutrosophic interval matrices with intervals of the form [0, $a_i$]; $a_i$ ∈ $Z_{10}$I} be a pure neutrosophic interval matrix semigroup.
    a. Find the order of W.
    b. Find ideals in W.
    c. Does there exists a pure neutrosophic interval matrix subsemigroup in W which is not an ideal of W?
    d. Is W a S-semigroup?
    e. Can W have zero divisors?
    f. Find idempotents in W.
    g. Does the order of every pure neutrosophic interval matrix subsemigroup divide the order of W?
    h. Is W a S-pure neutrosophic interval matrix semigroup?

      i.    Will the set of all non singular interval matrices of W be a pure neutrosophic interval matrix subgroup?

138. Let V = {all 6 × 6 neutrosophic interval matrices, with intervals of the form [0, $a_i$] where $a_i \in Z^+ I \cup \{0\}$} be a pure neutrosophic interval matrix semigroup under interval matrix multiplication.
    a. Is V a S-semigroup?
    b. Find atleast two ideals in V.
    c. Does V have zero divisors?
    d. Can V have idempotents?
    e. Can V have nilpotents?
    f. Does the set of all non singular interval neutrosophic matrices form a group? Justify!
    g. Does V have pure neutrosophic interval matrix subsemigroups which are not ideals?
    h. Can V be S-normal?

139. Let $P = \left\{ \sum_{i=0}^{4} [0, a_i]x^i \,\middle|\, a_i \in Z_{10} \right\}$; $x^5 = 1$, $x^6 = 1$ and so on} be a pure neutrosophic polynomial interval semigroup. Find the order of P.

140. Let $S = \left\{ \sum_{i=0}^{\infty} [0, a_i]x^3 \,\middle|\, a_i \in Z_{19}; ; x^4 = 1 \right\}$ be a pure neutrosophic polynomial interval semigroup. Find order of S. Can S have pure neutrosophic polynomial interval semigroups? Justify.

141. Let $T = \left\{ \sum_{i=0}^{12} [0, a_i]x^i \,\middle|\, a_i \in Z_{41}; x^{13} = 1 \right\}$ be a pure neutrosophic interval coefficient polynomial semigroup.
    a. Find the order of T.
    b. Does T have a S-Langrange subgroup?
    c. Does T contain an ideal?
    d. Find pure neutrosophic interval polynomial subsemigroups of T.
    e. Does T have S-Cauchy elements?

   f. Is T a S-semigroup?

142. Let V = {All n × n interval neutrosophic matrices with entries from $Z_{45}I$} be a pure neutrosophic interval matrix semigroup.
  a. Can V have S-pure neutrosophic interval matrix subsemigroup?
  b. Can V have S-zero divisors?
  c. Can V have S-nilpotents?
  d. Can V have S-idempotent?

143. Give an example of a neutrosophic interval matrix semigroup in which every neutrosophic interval matrix subsemigroup is a neutrosophic interval matrix ideal.

144. Give an example of a neutrosophic interval matrix semigroup in which there exists no neutrosophic interval matrix ideals.

145. Give an example of a neutrosophic interval matrix semigroup which has no zero divisors.

146. Give an example of a neutrosophic interval matrix semigroup which is a S-Lagrange semigroup.

147. Give an example of a neutrosophic matrix semigroup which is a S-Cauchy semigroup.

148. Does there exists a pure neutrosophic matrix interval semigroup S of finite order such that the order of every pure neutrosophic matrix ideal subsemigroup divide the order of S?

149. Give an example a pure neutrosophic matrix interval semigroup of order 41.

150. Give an example of a pure neutrosophic matrix interval semigroup of order 19.

151. Give an example of a pure neutrosophic interval matrix semigroup which has no proper pure neutrosophic interval matrix subsemigroup.

152. Obtain some interesting properties about neutrosophic interval matrix semigroups under multiplication.

153. Can we have a neutrosophic matrix interval semigroup in which no element is of finite order?

154. Does there exists a neutrosophic interval matrix semigroup in which every element is of infinite order?

155. Can we say every matrix in the pure neutrosophic matrix interval semigroup constructed using $Z^+I \cup \{0\}$ is of infinite order? Justify your answer.

156. Enumerate some interesting properties about the pure neutrosophic interval matrices.

157. Give some interesting applications of neutrosophic interval matrices.

158. Let V be any pure neutrosophic interval polynomial semigroup of order p, p a prime, prove V can in general have non trivial pure neutrosophic polynomial subsemigroups.

159. Construct a pure neutrosophic polynomial interval semigroup of order 25.

160. How many such polynomial interval semigroups of order 25 can be constructed (problem 159)?

161. Let $S = \left\{ \sum_{i=0}^{2}[0, a_i]x^i \mid x^3 = I, a_i \in Z_7I \right\}$ be a pure neutrosophic interval polynomial semigroup.
   a. Find order of S.
   b. Does S have non trivial pure neutrosophic interval polynomial subsemigroups?

      c. Is S a S-pure neutrosophic interval polynomial semigroup?
      d. Can S have ideals?
      e. Does S contain atleast one pure neutrosophic interval polynomial subsemigroup P so that o(P) / o(S)? Justify your answer.

162. Prove using $Z_pI$ (p-a prime) one can construct pure neutrosophic interval polynomial semigroup of non prime order.

163. Prove every pure neutrosophic interval polynomial semigroup constructed using $Z_nI$ (n prime or otherwise) is a S-pure neutrosophic interval polynomial semigroup.

164. Can any pure neutrosophic interval polynomial semigroup constructed using $Z^+ I \cup \{0\}$ or $Q^+ I \cup \{0\}$ or $R^+ I \cup \{0\}$ be S-pure neutrosophic interval polynomial semigroup? Substantiate your claim.

165. Construct a pure neutrosophic polynomial interval semigroup which has S- Lagrange subgroup.

166. Give an example of a pure neutrosophic polynomial interval semigroup which has S-Cauchy elements.

167. Study the problem (166) in case of neutrosophic polynomial interval semigroup.

168. Does there exist a pure neutrosophic interval polynomial semigroup which has S-p-Sylow subgroups?

169. Study the problem (168) in case of neutrosophic interval polynomial semigroup and just interval polynomial semigroup.

170. Give an example of a Smarandandache pure neutrosophic weakly Lagrange semigroup which is not a S-Lagrange semigroup.

171. Give an example of a Smarandache Lagrange pure neutrosophic interval polynomial semigroup.

172. Give an example of a pure neutrosophic interval polynomial semigroup in which no element is a Smarandache Cauchy element of S.

173. Give an example of a pure neutrosophic interval polynomial semigroup S in which every element is a S-Cauchy element of S.

174. Can a S-Lagrange pure neutrosophic interval polynomial semigroup be constructed using $Z_{18}I$?

175. Study the concept of Smarandache cosets in case of pure neutrosophic polynomial interval semigroup.

176. Give an example of a Smarandache coset H in a pure neutrosophic interval polynomial semigroup using $Z_{12}I$.

177. Give an example of a Smarandache weakly Lagrange pure neutrosophic interval polynomial semigroup using $Z_{10}I$.

178. Can the same problem be true if $Z_{10}I$ is replaced by $Z_{11}I$?

179. Study problems (177) and (178) in case $Z_{10}I$ and $Z_{11}I$ are replaced by $N(Z_{10})$ and $N(Z_{11})$ respectively.

180. Obtain some interesting applications of these interval semigroups.

181. Can interval semigroups be used in finite automaton?

182. Study the concept of Smarandache double cosets in S-pure neutrosophic interval polynomial semigroups.

183. Study the concept in problem (182) in case of S pure neutrosophic interval matrix semigroups.

184. Does the double coset relation on S-interval semigroups S built using $Z_n$ in general partition S for all subgroups in S?

185. Study the problem (184) in case of S-symmetric interval semigroups.

186. Does there exist a S-weakly Lagrange interval matrix semigroup?

187. Study problem (186) in case of pure neutrosophic interval matrix semigroup.

188. Find all Smarandache p-Sylow subgroups of the symmetric interval semigroup using $\{[a_1, b_1], [a_2, b_2], [a_3, b_3], [a_4, b_4]\}$.

189. Find all Smarandache p-Sylow subgroups of the special symmetric interval semigroup $\langle\{[a_1, b_1], [a_2, b_2], [a_3, b_3], [a_4, b_4]\}\rangle$.

190. Find S-Cauchy elements of the semigroups in problems (188) and (189).

191. Give an example of pure neutrosophic interval polynomial semigroup of order 51 using $Z_{11}I$ which has no S-Cauchy elements.

192. Can $S([a_1, b_1], [a_2, b_2], [a_3, b_3], [a_4, b_4], [a_5, b_5])$ the symmetric interval semigroup contain a subgroup of order 12?

193. Study problem (192) in case of the special symmetric interval semigroup $S(\langle[a_1, b_1], [a_2, b_2], \ldots, [a_n, b_n]\}\rangle)$.

194. Find all subgroups of the special symmetric interval semigroup $S(\langle[a_1, b_1], [a_2, b_2], \ldots, [a_9, b_9]\}\rangle)$.

195. Does there exists a pure neutrosophic interval polynomial semigroup of order 24 in which no element is a S-Cauchy element?

196. Study the problem (195) in case of the semigroup of order 100.

197. Let $S = \{[0, a_i] x_i / x^2 = I, a_i \in Z_{11}I\}$ be a pure neutrosophic interval polynomial semigroup. Can S have S-Cauchy elements? Justify your answer.

198. Let $S = \{[0, a_i] x_i / x^2 = I, a_i \in Z_{19}I\}$ be a pure neutrosophic interval polynomial semigroup.
    a. Prove S has no S-Cauchy element.
    b. Prove S is a S-pure neutrosophic interval polynomial semigroup.
    c. Prove S has non trivial S-pure neutrosophic interval polynomial subsemigroups.

199. Study problem (198) in which $Z_{19}I$ is replaced by $N(Z_{19}) = \{x + yI / x, y \in Z_{19}\}$.

200. Prove $H = S([a_1, b_1], [a_2, b_2], [a_3, b_3]) \times S([a_1, b_1], [a_2, b_2])$ is a S-interval symmetric semigroup and is a S-weakly Lagrange semigroup.

201. Find all S-Cauchy elements of H given in problem (200).

202. Find all S-p-Sylow subgroups of H given problem (200).

203. Let $M = \{[0, a] x^i / x^3 = 1 \; a \in N(Z_6)\}$ be a special neutrosophic interval polynomial semigroup.
    a. Find the number of elements in M.
    b. Is M a S-special neutrosophic interval polynomial semigroup?
    c. Find atleast 3 special neutrosophic interval polynomial subsemigroups of M.
    d. Does M contain ideals?

204. Let $G = \{[0, a] x^i / x^3 = 1 \; a \in N(Z_7)\}$ be a special neutrosophic interval polynomial semigroup.
    a. Find the order of G.
    b. Is G a S-special neutrosophic interval polynomial semigroup?
    c. Is G a S-Lagrange semigroup?
    d. Does G have S-p-Sylow subsemigroups?
    e. Can G have S-Cauchy elements?

205. Let $W = \{[0, a] x^i / x^3 = 1 \; a \in N(Z_2)\}$ be a special neutrosophic interval polynomial semigroup.

a. Find special neutrosophic interval polynomial subsemigroups which are ideals?
b. What is the order of W?
c. Does W have S-Lagrange semigroups?
d. Can W have S-Cauchy elements?

206. Obtain some interesting properties about the neutrosophic interval polynomial semigroup $S = \left\{ \sum_{i=0}^{\infty} [0, a_i] x^i \mid a_i \in N(Z_{10}) \right\}$ under addition.

    a. What is order of S?
    b. Can S have ideals?
    c. Can S have S-Cauchy elements?
    d. Is S a S-neutrosophic interval polynomial semigroup?

207. Prove a neutrosophic interval polynomial semigroup under addition cannot have ideals.

208. Let $S = \left\{ \sum_{i=0}^{\infty} [0, a_i] x^i \mid a_i \in N(Z^+ I \cup \{0\}) \right\}$ be a neutrosophic interval polynomial semigroup under multiplication. Prove S has non trivial neutrosophic interval polynomial ideals.

209. Let $S = \left\{ \sum_{i=0}^{5} [0, a_i] x^i \mid x^6 = 1; a_i \in N(Z_4) \right\}$ be a neutrosophic interval polynomial semigroup under multiplication.
    a. Find the order of S.
    b. Can S have ideals?
    c. Does S contain S Cauchy elements?
    d. Can S contain zero divisors?
    e. Can S have idempotents,
    f. Can S have S-p-Sylow semigroups?

210. Let $S = \left\{ \sum_{i=0}^{7} [0, a_i] x^i \mid x^8 = 1; a_i \in N(Z_6) \right\}$ be a neutrosophic interval polynomial semigroup.

a. Find the order of S.
  b. Is S a S-neutrosophic interval polynomial semigroup?
  c. Does S contain proper ideal?
  d. Find zero divisors and idempotents in S.

211. Obtain some interesting properties about $V = \left\{ \sum_{i=0}^{m} [0, a_i] x^i \mid a_i \in N(Z_n); m < \infty \text{ and } x^{m+1} = 1 \right\}$, the neutrosophic interval polynomial semigroup.

212. Does there exist a neutrosophic interval polynomial semigroup which is a S-Lagrange semigroup?

213. Does there exist a neutrosophic interval polynomial semigroup which is a S-weakly Lagrange semigroup?

214. Obtain some interesting applications of neutrosophic interval polynomial semigroups.

215. Give an example of a neutrosophic interval polynomial semigroup which has S-p-Sylow subgroups.

216. Give an example of a neutrosophic interval polynomial semigroup which has no S-p-Sylow subgroups.

217. Give an example of a neutrosophic polynomial interval semigroup which has no S-hyper subsemigroup.

218. Does there exists a neutrosophic polynomial interval semigroup which has S-hyper subsemigroup?

219. Does there exist a neutrosophic polynomial interval semigroup which is S-normal?

220. Give an example of a neutrosophic interval polynomial semigroup which has S-Cauchy elements.

221. Give an example of a neutrosophic interval polynomial polynomial semigroup which has no S-Cauchy elements and Cauchy elements.

222. Does there exist a special neutrosophic polynomial interval semigroup in which every element is a S-Cauchy element?

223. Obtain some interesting properties about special neutrosophic interval polynomial semigroup.

224. Can we always prove a special neutrosophic interval polynomial semigroup built using $N(Z_p)$, p a prime is always of a prime order?

225. Does there exist a special neutrosophic interval polynomial semigroup which is a S-Smarandache semigroup?

226. Study Smarandache coset properties in a special neutrosophic interval polynomial semigroup $S = \{[0, a] x^i / x^3 = 1, a \in N(Z_8)\}$.

227. Does there exist a Smarandache Cauchy special neutrosophic interval polynomial semigroup of order 256?

228. Enumerate some interesting properties about fuzzy interval semigroups.

229. Obtain some interesting properties about fuzzy interval matrix semigroups.

230. Analyse all the properties related with fuzzy interval polynomial semigroups.

231. Suppose $S = \{$all $10 \times 10$ interval matrices with special intervals of the form $[0, a_i]$ $a_i \in Q^+ \cup \{0\}\}$ be a interval $10 \times 10$ matrix semigroup.
Define $\eta : S \to [0,1]$ so that $S_\eta$ or $(\eta, S)$ is a fuzzy interval $10 \times 10$ matrix semigroup.

232. Let $T = \{([0, a_1], [0, a_2], \ldots, [0, a_{10}]) / a_i \in Z_{20}; 1 \leq i \leq 10\}$ be a interval row matrix semigroup under multiplication.
Define $\eta : T \to [0,1]$ so that $T_\eta$ or $(T, \eta)$ is a special fuzzy interval row matrix semigroup.

233. Define special fuzzy interval row matrix ideal of T given in problem (232)

234. Let $S = \left\{ \sum_{i=0}^{\infty} [0, a_i] x^i \,\middle|\, a_i \in R^+ \right\}$ be a interval polynomial semigroup. Define $\eta : S \to [0,1]$ so that $S_\eta$ or $(S, \eta)$ is a special fuzzy interval polynomial semigroup.

235. For S given in problem (234) define fuzzy interval polynomial semigroup.

236. Determine conditions for special fuzzy interval polynomial semigroup to be a S-special fuzzy interval polynomial semigroup.

237. Given $S = \left\{ \sum_{i=0}^{\infty} [0, a_i] x^i \,\middle|\, a_i \in Z^+ \right\}$ is a interval polynomial semigroup. Find the special fuzzy interval polynomial semigroup $S_\eta$ or $(S, \eta)$.

238. Give some interesting properties about fuzzy interval matrix semigroups.

239. Enumerate the properties enjoyed by S-fuzzy interval polynomial semigroups.

240. Does there exist a finite fuzzy interval semigroup? Give some examples.

241. Let $V = \left\{ \sum_{i=0}^{7} [0, a_i] x^i \,\middle|\, a_i = \frac{1}{2}, \frac{1}{4}, \frac{1}{8}, \frac{1}{16}, \frac{1}{32}, \frac{1}{64}, \frac{1}{128}, \frac{1}{256} \right\}$ under usual product for powers of x and min for $a_i$'s and $x^8 = 1$, i.e., if $\left[0, \frac{1}{8}\right] x^6 \cdot \left[1, \frac{1}{128}\right] x^4 = \left[0, \frac{1}{128}\right] x^2$. This is the operation used on V: V is a fuzzy interval polynomial semigroup of finite order prove!

# FURTHER READING

1. Ashbacher, C., *Introduction to Neutrosophic Logic*, American Research Press, Rehoboth, (2002). http://www.gallup.unm.edu/~smarandache/IntrodNeutLogic.pdf

2. Kuperman, I.B., *Approximate Linear Algebraic Equations*, The new University Mathematics Series, Van Nostrand Reinhold Company, London (1971).

3. Liu, F., and Smarandache, F., Logic: A Misleading Concept. A Contradiction Study toward Agent's Logic, in *Proceedings of the First International Conference on Neutrosophy, Neutrosophic Logic, Neutrosophic Set, Neutrosophic Probability and Statistics,* Florentin Smarandache editor, Xiquan, Phoenix, ISBN: 1-931233-55-1, 147 p., 2002, *also published in* "Libertas Mathematica", University of Texas at Arlington, 22 (2002) 175-187. http://lanl.arxiv.org/ftp/math/papers/0211/0211465.pdf

4. Siraj, A., Bridges, S.M., and Vaughn, R.B., *Fuzzy cognitive maps for decision support in an intelligent intrusion detection systems*, www.cs.msstate.edu/~bridges/papers/nafips2001.pdf

5. Smarandache, F., (editor), *Proceedings of the First International Conference on Neutrosophy, Neutrosophic Set, Neutrosophic Probability and Statistics*, Univ. of New Mexico – Gallup, 2001.

http://www.gallup.unm.edu/~smarandache/NeutrosophicProceedings.pdf

6. Smarandache, F., *Collected Papers III,* Editura Abaddaba, Oradea, (2000). http://www.gallup.unm.edu/~smarandache/CP3.pdf

7. Smarandache, F., Neutrosophic Logic - Generalization of the Intuitionistic Fuzzy Logic, To be presented at the *Special Session on Intuitionistic Fuzzy Sets and Related Concepts*, of International EUSFLAT Conference, Zittau, Germany, 10-12 September 2003. http://lanl.arxiv.org/ftp/math/papers/0303/0303009.pdf

8. Smarandache, F., *A Unifying Field in Logics: Neutrosophic Logic. Neutrosophy, Neutrosophic Set, Neutrosophic Probability and Statistics*, third edition, Xiquan, Phoenix, (2003).

9. Vasantha Kandasamy, W.B., *Smarandache Semigroups*, American Research Press, Rehoboth, (2002).

10. Vasantha Kandasamy, W.B., and Florentin Smarandache, *Fuzzy Interval Matrices, Neutrosophic Interval Matrices and their application*, Hexis, Phoenix, Arizona, (2006).

11. Zimmermann, H.J., *Fuzzy Set Theory and its Applications*, Kluwer, Boston, (1988).

# INDEX

# ABOUT THE AUTHORS

**Dr.W.B.Vasantha Kandasamy** is an Associate Professor in the Department of Mathematics, Indian Institute of Technology Madras, Chennai. In the past decade she has guided 13 Ph.D. scholars in the different fields of non-associative algebras, algebraic coding theory, transportation theory, fuzzy groups, and applications of fuzzy theory of the problems faced in chemical industries and cement industries. She has to her credit 646 research papers. She has guided over 68 M.Sc. and M.Tech. projects. She has worked in collaboration projects with the Indian Space Research Organization and with the Tamil Nadu State AIDS Control Society. She is presently working on a research project funded by the Board of Research in Nuclear Sciences, Government of India. This is her $52^{nd}$ book.

On India's 60th Independence Day, Dr.Vasantha was conferred the Kalpana Chawla Award for Courage and Daring Enterprise by the State Government of Tamil Nadu in recognition of her sustained fight for social justice in the Indian Institute of Technology (IIT) Madras and for her contribution to mathematics. The award, instituted in the memory of Indian-American astronaut Kalpana Chawla who died aboard Space Shuttle Columbia, carried a cash prize of five lakh rupees (the highest prize-money for any Indian award) and a gold medal.
She can be contacted at vasanthakandasamy@gmail.com
Web Site: http://mat.iitm.ac.in/home/wbv/public_html/

**Dr. Florentin Smarandache** is a Professor of Mathematics at the University of New Mexico in USA. He published over 75 books and 150 articles and notes in mathematics, physics, philosophy, psychology, rebus, literature.

In mathematics his research is in number theory, non-Euclidean geometry, synthetic geometry, algebraic structures, statistics, neutrosophic logic and set (generalizations of fuzzy logic and set respectively), neutrosophic probability (generalization of classical and imprecise probability). Also, small contributions to nuclear and particle physics, information fusion, neutrosophy (a generalization of dialectics), law of sensations and stimuli, etc. He can be contacted at smarand@unm.edu